\documentclass[reqno,12pt]{amsart}
\input texdraw
\catcode`\@=11
\def\sideset#1#2#3{%
  \@mathmeasure\z@\displaystyle{#3}%
  \global\setbox\@ne\vbox to\ht\z@{}\dp\@ne\dp\z@
  \setbox\tw@\box\@ne
  \@mathmeasure4\displaystyle{\copy\tw@#1}%
  \@mathmeasure6\displaystyle{#3{#2}}%
  \dimen@-\wd6 \advance\dimen@\wd4 \advance\dimen@\wd\z@
  \hbox to\dimen@{}\mathop{\kern-\dimen@\box4\box6}%
}
\catcode`\@=12

\newfont{\fourteenpoint}{cmr10 scaled\magstep3}
\newfont{\fourteenit}{cmti10 scaled\magstep2}
\newfont{\fourteensl}{cmsl10 scaled\magstep2}
\newfont{\fourteensmc}{cmcsc10 scaled\magstep2}
\newfont{\fourteentt}{cmtt10 scaled\magstep2}
\newfont{\fourteenbf}{cmbx10 scaled\magstep2}
\newfont{\fourteeni}{cmmi10 scaled\magstep2}
\newfont{\fourteensy}{cmsy10 scaled\magstep2}
\newfont{\fourteenex}{cmex10 scaled\magstep2}
\newfont{\fourteenmsa}{msam10 scaled\magstep2}
\newfont{\fourteeneufm}{eufm10 scaled\magstep2}
\newfont{\fourteenmsb}{msbm10 scaled\magstep2}

\catcode`\@=11

\thinlines
\newskip\Einheit \Einheit=0.6cm
\newcount\xcoord \newcount\ycoord
\newdimen\xdim \newdimen\ydim \newdimen\PfadD@cke \newdimen\Pfadd@cke
\def\PfadDicke#1{\PfadD@cke#1 \divide\PfadD@cke by2 \Pfadd@cke\PfadD@cke \multiply\PfadD@cke by2}
\long\def\LOOP#1\REPEAT{\def\BODY{#1}\ITERATE}
\def\ITERATE{\BODY \let\next\ITERATE \else\let\next\relax\fi \next}
\let\REPEAT=\fi
\def\Punkt{\hbox{\raise-2pt\hbox to0pt{\hss\scriptsize$\bullet$\hss}}}
\def\DuennPunkt(#1,#2){\unskip
  \raise#2 \Einheit\hbox to0pt{\hskip#1 \Einheit
          \raise-2.5pt\hbox to0pt{\hss\normalsize$\bullet$\hss}\hss}}
\def\NormalPunkt(#1,#2){\unskip
  \raise#2 \Einheit\hbox to0pt{\hskip#1 \Einheit
          \raise-3pt\hbox to0pt{\hss\large$\bullet$\hss}\hss}}
\def\DickPunkt(#1,#2){\unskip
  \raise#2 \Einheit\hbox to0pt{\hskip#1 \Einheit
          \raise-4pt\hbox to0pt{\hss\Large$\bullet$\hss}\hss}}
\def\Kreis(#1,#2){\unskip
  \raise#2 \Einheit\hbox to0pt{\hskip#1 \Einheit
          \raise-4pt\hbox to0pt{\hss\Large$\circ$\hss}\hss}}
\def\Diagonale(#1,#2)#3{\unskip\leavevmode
  \xcoord#1\relax \ycoord#2\relax
      \raise\ycoord \Einheit\hbox to0pt{\hskip\xcoord \Einheit
         \unitlength\Einheit
         \line(1,1){#3}\hss}}
\def\AntiDiagonale(#1,#2)#3{\unskip\leavevmode
  \xcoord#1\relax \ycoord#2\relax \advance\xcoord by -0.05\relax
      \raise\ycoord \Einheit\hbox to0pt{\hskip\xcoord \Einheit
         \unitlength\Einheit
         \line(1,-1){#3}\hss}}
\def\Pfad(#1,#2),#3\endPfad{\unskip\leavevmode
  \xcoord#1 \ycoord#2 \thicklines\ZeichnePfad#3\endPfad\thinlines}
\def\ZeichnePfad#1{\ifx#1\endPfad\let\next\relax
  \else\let\next\ZeichnePfad
    \ifnum#1=1
      \raise\ycoord \Einheit\hbox to0pt{\hskip\xcoord \Einheit
         \vrule height\Pfadd@cke width1 \Einheit depth\Pfadd@cke\hss}%
      \advance\xcoord by 1
    \else\ifnum#1=2
      \raise\ycoord \Einheit\hbox to0pt{\hskip\xcoord \Einheit
        \hbox{\hskip-1pt\vrule height1 \Einheit width\PfadD@cke depth0pt}\hss}%
      \advance\ycoord by 1
    \else\ifnum#1=3
      \raise\ycoord \Einheit\hbox to0pt{\hskip\xcoord \Einheit
         \unitlength\Einheit
         \line(1,1){1}\hss}
      \advance\xcoord by 1
      \advance\ycoord by 1
    \else\ifnum#1=4
      \raise\ycoord \Einheit\hbox to0pt{\hskip\xcoord \Einheit
         \unitlength\Einheit
         \line(1,-1){1}\hss}
      \advance\xcoord by 1
      \advance\ycoord by -1
    \fi\fi\fi\fi
  \fi\next}
\def\hSSchritt{\leavevmode\raise-.4pt\hbox to0pt{\hss.\hss}\hskip.2\Einheit
  \raise-.4pt\hbox to0pt{\hss.\hss}\hskip.2\Einheit
  \raise-.4pt\hbox to0pt{\hss.\hss}\hskip.2\Einheit
  \raise-.4pt\hbox to0pt{\hss.\hss}\hskip.2\Einheit
  \raise-.4pt\hbox to0pt{\hss.\hss}\hskip.2\Einheit}
\def\vSSchritt{\vbox{\baselineskip.2\Einheit\lineskiplimit0pt
\hbox{.}\hbox{.}\hbox{.}\hbox{.}\hbox{.}}}
\def\DSSchritt{\leavevmode\raise-.4pt\hbox to0pt{%
  \hbox to0pt{\hss.\hss}\hskip.2\Einheit
  \raise.2\Einheit\hbox to0pt{\hss.\hss}\hskip.2\Einheit
  \raise.4\Einheit\hbox to0pt{\hss.\hss}\hskip.2\Einheit
  \raise.6\Einheit\hbox to0pt{\hss.\hss}\hskip.2\Einheit
  \raise.8\Einheit\hbox to0pt{\hss.\hss}\hss}}
\def\dSSchritt{\leavevmode\raise-.4pt\hbox to0pt{%
  \hbox to0pt{\hss.\hss}\hskip.2\Einheit
  \raise-.2\Einheit\hbox to0pt{\hss.\hss}\hskip.2\Einheit
  \raise-.4\Einheit\hbox to0pt{\hss.\hss}\hskip.2\Einheit
  \raise-.6\Einheit\hbox to0pt{\hss.\hss}\hskip.2\Einheit
  \raise-.8\Einheit\hbox to0pt{\hss.\hss}\hss}}
\def\SPfad(#1,#2),#3\endSPfad{\unskip\leavevmode
  \xcoord#1 \ycoord#2 \ZeichneSPfad#3\endSPfad}
\def\ZeichneSPfad#1{\ifx#1\endSPfad\let\next\relax
  \else\let\next\ZeichneSPfad
    \ifnum#1=1
      \raise\ycoord \Einheit\hbox to0pt{\hskip\xcoord \Einheit
         \hSSchritt\hss}%
      \advance\xcoord by 1
    \else\ifnum#1=2
      \raise\ycoord \Einheit\hbox to0pt{\hskip\xcoord \Einheit
        \hbox{\hskip-2pt \vSSchritt}\hss}%
      \advance\ycoord by 1
    \else\ifnum#1=3
      \raise\ycoord \Einheit\hbox to0pt{\hskip\xcoord \Einheit
         \DSSchritt\hss}
      \advance\xcoord by 1
      \advance\ycoord by 1
    \else\ifnum#1=4
      \raise\ycoord \Einheit\hbox to0pt{\hskip\xcoord \Einheit
         \dSSchritt\hss}
      \advance\xcoord by 1
      \advance\ycoord by -1
    \fi\fi\fi\fi
  \fi\next}
\def\Koordinatenachsen(#1,#2){\unskip
 \hbox to0pt{\hskip-.5pt\vrule height#2 \Einheit width.5pt depth1 \Einheit}%
 \hbox to0pt{\hskip-1 \Einheit \xcoord#1 \advance\xcoord by1
    \vrule height0.25pt width\xcoord \Einheit depth0.25pt\hss}}
\def\Koordinatenachsen(#1,#2)(#3,#4){\unskip
 \hbox to0pt{\hskip-.5pt \ycoord-#4 \advance\ycoord by1
    \vrule height#2 \Einheit width.5pt depth\ycoord \Einheit}%
 \hbox to0pt{\hskip-1 \Einheit \hskip#3\Einheit 
    \xcoord#1 \advance\xcoord by1 \advance\xcoord by-#3 
    \vrule height0.25pt width\xcoord \Einheit depth0.25pt\hss}}
\def\Gitter(#1,#2){\unskip \xcoord0 \ycoord0 \leavevmode
  \LOOP\ifnum\ycoord<#2
    \loop\ifnum\xcoord<#1
      \raise\ycoord \Einheit\hbox to0pt{\hskip\xcoord \Einheit\Punkt\hss}%
      \advance\xcoord by1
    \repeat
    \xcoord0
    \advance\ycoord by1
  \REPEAT}
\def\Gitter(#1,#2)(#3,#4){\unskip \xcoord#3 \ycoord#4 \leavevmode
  \LOOP\ifnum\ycoord<#2
    \loop\ifnum\xcoord<#1
      \raise\ycoord \Einheit\hbox to0pt{\hskip\xcoord \Einheit\Punkt\hss}%
      \advance\xcoord by1
    \repeat
    \xcoord#3
    \advance\ycoord by1
  \REPEAT}
\def\Label#1#2(#3,#4){\unskip \xdim#3 \Einheit \ydim#4 \Einheit
  \def\lo{\advance\xdim by-.5 \Einheit \advance\ydim by.5 \Einheit}%
  \def\llo{\advance\xdim by-.25cm \advance\ydim by.5 \Einheit}%
  \def\loo{\advance\xdim by-.5 \Einheit \advance\ydim by.25cm}%
  \def\o{\advance\ydim by.25cm}%
  \def\ro{\advance\xdim by.5 \Einheit \advance\ydim by.5 \Einheit}%
  \def\rro{\advance\xdim by.25cm \advance\ydim by.5 \Einheit}%
  \def\roo{\advance\xdim by.5 \Einheit \advance\ydim by.25cm}%
  \def\l{\advance\xdim by-.30cm}%
  \def\r{\advance\xdim by.30cm}%
  \def\lu{\advance\xdim by-.5 \Einheit \advance\ydim by-.6 \Einheit}%
  \def\llu{\advance\xdim by-.25cm \advance\ydim by-.6 \Einheit}%
  \def\luu{\advance\xdim by-.5 \Einheit \advance\ydim by-.30cm}%
  \def\u{\advance\ydim by-.30cm}%
  \def\ru{\advance\xdim by.5 \Einheit \advance\ydim by-.6 \Einheit}%
  \def\rru{\advance\xdim by.25cm \advance\ydim by-.6 \Einheit}%
  \def\ruu{\advance\xdim by.5 \Einheit \advance\ydim by-.30cm}%
  #1\raise\ydim\hbox to0pt{\hskip\xdim
     \vbox to0pt{\vss\hbox to0pt{\hss$#2$\hss}\vss}\hss}%
}
\catcode`\@=12

\def\ldreieck{\bsegment
  \rlvec(0.866025403784439 .5) \rlvec(0 -1)
  \rlvec(-0.866025403784439 .5)  
  \savepos(0.866025403784439 -.5)(*ex *ey)
        \esegment
  \move(*ex *ey)
        }
\def\rdreieck{\bsegment
  \rlvec(0.866025403784439 -.5) \rlvec(-0.866025403784439 -.5)  \rlvec(0 1)
  \savepos(0 -1)(*ex *ey)
        \esegment
  \move(*ex *ey)
        }
\def\rhombus{\bsegment
  \rlvec(0.866025403784439 .5) \rlvec(0.866025403784439 -.5) 
  \rlvec(-0.866025403784439 -.5)  \rlvec(0 1)        
  \rmove(0 -1)  \rlvec(-0.866025403784439 .5) 
  \savepos(0.866025403784439 -.5)(*ex *ey)
        \esegment
  \move(*ex *ey)
        }
\def\RhombusA{\bsegment
  \rlvec(0.866025403784439 .5) \rlvec(0.866025403784439 -.5) 
  \rlvec(-0.866025403784439 -.5) \rlvec(-0.866025403784439 .5) 
  \savepos(0.866025403784439 -.5)(*ex *ey)
        \esegment
  \move(*ex *ey)
        }
\def\RhombusB{\bsegment
  \rlvec(0.866025403784439 .5) \rlvec(0 -1)
  \rlvec(-0.866025403784439 -.5) \rlvec(0 1) 
  \savepos(0 -1)(*ex *ey)
        \esegment
  \move(*ex *ey)
        }
\def\RhombusC{\bsegment
  \rlvec(0.866025403784439 -.5) \rlvec(0 -1)
  \rlvec(-0.866025403784439 .5) \rlvec(0 1) 
  \savepos(0.866025403784439 -.5)(*ex *ey)
        \esegment
  \move(*ex *ey)
        }
\def\RhombusAsh{\bsegment
  \rlvec(0.866025403784439 .5) \rlvec(0.866025403784439 -.5) 
  \rlvec(-0.866025403784439 -.5) \rlvec(-0.866025403784439 .5) 
  \lfill f:.8
  \savepos(0.866025403784439 -.5)(*ex *ey)
        \esegment
  \move(*ex *ey)
        }
\def\RhombusBsh{\bsegment
  \rlvec(0.866025403784439 .5) \rlvec(0 -1)
  \rlvec(-0.866025403784439 -.5) \rlvec(0 1) 
  \lfill f:0.2
  \savepos(0 -1)(*ex *ey)
        \esegment
  \move(*ex *ey)
        }
\def\RhombusCsh{\bsegment
  \rlvec(0.866025403784439 -.5) \rlvec(0 -1)
  \rlvec(-0.866025403784439 .5) \rlvec(0 1) 
  \savepos(0.866025403784439 -.5)(*ex *ey)
        \esegment
  \move(*ex *ey)
        }
\def\hdSchritt{\bsegment
  \lpatt(.05 .13)
  \rlvec(0.866025403784439 -.5) 
  \savepos(0.866025403784439 -.5)(*ex *ey)
        \esegment
  \move(*ex *ey)
        }
\def\vdSchritt{\bsegment
  \lpatt(.05 .13)
  \rlvec(0 -1) 
  \savepos(0 -1)(*ex *ey)
        \esegment
  \move(*ex *ey)
        }

\def\odSchritt{\bsegment
  \lpatt(.05 .13)
  \rlvec(-0.866025403784439 -.5) 
  \savepos(-0.866025403784439 -.5)(*ex *ey)
        \esegment
  \move(*ex *ey)
        }

\def\hdaSchritt{\bsegment
  \lpatt(.05 .13)
  \rlvec(-0.866025403784439 .5) 
  \savepos(-0.866025403784439 .5)(*ex *ey)
        \esegment
  \move(*ex *ey)
        }
\def\vdaSchritt{\bsegment
  \lpatt(.05 .13)
  \rlvec(0 1) 
  \savepos(0 1)(*ex *ey)
        \esegment
  \move(*ex *ey)
        }

\def\odaSchritt{\bsegment
  \lpatt(.05 .13)
  \rlvec(0.866025403784439 .5) 
  \savepos(0.866025403784439 .5)(*ex *ey)
        \esegment
  \move(*ex *ey)
        }

\def\ringerl(#1 #2){\move(#1 #2)\fcir f:0 r:.1}

\def\({\left(}
\def\){\right)}
\def\[{\left[}
\def\]{\right]}
\def\fl#1{\left\lfloor#1\right\rfloor}
\def\cl#1{\left\lceil#1\right\rceil}
\def\P{\mathcal P}
\def\om{\omega}
\def\po#1#2{(#1)_#2}

\numberwithin{equation}{section}

\newtheorem{theorem}{Theorem}
\newtheorem{lemma}[theorem]{Lemma} 
\newtheorem{cor}[theorem]{Corollary}
\newtheorem{conj}{Conjecture}

\theoremstyle{remark}
\newtheorem*{remark}{Remark}

\DeclareMathOperator{\h}{H}
\DeclareMathOperator{\sgn}{sgn}

\begin{document}

\newbox\Adr
\setbox\Adr\vbox{
\centerline{Department of Mathematics, Georgia Institute of Technology}
\centerline{Atlanta, GA 30332-0160}
%\centerline{E-mail: \footnotesize{\tt ciucu@math.gatech.edu}}
%\centerline{WWW: \footnotesize{\tt
%http://www.math.gatech.edu/\~{}ciucu}}
\vskip6pt
\centerline{Institut f\"ur Mathematik der Universit\"at Wien}
\centerline{Strudlhofgasse 4, A-1090 Wien, Austria}
%\centerline{E-mail: \footnotesize{\tt teisenko@radon.mat.univie.ac.at},
%{\tt kratt@pap.univie.ac.at}}
%\centerline{WWW: \footnotesize\tt http://radon.mat.univie.ac.at/People/kratt}
\vskip6pt
\centerline{Department of Mathematics, Columbia University}
\centerline{New York, NY 10027}
%\centerline{E-mail: \footnotesize{\tt zare@math.columbia.edu}}
%\centerline{WWW: \footnotesize{\tt http://www.math.columbia.edu/\~{}zare}}
}

\title[Lozenge tilings of hexagons with a central 
triangular hole]{Enumeration of lozenge tilings of hexagons with a central 
triangular hole}

\author[M. Ciucu, T. Eisenk\"olbl, C.~Krattenthaler and D. Zare]{%
M. Ciucu$^\dagger$, T. Eisenk\"olbl, C.~Krattenthaler$^\ddagger$, D. Zare\\
[18pt]\box\Adr}

\address{Department of Mathematics, Georgia Institute of Technology.
\newline
E-mail: {\tt ciucu@math.gatech.edu}\\
WWW: {\tt http://www.math.gatech.edu/\~{}ciucu}}

\address{Institut f\"ur Mathematik der Universit\"at Wien,
Strudlhofgasse 4, A-1090 Wien, Austria.
E-mail: {\tt teisenko@radon.mat.univie.ac.at,
KRATT@Pap.Univie.Ac.At}\newline
WWW: \tt http://radon.mat.univie.ac.at/People/kratt}

\address{Department of Mathematics, Columbia University, New York, 
NY 10027.\newline
E-mail: {\tt zare@math.columbia.edu},
WWW: {\tt http://www.math.columbia.edu/\~{}zare}}

\thanks{$^\dagger$ Research partially supported by a membership at the
Institute for Advanced Study and by NSF grant DMS 9802390}

\thanks{$^\ddagger$ Research partially supported by the Austrian
Science Foundation FWF, grant P12094-MAT and P13190-MAT}

\subjclass{Primary 05A15;
 Secondary 05A17 05A19 05B45 33C20 52C20}

\keywords{lozenge tilings, rhombus tilings, plane partitions,
determinants, nonintersecting lattice paths}

\begin{abstract}
We deal with unweighted and weighted
enumerations of lozenge tilings of a hexagon with side
lengths $a,b+m,c,a+m,b,c+m$, where an equilateral triangle of
side length $m$ has been removed from the center. We give closed formulas
for the plain enumeration and for a certain $(-1)$-enumeration of
these lozenge tilings. In the case that $a=b=c$,
we also provide closed formulas for certain
weighted enumerations of those lozenge tilings that are
cyclically symmetric. 
For $m=0$, the latter formulas specialize to statements 
about weighted enumerations of cyclically symmetric plane partitions.
One such specialization gives a proof of a conjecture of Stembridge
on a certain weighted count of
cyclically symmetric plane partitions. 
The tools employed in our proofs are nonstandard applications of 
the theory of
nonintersecting lattice paths and determinant
evaluations. In particular, we evaluate the determinants $\det_{0\le
i,j\le n-1}\big(\om \delta_{ij}+\binom {m+i+j}j\big)$, where $\om$ is any
$6$th root of unity. These determinant evaluations are variations of
a famous result due to Andrews 
(Invent\@. Math. {\bf 53} (1979), 193--225),
which corresponds to $\om=1$.
\end{abstract}

\maketitle

\begin{section}{Introduction}
Let $a$, $b$ and $c$ be positive integers, and consider a 
hexagon with side lengths $a,b,c,a,b,\break c$ (in cyclic order) and 
angles of $120^\circ$.
It is well-known that the total number of lozenge\footnote{\ Here and 
in the following, by a lozenge we mean a rhombus with side 
lengths 1 and angles of $60^\circ$ and $120^\circ$.}
tilings of such a hexagon equals
\begin{align} \label{box}
\frac {\h(a)\h(b)\h(c)\h(a+b+c)} {\h(a+b)\h(b+c)\h(c+a)},
\end{align}
where $\h(n)$ stands for the ``hyperfactorial" $\prod _{k=0} ^{n-1}k!$. 
This follows from a bijection (cf\@. \cite{David-Tomei}) 
between such lozenge tilings and plane partitions contained in 
an $a\times b\times c$ box, and from MacMahon's enumeration 
\cite[Sec.~429, $q\rightarrow 1$; proof in Sec.~494]{MM}
of the latter.

In \cite{Propp} (see also \cite{Propp2}), Propp posed several problems regarding
``incomplete'' hexagons.
For example, Problem~2 in \cite{Propp} (and \cite{Propp2})
asks for the number of lozenge tilings of a hexagon with side lengths
$n,n+1,n,n+1,n,n+1$
with the central unit triangle removed. 
This problem was solved in \cite[Theorem~1]{CiucAI}, 
\cite[Theorem~20]{HeGeAA} and \cite[Theorem~1]{OkKr} (the
most general result, for a hexagon with side lengths
$a,b+1,c,a+1,b,c+1$, being contained in \cite{OkKr}). In
\cite{CiucAH}, 
the first author considers the case when a larger triangle (in fact,
possibly several) is removed.
However, in contrast to \cite{OkKr}, 
the results in \cite{CiucAH} 
assume that the hexagon has a reflective symmetry, i.e., that $b=c$.

Continuing this line of research, in this paper we address the general
case, when no symmetry axis is required.
We consider hexagons of
sides $a,b+m,c,a+m,b,c+m$ (in clockwise order) with an equilateral 
triangle of
side $m$ removed from the center
(see Figures~\ref{hex} and \ref{hex2} for examples).
We call this triangle the 
{\it core}, and the leftover region, denoted $C_{a,b,c}(m)$, a 
{\it cored hexagon}.

To define $C_{a,b,c}(m)$ precisely,
we need to specify what position of the core is the ``central'' 
one. Let $s$ be a side of the core, and let $u$ and $v$ be the sides of the
hexagon parallel to it. The most natural definition (and the one that we are 
going to adopt) would require that the distance between $s$ and $u$ is the 
same as the distance between $v$ and the vertex of the core opposite $s$, for 
all three choices of $s$.

However, since the sides of the core have to be along lines of the 
underlying triangular lattice, it is easy to see that this can be 
achieved only if $a$, $b$ and $c$ have the same parity
(Figure~\ref{hex} illustrates  
such a case); in that case, we define this to be the position of the core. 
On the other hand, if for instance $a$
has parity different from that of $b$ and $c$, the triangle satisfying the
above requirements would only have one side along a lattice line, while
each of the remaining two extends midway between two consecutive lattice
lines (this can be seen from Figure~\ref{hex2}). To resolve this, we translate this
central triangle half a unit towards the side of the hexagon of length $b$,
in a direction parallel to the side of length $a$, and define this to be the
position of the core in this case.

Note that, when translating the central triangle, there is no ``natural''
reason to do it in the sense we chose: we could have just as well
chosen the opposite sense, obtaining an alternative (and not less central)
definition of the core. However, it is easy to see that the alternative
definition does not lead to new regions: it generates the same region that we
obtain by swapping $b$ and $c$ in our definition.  
(In fact, this ambiguity in choosing the center will be used effectively
in Section~\ref{sec:Comments}, see Theorem~\ref{thm:Watson3} 
and the paragraph preceding it.)

\begin{figure}
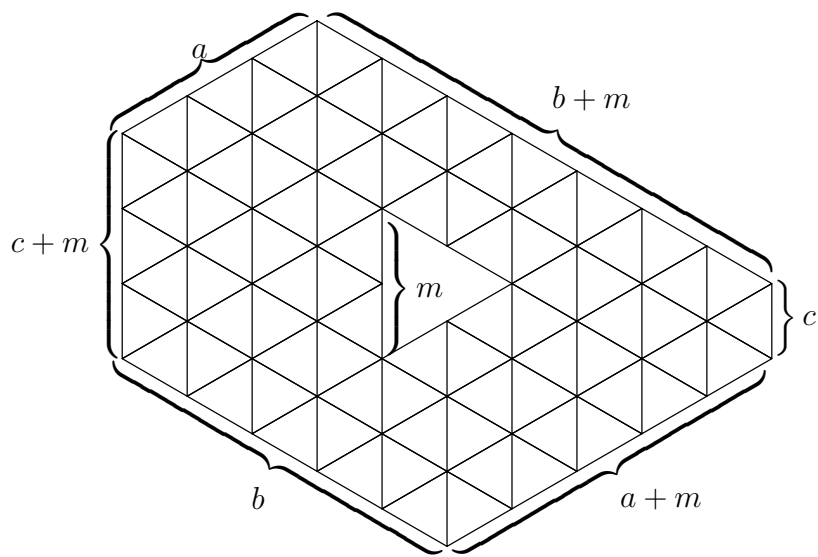
 
\centertexdraw{
\drawdim truecm \linewd.02
\rhombus \rhombus \rhombus \rhombus \rhombus \rhombus \rhombus
\ldreieck
\move(-.866025 -.5)
 \rhombus \rhombus \rhombus \rhombus \rhombus \rhombus \rhombus \rhombus 
\move(-1.7305 -1)
\rhombus \rhombus \rhombus \ldreieck \rmove(.866025 -.5)
\rhombus \rhombus \rhombus 
\move(-1.7305 -1)
\rdreieck \rhombus \rhombus \rhombus \ldreieck \rhombus \rhombus \rhombus 
\move(-1.7305 -2)
\rdreieck \rhombus \rhombus \rhombus \rhombus \rhombus \rhombus 
\move(-1.7305 -3)
\rdreieck \rhombus \rhombus \rhombus \rhombus \rhombus 
\htext(-.8 0){$a$}
\htext(4 -.7){$b+m$}
\htext(6.9 -3.95){$\left. \vbox{\vskip.6cm} \right\} c$}
\htext(-3.2 -4){$c+m \left\{ \vbox{\vskip1.6cm} \right.$}
\htext(0 -6){$b$}
\htext(4.9 -6){$a+m$}
\htext(1.7 -3.95){$\left. \vbox{\vskip1cm} \right\}m$}
\rtext td:60 (4 -1.3) {$\left. \vbox{\vskip3.6cm} \right\} $}
\rtext td:-60 (-.8 0){$\left\{ \vbox{\vskip1.6cm} \right. $}
\rtext td:-60 (4.6 -5.3) {$\left. \vbox{\vskip2.5cm} \right\} $}
\rtext td:60 (0.3 -5.6){$\left\{ \vbox{\vskip2.7cm}\right. $}
}
\caption{\small
Position of the core when $a$, $b$ and $c$ have the same parity: $C_{3,5,1}(2)$}
\label{hex}
\end{figure}

\newbox\hexshbox
\setbox\hexshbox\hbox{\small 
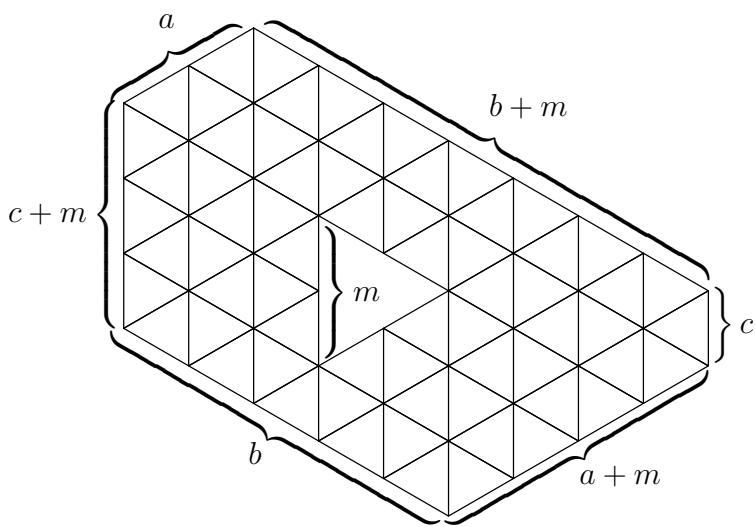
Position of the core when $a$, $b$ and $c$ have mixed parities: $C_{2,5,1}(2)$}
\begin{figure} 
\centertexdraw{
\drawdim truecm \linewd.02
\rhombus \rhombus \rhombus \rhombus \rhombus \rhombus \rhombus
\ldreieck
\move(-.866025 -.5)
 \rhombus \rhombus \rhombus \rhombus \rhombus \rhombus \rhombus \rhombus 
\move(-.866025 -.5)
\rdreieck
\rhombus \rhombus \ldreieck \rmove(.866025 -.5)
\rhombus \rhombus \rhombus 
\move(-.866025 -1.5)
\rdreieck \rhombus \rhombus \ldreieck \rhombus \rhombus \rhombus 
\move(-.866025 -2.5)
\rdreieck  \rhombus \rhombus \rhombus \rhombus \rhombus 
\move(-1.7305 -3)
\htext(-.4 .5){$a$}
\htext(4 -.7){$b+m$}
\htext(6.9 -3.95){$\left. \vbox{\vskip.6cm} \right\} c$}
\htext(-2.4 -3.5){$c+m \left\{ \vbox{\vskip1.6cm} \right.$}
\htext(.8 -5.3){$b$}
\htext(5.2 -5.6){$a+m$}
\htext(1.7 -3.9){$\left. \vbox{\vskip1cm} \right\}m$}
\rtext td:60 (4 -1.3) {$\left. \vbox{\vskip3.6cm} \right\} $}
\rtext td:-60 (-.4 0.2){$\left\{ \vbox{\vskip1.1cm} \right. $}
\rtext td:-60 (5 -5) {$\left. \vbox{\vskip2.2cm} \right\} $}
\rtext td:60 (1.1 -5.1){$\left\{ \vbox{\vskip2.7cm}\right. $}
}
\caption{\unhbox\hexshbox}
\label{hex2}
\end{figure}

Our main results, given in Theorems~\ref{enum} and
\ref{shifted} below, provide explicit formulas for the total number of lozenge
tilings of such a cored hexagon
(see Figures~\ref{-1fig} and \ref{tiling}.a for examples of such
tilings).
Remarkably, the results can be expressed in {\em closed form},
more precisely, as quotients of products
of hyperfactorials (completely
analogous to formula \eqref{box}), thus providing an infinite family
of enumerations which contains
MacMahon's ``box formula" \eqref{box} as
a special case.
For the statement of the theorems, it is convenient to extend the
definition of hyperfactorials to half-integers (i.e.,
odd integers divided by~2):
$$\h(n):=\begin{cases}
\prod _{k=0} ^{n-1}{\Gamma(k+1)}\quad &\text {for $n$ an integer,}\\
\prod _{k=0} ^{n-\frac {1} {2}}{\Gamma(k+\frac {1} {2})} \quad &\text 
{for $n$ a half-integer}.
\end{cases}
$$ 

Now we are able to state our theorems. The first result addresses the
case that $a$, $b$ and $c$ have the same parity. Let $L(R)$ stand for the
number of lozenge tilings of the region $R$.

\begin{theorem} \label{enum}
Let $a,b,c,m$ be nonnegative integers, $a,b,c$ having the same parity. 
The number of lozenge tilings of a hexagon with sides
$a,b+m,c,a+m,b,c+m$, with an equilateral triangle of
side $m$ removed from its center (see Figure~\ref{hex} for an example)
is given by
\begin{multline} \label{eq:enum}
L(C_{a,b,c}(m))=\\
\frac {\h(a + m)\h(b + m)\h(c + m)\h(a + b + c + m)
} 
{\h(a + b + m)\h(a + c + m)\h(b + c + m)
}
\frac {\h(m + \left \lceil {\frac{a + b + c}{2}} \right \rceil)
\h(m + \left \lfloor {\frac{a + b + c}{2}} \right \rfloor)
} {\h({\frac{a + b}{2}} + m)     \h({\frac{a + c}{2}} + m)\h({\frac{b + c}{2}} + m)
}
\\
\times\frac {\h(\left \lceil {\frac{a}{2}} \right \rceil)
\h(\left \lceil {\frac{b}{2}} \right \rceil)
     \h(\left \lceil {\frac{c}{2}} \right \rceil)
     \h(\left \lfloor {\frac{a}{2}} \right \rfloor)\,
     \h(\left \lfloor {\frac{b}{2}} \right \rfloor)\,
     \h(\left \lfloor {\frac{c}{2}} \right \rfloor)\,
} 
{\h({\frac{m}{2}} + \left \lceil {\frac{a}{2}} \right \rceil)\,
     \h({\frac{m}{2}} + \left \lceil {\frac{b}{2}} \right \rceil)\,
     \h({\frac{m}{2}} + \left \lceil {\frac{c}{2}} \right \rceil)\,
\h({\frac{m}{2}} + \left \lfloor {\frac{a}{2}} \right \rfloor)\,
     \h({\frac{m}{2}} + \left \lfloor {\frac{b}{2}} \right \rfloor)\,
     \h({\frac{m}{2}} + \left \lfloor {\frac{c}{2}} \right \rfloor)\,
     }\\
\times
\frac {\h(\frac{m}{2})^2 \h({\frac{a + b + m}{2}})^2 
\h({\frac{a + c + m}{2}})^2 \h({\frac{b + c +
m}{2}})^2
} 
{\h({\frac{m}{2}} + \left \lceil {\frac{a + b + c}{2}} \right \rceil)
\h({\frac{m}{2}} + \left \lfloor {\frac{a + b + c}{2}} \right \rfloor)
\h({\frac{a + b}{2}})\h({\frac{a + c}{2}})\h({\frac{b + c}{2}})
               }.
\end{multline}
\end{theorem}

Clearly, formula \eqref{eq:enum} reduces to \eqref{box} for $m=0$
(as it should).
The special case $m=1$ has been obtained earlier in \cite[Theorem~1]{OkKr}.

The corresponding result for the case when
$a$, $b$ and $c$ do not
have the same parity reads as follows.

\begin{theorem} \label{shifted}
Let $a,b,c,m$ be nonnegative integers, with
$a$ of parity different from the parity of $b$ and $c$.
The number of lozenge tilings of a hexagon with sides
$a,b+m,c,a+m,b,c+m$, with 
the ``central" (in the sense described above) triangle of side $m$ removed
(see Figure~\ref{hex2} for an example) is given by
\begin{multline} \label{eq:shifted}
L(C_{a,b,c}(m))=\\
\frac 
{\h\(a + m\)\h\(b + m\)\h\(c + m\)\h\(a + b + c + m\)
\h\(m + \cl{\frac{a + b + c} {2}}\)\h\(m + \fl{\frac{a + b + c} {2}}\)} 
{\h\(a + b + m\)\h\(a + c + m\)\h\(b + c + m\)
\h\(\fl{\frac{a + c} {2}}+m\)\h\(\frac{b + c} {2} +
m\)\h\(\cl{\frac{a+b}{2}}+m\)
}\\
\times
\frac {\h\(\frac m 2\)^2\h\(\cl{\frac{a} {2}}\)\h\(\cl{\frac{b} {2}}\)
    \h\(\cl{\frac{c} {2}}\)\h\(\fl{\frac{a} {2}}\)
\h\(\fl{\frac{b} {2}}\)\h\(\fl{\frac{c} {2}}\)  } 
{\h\(\frac{m} {2} + \cl{\frac{a} {2}}\)\h\(\frac{m} {2} +
\cl{\frac{b} {2}}\)\h\(\frac{m} {2} + \cl{\frac{c} {2}}\)
\h\(\frac{m} {2} + \fl{\frac{a} {2}}\)\h\(\frac{m} {2} +
\fl{\frac{b} {2}}\)\h\(\frac{m} {2} + \fl{\frac{c} {2}}\)
}\\
\times
\frac 
{\h\(\cl{\frac {a + b} 2} + \frac{m} {2}\)
\h\(\fl{\frac{a + b}{2}} + \frac{m} {2}\)
\h\(\fl{\frac{a + c} {2}} + \frac{m} {2}\)\h\(\cl{\frac{a + c}
{2}} + \frac{m} {2}\)
    \h\(\frac{b + c} {2} + \frac{m} {2}\)^2
} 
{\h\(\frac{m} {2} + \cl{\frac{a + b + c} {2}}\)
\h\(\frac{m} {2} + \fl{\frac{a + b + c} {2}}\)
\h\(\fl{\frac{a + b} {2}}\)\h\(\cl{\frac{a + c}
{2}}\) \h\(\frac{b + c} {2}\)
}.
\end{multline}
\end{theorem}

Again, formula \eqref{eq:enum} reduces to \eqref{box} for $m=0$.
The special case $m=1$ has been obtained earlier in \cite[Theorem~4]{OkKr}.

Given the explicit results in Theorems~\ref{enum} and \ref{shifted},
it is routine to determine, using the Euler--MacLaurin summation
formula, the asymptotic behavior of
the number of lozenge tilings of a cored hexagon.
For instance, when $a$, $b$ and $c$ have the same parity we obtain the 
following result.

\begin{cor} \label{asy1}
Let $a,b,c,m,n$ be nonnegative integers, $a,b,c$ having the same parity. 
The number of lozenge tilings of a hexagon with sides
$an,(b+m)n,cn,(a+m)n,bn,(c+m)n$, with an equilateral triangle of
side $mn$ removed from its center, is asymptotically 
given by
$$L(C_{an,bn,cn}(mn))\sim e^{kn^2},\ \ \ n\to\infty,$$
where
\begin{multline} \label{eq:asy1}
k= \tfrac {(a + m)^2}2\log(a + m) + \tfrac{(b + m)^2} 2\log(b + m) + 
\tfrac{(c + m)^2}2\log(c + m)\\
+ \tfrac {(a + b + c + m)^2} 2 \log(a + b + c + m)
 +2(m + {\tfrac{a + b + c}{2}} )^2\log(m + {\tfrac{a + b + c}{2}} )
\\
+ 2({\tfrac{b}{2}} )^2\log(\tfrac b 2)+ 2({\tfrac{c}{2}} )^2\log(\tfrac c 2)
+ 2({\tfrac{a}{2}} )^2\log(\tfrac a 2)+ (\tfrac{m}{2})^2\log(m) \\
-
\(
\tfrac {3} {4}(a + b + m)^2\log(a + b + m) + 
\tfrac {3} {4}(a + c + m)^2\log(a + c + m) \right.\\
+\tfrac {3} {4}(b + c + m)^2\log(b + c + m) \\
+ ({\tfrac{a + b}{2}} + m)^2\log({\tfrac{a + b}{2}} + m) + 
({\tfrac{a + c}{2}} +m)^2\log({\tfrac{a + c}{2}} + m) + 
({\tfrac{b + c}{2}} + m)^2\log({\tfrac{b + c}{2}} + m)\\
\left.
+ ({\tfrac{a + b}{2}})^2\log(a + b) 
+ ({\tfrac{a +c}{2}})^2\log(a + c) 
+({\tfrac{b + c}{2}})^2\log(b + c) 
\)\\
+(m^2+ a^2+b^2+c^2+\tfrac {3m(a+b+c)}2 +ab+bc+ca)\log 2
.
\end{multline}
\end{cor}
%Using the Euler summation formula, we obtain
%$$\log \h(n)=\frac {n^2} {2}\log n -\frac {3} {4} n^2 + \O(n\log n).$$

In addition to plain counts,
$(-1)$-enumerations of plane partitions, i.e., enumerations where
plane partitions are given a weight of $1$ or $-1$, according to
certain rules, have been found to possess remarkable properties (see
\cite{StemAL,StemAI}). 
Motivated in part by a conjectured $(-1)$-enumeration on cyclically symmetric
plane partitions due to Stembridge \cite{StemAZ},
in Section~\ref{sec:-1} we consider a 
$(-1)$-enumeration of the lozenge tilings of Theorems~\ref{enum} and
\ref{shifted}. 
The corresponding results are given in Theorems~\ref{-enum} and 
\ref{-shifted}.

In Section~\ref{sec:cyclically}, we restrict our attention to
cyclically symmetric lozenge tilings 
(i.e., tilings invariant under rotation by $120^\circ$)
of cored hexagons. Clearly, this makes 
sense only if $a=b=c$, i.e., for cored hexagons of the form $C_{a,a,a}(m)$.
The plain enumeration of such cyclically symmetric lozenge
tilings had already been considered in \cite[Theorem~3.2 and
Corollary~3.3]{CiKrAB}. 
We restate the result here as Theorem~\ref{thm:cycl}. 
We provide several additional results. 
Theorem~\ref{thm:cycl-1} concerns the $(-1)$-enumeration of such cyclically
symmetric lozenge tilings and some additional weighted enumerations of them, 
where each
lozenge tiling is weighted by some $6$th root of unity, according to
a certain rule (see the paragraph before Theorem~\ref{thm:cycl-1}
for the precise definition). 
In the special case 
$m=0$ we obtain results about
weighted enumerations of cyclically symmetric plane partitions 
(see Corollary~\ref{cor:cycl-1}). A particular case of
Corollary~\ref{cor:cycl-1} proves a conjecture of Stembridge
\cite[Case~9 on p.~6]{StemAZ} 
about a certain $(-1)$-enumeration of cyclically
symmetric plane partitions. (The first proof of this conjecture, by
totally different means, is due to Kuperberg 
\cite[last displayed equation on p.~27]{KupeAG}.) Our results 
also allow us to prove another conjecture 
on $(-1)$-enumeration of
cyclically symmetric plane partitions due to Stembridge
\cite[Case~10 on p.~7]{StemAZ}. In fact, we again prove a more
general result, namely a result on
cyclically symmetric lozenge tilings (see Theorem~\ref{thm:Case10}).

The remaining sections, Sections~\ref{sec:outline}--\ref{omsec}, are
devoted to the proofs of these results. For the proofs of 
Theorems~\ref{enum}--\ref{-shifted}, the enumeration
results for lozenge tilings without symmetry, we proceed as follows. 
First, we identify tilings with
certain families of nonintersecting
lattice paths (see Section~\ref{sec:nonint}). 
Then, a nonstandard application of 
the main theorem on nonintersecting lattice paths
\cite[Lemma~1]{LindAA}, \cite[Theorem~1]{gv} 
(restated here in Lemma~\ref{gv}) provides a
determinant for the weighted count of lozenge tilings (see
\eqref{mat1}, respectively \eqref{mat2}).
To be precise, the determinant gives the correct weighted count either 
only for even $m$ ($m$ being the side of the core)
or only for odd $m$, depending on whether we are considering
plain enumeration or $(-1)$-enumeration.
To cover the other case as well, we prove that the weighted 
count of lozenge tilings that we are
interested in is polynomial in $m$, so that it suffices to
determine this number only for one of the two possibilities, either
for even $m$ or for odd $m$. This is in turn achieved
by evaluating the aforementioned determinant (see 
Lemmas~\ref{lem:amgerade}--\ref{lem:shaomo}). 

The results on weighted enumerations of cyclically symmetric
lozenge tilings in Section~\ref{sec:cyclically} can be obtained in a
similar way.
We phrase the problem in terms of nonintersecting lattice paths,
and thus find determinants for these enumerations. The determinants 
have the form 
\begin{equation} \label{eq:omega}
\det_{0\le
i,j\le a-1}\(\om \delta_{ij}+\binom {m+i+j}j\),
\end{equation} 
where $\om$ is any
$6$th root of unity. These
determinants are remarkable. The case $\om=1$ occured first in the work
of Andrews on plane partitions. He evaluated the determinant
\eqref{eq:omega} in that case \cite[Theorem~8]{AndrAN} (restated here as 
Theorem~\ref{thm:Andrews}) in order to prove the ``weak Macdonald 
conjecture'' on counting cyclically symmetric
plane partitions. It had already been
observed in \cite[Sec.~3]{CiKrAB} that
Andrews' evaluation of \eqref{eq:omega} with $\om=1$ gives the number
of cyclically symmetric lozenge tilings of 
the cored hexagon $C_{a,a,a}(m)$.
We prove our {\it weighted\/}
enumerations of these lozenge tilings by evaluating the determinant
\eqref{eq:omega} when $\om$ is {\em any} $6$th root of unity (see
Theorems~\ref{thm:Zare1}--\ref{thm:om6}).

Our paper is structured as follows. In Section~\ref{sec:-1} we 
give the precise definition of our $(-1)$-enumeration of lozenge tilings, 
and we state the corresponding results (see Theorems~\ref{-enum} and
\ref{-shifted}).
In Section~\ref{sec:cyclically} we define
precisely our unusual weightings of cyclically symmetric lozenge
tilings. Theorems~\ref{thm:cycl} and \ref{thm:cycl-1}, 
Corollary~\ref{cor:cycl-1} and Theorem~\ref{thm:Case10}
state the corresponding results. The subsequent section,
Section~\ref{sec:outline}, gives the proofs of our enumeration
results in Theorems~\ref{enum}--\ref{thm:Case10}, 
leaving out, however, several details. 
These details are then worked out in later
sections. First of all, in Section~\ref{sec:nonint}, it is explained
how lozenge tilings correspond, in a one-to-one fashion, 
to families of nonintersecting lattice paths. 
We then employ the result of Lemma~\ref{gv}
to obtain, at least for 
every other value of $m$, a determinant for the weighted
count of lozenge tilings that we are interested in (see
Lemmas~\ref{lem:nonint} and \ref{lem:nonint2}). It is then
argued in Section~\ref{sec:poly} that this number is in fact
polynomial in $m$, so that the evaluation of the determinant in
Lemma~\ref{lem:nonint}, respectively Lemma~\ref{lem:nonint2}, 
suffices. The precise form of the
evaluation of the determinant in Lemma~\ref{lem:nonint} 
(again, a case-by-case analysis is necessary, depending on
the parity of $a$) is stated and proved in 
Section~\ref{sec:det} (see
Lemmas~\ref{lem:amgerade}--\ref{lem:-aung}), 
while the precise form of the
evaluation of the determinant in Lemma~\ref{lem:nonint2} 
is stated and proved in 
Section~\ref{sec:det2} (see
Lemmas~\ref{lem:shamgerade}--\ref{lem:shaomo}). 
Finally, in Section~\ref{sec:Zare1} 
we prove the determinant evaluation of Theorem~\ref{thm:Zare1},
in Section~\ref{corsec} the one in
Theorem~\ref{thm:om3}, and in Section~\ref{omsec} the one in
Theorem~\ref{thm:om6}. 
We conclude the article with some comments
concerning connections of this work with multiple
hypergeometric series and some open problems. These are the subject
of Section~\ref{sec:Comments}.

\end{section}

\begin{section}{$(-1)$-enumerations of lozenge tilings of cored
hexagons} \label{sec:-1}

\newbox\minusbox
\setbox\minusbox\hbox{\small
A lozenge tiling of a hexagon with sides $a=5$, $b=3$, $c=1$ and removed triangle 
of side length $m=2$.}
\begin{figure} 
\centertexdraw{
\drawdim truecm \linewd.02
\RhombusA \RhombusA \RhombusA \RhombusA \RhombusB \RhombusA
\move(-.866025 -.5)
\RhombusA \RhombusA \RhombusB \RhombusA \RhombusA \RhombusA 
\move(-1.7305 -1)
\RhombusB \RhombusA \RhombusB \RhombusB \RhombusA \RhombusA 
\move(-2.596 -1.5)
\RhombusB \RhombusB \RhombusB \RhombusA \RhombusA \RhombusA 
\move(-3.462 -2)
\RhombusB \RhombusB \RhombusB \RhombusA \RhombusA \RhombusA 
\move(0 -4)
\RhombusA \RhombusA 
\move(.866025 -3.5)
\RhombusA \RhombusA 
\move(-1.7305 -3)
\rlvec(.866025 -.5)
\move(0 -1)
\rlvec(0 -1)
\move(2.596 -1.5)
\rlvec(0 -1)
\move(4.33 -1.5)
\RhombusC
\move(0 -2)
\rlvec(.866025 -.5)

%\move(-2.165 -1.25)
%\vdSchritt \vdSchritt \vdSchritt \hdSchritt \hdSchritt \hhdSchritt
%\hodaSchritt \odaSchritt \odaSchritt \odaSchritt \odaSchritt
%\odaSchritt \vdaSchritt
%\move(-3.031 -1.75)
%\vdSchritt \vdSchritt \vdSchritt \hdSchritt \hdSchritt \hdSchritt
%\move(-1.298 -1.25)
%\hvdSchritt \hdSchritt \vdSchritt \vdSchritt \hhdSchritt \hodaSchritt 
%\odaSchritt
%\odaSchritt \odaSchritt \hvdaSchritt \hhdaSchritt \hdaSchritt \odSchritt 
%\hdaSchritt \hdaSchritt  \hodSchritt 

\linewd.15
%\move(0 -4)
%\rlvec(0 -2)
\move(1.73205 -3)
\rlvec(0.866 0.5)
\move(3.4641 -2)
\rlvec(0.866 0.5)
%\move(0 -2)
%\rlvec(-1.73205 1)
\linewd.02

\htext(-1.8 -.4){$a$}
\htext(3 -.1){$b+m$}
\htext(5.15 -3){$\left. \vbox{\vskip.5cm} \right\} c$}
\htext(-4.85 -5){$c+m \left\{ \vbox{\vskip1.6cm} \right.$}
\htext(-2.9 -6.4){$b$}
\htext(2.6 -5.5){$a+m$}
\htext(0 -4){$\left. \vbox{\vskip1.1cm} \right\}m$}
\rtext td:60 (3.1 -.9) {$\left. \vbox{\vskip2.75cm} \right\} $}
\rtext td:-60 (-1.6 -.6){$\left\{ \vbox{\vskip2.7cm} \right. $}
\rtext td:-60 (2 -4.7) {$\left. \vbox{\vskip3.7cm} \right\} $}
\rtext td:60 (-2.4 -6){$\left\{ \vbox{\vskip1.7cm}\right. $}
}
\caption{\unhbox\minusbox}
\label{-1fig}
\end{figure}

In this section we enumerate lozenge tilings of a cored hexagon with 
respect to a certain weight that assigns to each
lozenge tiling the value $1$ or $-1$. More precisely, fix a lozenge
tiling $T$ of the cored hexagon $C_{a,b,c}(m)$ (see Figures~\ref{hex} and
\ref{hex2} for examples of such regions, and
Figure~\ref{-1fig} for an example of a tiling; at this point, the 
thickness of edges is without significance).
Consider the side of the core which is parallel to the sides of
the hexagon of lengths $a$ and $a+m$ (in the figure this is the
bottommost side of the core). Extend this side of the triangle to
the right. 
Let $n(T)$ be the number of edges of lozenges of the tiling $T$ 
contained in the extended side (in Figure~\ref{-1fig} there are two
such edges, marked as thick segments). 
The statistic $n(T)$ becomes most transparent in 
the lattice path interpretation of lozenge tilings
that is going to be explained in Section~\ref{sec:nonint}, as it
counts exactly the number of paths which pass the core 
{\em on the right\/}. Furthermore, we shall see  
in Section~\ref{sec:cyclically} that in the plane partitions case,
i.e., in the case $m=0$ (when the core shrinks to a point), the
statistic $n(T)$ has a very natural meaning as well (see the remarks
after Theorem~\ref{thm:cycl-1}).

In the $(-1)$-enumeration, which is the
subject of the following two theorems, each lozenge tiling $T$ is
weighted by $(-1)^{n(T)}$. Let $L^{-1}(R)$ be the weighted count of
lozenge tilings of region $R$ under the above weight.

\begin{theorem} \label{-enum}
Let $a,b,c,m$ be nonnegative integers.
If all of $a$, $b$ and $c$ are even, then
the weighted count $\sum (-1)^{n(T)}$, summed over all 
lozenge tilings $T$ of a hexagon with sides
$a,b+m,c,a+m,b,c+m$, with an equilateral triangle of
side length $m$ removed from its center (see Figure~\ref{hex}) is given by
\begin{multline} \label{eq:-enum}
L^{-1}(C_{a,b,c}(m))=\\
(-1)^{a/2}\frac 
{\h(a + m)\,\h(b + m)\,\h(c + m)\,\h(a + b + c + m)\,} 
{\h(a + b + m)\,\h(a + c + m)\,\h(b + c + m)\,
} 
\\
\times
\frac 
{\h({\frac a 2})^2\,{{\h({\frac b 2})}^2}\,{{\h({\frac c 2})}^2}\,
\h({\frac {m-1} 2})\,\h({\frac {m+1} 2})
} 
{\h({\frac a 2} + {\frac {m-1} 2})\,\h({\frac b 2} + {\frac {m-1} 2})\,
     \h({\frac c 2} + {\frac {m-1} 2})\,\h({\frac a 2} + {\frac
{m+1} 2})\,\h({\frac b 2} + {\frac {m+1} 2})\,
     \h({\frac c 2} + {\frac {m+1} 2})}
\\
\times
\frac {\h(\frac {a + b+m-1} 2)\,
\h(\frac {a + b + m + 1} 2)\,
     \h({\frac {a + c + m - 1} 2})\,\h({\frac {a + c + m +1} 2})\,
     \h({\frac {b + c + m-1} 2})\,\h({\frac {b + c + m+1} 2})} 
   {\h({\frac {a + b} 2})\,\h({\frac {a + c} 2})\,\h({\frac {b + c} 2})\,
     \h({\frac {a + b} 2} + m)\,     \h({\frac {a + c} 2} + m)\,
   \h({\frac {b + c} 2} + m)}\\
\times
\frac {     \h({\frac {a + b + c} 2} + m)^2} 
{  \h({\frac {a + b + c} 2} + {\frac {m-1} 2})\,
   \h({\frac {a + b + c} 2} + {\frac {m+1} 2})}.
\end{multline}
For $a,b,c$ all odd, the $(-1)$-enumeration equals zero.
\end{theorem}

The analogous theorem for the case when $a$ has a parity different
from the parity of $b$ and $c$ reads as follows.

\begin{theorem} \label{-shifted}
Let $a,b,c,m$ be nonnegative integers, 
$a$ of parity different from the parity of $b$ and $c$.
The weighted count $\sum {(-1)^{n(T)}} $, summed over all
lozenge tilings $T$ of a hexagon 
with sides
$a,b+m,c,a+m,b,c+m$, with an equilateral triangle of
side length $m$ removed that is ``central" in the sense 
that was described in the Introduction (see Figure~\ref{hex2}), equals
\begin{multline}
L^{-1}(C_{a,b,c}(m))=\\
(-1)^{\cl{a/2}}\frac {\h(a + m)\h(b + m)\h(c + m)\h(a + b + c + m) 
} 
{\h(a + b + m)\h(a + c + m)\h(b + c + m)}\\
\times
\frac {\h(\fl{\frac {a + b + c} 2} +m)
\h(\cl{\frac {a + b + c} 2} +m)} 
{\h(\frac {a + b+1} 2 + m)
\h(\frac {a + c-1} 2 + m)\h(\frac {b + c} 2 + m)}\\
\times
\frac 
{\h(\fl{\frac a 2}) \h(\cl{\frac a 2}) 
 \h(\fl{\frac b 2)} \h(\cl{\frac b 2})
 \h(\fl{\frac c 2}) \h(\cl{\frac c 2}) 
 \h(\frac {m-1} 2)  \h(\frac {m+1} 2) }  
{\h(\frac {m-1} {2}+\fl{\frac {a + 1} 2})
\h(\frac {m+1} {2}+\cl{\frac {a -1} 2})
\h(\frac {m-1} {2}+\fl{\frac {b + 1} 2}) 
\h(\frac {m+1} {2}+\cl{\frac{b-1}2})
\h(\frac {m-1} {2}+\fl{\frac {c + 1} 2}) 
%(c-1)-Term ist in der naechsten Zeile
}\\
\times
\frac{\h(\frac {a + b + m} 2)^2
     \h(\frac {a + c + m} 2)^2 \h(\frac {b + c + m - 1} 2)
     \h(\frac {b + c + m+1} 2)
     }
   {\h(\frac {m+1} {2}+\cl{\frac{c-1}2})
\h({\frac {a + b-1} 2})\,\h({\frac {a + c+1} 2}) \h({\frac {b + c} 2})
    \h(\frac {m-1} {2}+\fl{\frac {a + b + c + 1} 2}) 
    \h(\frac {m+1} {2}+\cl{\frac {a + b + c - 1} 2})
   }.
\end{multline}
\end{theorem}

\end{section}

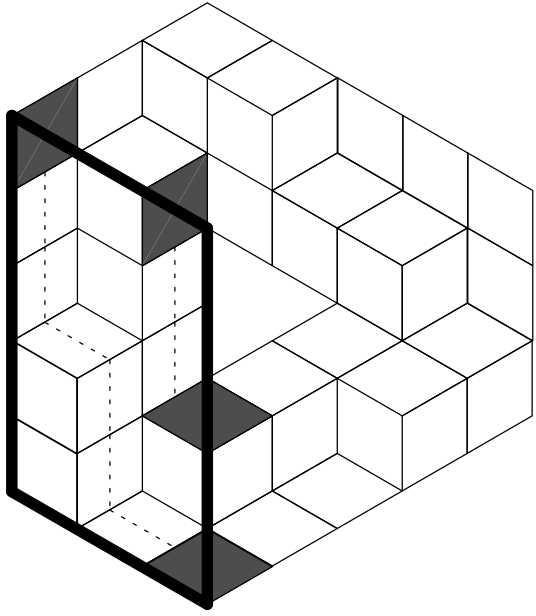
\begin{section}{Enumeration of cyclically symmetric lozenge tilings}
\label{sec:cyclically}

\newbox\cycbox
\setbox\cycbox\hbox{\small
A cyclically symmetric lozenge tiling of a hexagon with sides 
3, 5, 3, 5, 3, 5 
and core of size $m=2$.}
\begin{figure} 
\centertexdraw{
\drawdim truecm \linewd.02
\move(.866025 -2.5)
%\rlvec(0 -5) \rlvec(-2.598 1.5) \lfill f:.5
%\rlvec(0 5) \rlvec (2.598 -1.5) \lfill f:.5
\move(0 0)
\RhombusA \RhombusA \RhombusB\RhombusA \RhombusA \RhombusB\RhombusA \RhombusB
\move(-.866025 -.5)
\RhombusB\RhombusA \RhombusB\RhombusB\RhombusB\RhombusA \RhombusB\RhombusA 
\move(-1.7305 -1)
\RhombusB\RhombusB\RhombusB\RhombusA \RhombusB\RhombusB\RhombusA \RhombusA 
\move(.866025 -4.5)
\RhombusA \RhombusB\RhombusA 
\move(1.7305 -4)
\RhombusA \RhombusA \RhombusB
\move(0 0)  \RhombusC
\move(0.866025 -1.5) \RhombusC \RhombusC \RhombusC
\move(2.598 -.5) \RhombusC \RhombusC \RhombusC
\move(4.3301 -2.5) \RhombusC
\move(-1.73205 -4) \RhombusC
\move(-1.73205 -5) \RhombusC
\move(0 -3) \rlvec(-.866025 .5)

\move(.433012 -2.75)
%\hodaSchritt \odaSchritt \odaSchritt \hvdaSchritt \hhdaSchritt
%\hdaSchritt \hdaSchritt \hodSchritt \hvdSchritt 
\vdSchritt
\vdSchritt %\hhdSchritt
\move(-1.29904 -1.75)
%\hodaSchritt \odaSchritt \odaSchritt \vdaSchritt \odaSchritt 
%\odaSchritt \vdaSchritt \hvdaSchritt \hhdaSchritt \hdaSchritt
%\hdaSchritt \odSchritt \hdaSchritt \hdaSchritt \odSchritt \hodSchritt
%\hvdSchritt 
\vdSchritt \vdSchritt \hdSchritt \vdSchritt \vdSchritt
\hdSchritt %\hhdSchritt

\move(.866025 -4.5)  \rlvec(-.866025 -.5) \rlvec(.866025 -.5)
\lfill f:.3
\rlvec(.866025 .5) \rlvec(-.866025 .5)
\lfill f:.3

\move(.866025 -6.5)  \rlvec(-.866025 -.5) \rlvec(.866025 -.5)
\lfill f:.3
\rlvec(.866025 .5) \rlvec(-.866025 .5)
\lfill f:.3

\move(.866025 -1.5)  \rlvec(-.866025 -.5) \rlvec(0 -1)
\lfill f:.3
\rlvec(.866025 .5) \rlvec(0 1)
\lfill f:.3

\move(-.866025 -.5)  \rlvec(-.866025 -.5) \rlvec(0 -1)
\lfill f:.3
\rlvec(.866025 .5) \rlvec(0 1)
\lfill f:.3

\linewd.15
\move(.866025 -2.5) \rlvec(0 -5)
\move(.866025 -2.5) \rlvec(-2.598 1.5)
\move(-1.73205 -6) \rlvec(2.598 -1.5)
\move(-1.73205 -6) \rlvec(0 5)
%\move(2.598 -3.5) \rlvec(2.598 1.5)

%\htext(-.8 0){$a$}
%\htext(4 -1){$b+m$}
%\htext(6.9 -3.9){$\left. \vbox{\vskip.5cm} \right\} c$}
%\htext(-3.1 -4){$c+m \left\{ \vbox{\vskip1.6cm} \right.$}
%\htext(0 -6){$b$}
%\htext(5 -6){$a+m$}
%\htext(1.7 -4){$\left. \vbox{\vskip1.1cm} \right\}m$}
%\rtext td:60 (4 -1.3) {$\left. \vbox{\vskip3.6cm} \right\} $}
%\rtext td:-60 (-.8 0){$\left\{ \vbox{\vskip1.6cm} \right. $}
%\rtext td:-60 (4.6 -5.3) {$\left. \vbox{\vskip2.5cm} \right\} $}
%\rtext td:60 (0.3 -5.6){$\left\{ \vbox{\vskip2.7cm}\right. $}
}
\caption{\unhbox\cycbox}
\label{fig:cycl}
\end{figure}

In this section we enumerate {\em cyclically symmetric lozenge
tilings} of the cored hexagon $C_a(m):=C_{a,a,a}(m)$
with respect to certain weights.
By a cyclically symmetric lozenge tiling we mean a lozenge tiling
which is invariant under rotation by $120^\circ$. See
Figure~\ref{fig:cycl} for an example. (At this point, all shadings,
thick and dotted lines should be ignored.)
The unweighted
enumeration of these lozenge tilings was given earlier in
\cite[Theorem~3.2 and Corollary~3.3]{CiKrAB}. 
We restate the result below. 
Let $L_c(R)$ denote the number of cyclically symmetric lozenge tilings 
of region $R$.

\begin{theorem} \label{thm:cycl}
Let $a$ be a nonnegative integer.
The number $L_c(C_a(m))$ of cyclically symmetric lozenge tilings of a 
hexagon with side lengths $a,a+m,a,a+m,a,a+m$, with an equilateral triangle of
side length $m$ removed from the center, equals the right-hand
side in \eqref{eq:Andrews}. \qed
\end{theorem}

Let us now associate certain weights to each such lozenge tiling $T$.
These weights depend again on the number $n(T)$ of edges of 
lozenges of the tiling $T$ which are incident to the extension to the
right
of the bottommost side of the core. (Since we are now dealing with
cyclically symmetric tilings, it does, in fact, not matter which side is
considered, and the weighted count is not even affected by the choice 
of direction.) 
In the following three theorems,
each lozenge tiling $T$ is assigned the weight $\om^{n(T)}$, where
$\om$ is some fixed $6$th root of unity. 
Denote by $L_c^\om(R)$ the corresponding
weighted count of cyclically symmetric lozenge tilings of region $R$.

\begin{theorem} \label{thm:cycl-1} 
Let $a\geq0$ and $m\geq0$ be integers.
Then the weighted count $L_c^\om(C_a(m)):=\sum \om^{n(T)}$, summed over all
cyclically symmetric lozenge tilings $T$ of a 
hexagon with side lengths $a,a+m,a,a+m,a,a+m$, with an equilateral triangle of
side length $m$ removed from the center, equals the right-hand
side in \eqref{eq:Zare1} if $\om=-1$,
it equals the right-hand side in \eqref{eq:om3} if
$\om$ is a primitive third root of unity, and it
equals the right-hand side in \eqref{eq:om6} if
$\om$ is a primitive sixth root of unity.
\end{theorem}

If we specialize these results to $m=0$, i.e., to the case where
there exists {\em no} core, we obtain enumeration results for {\em
cyclically symmetric plane partitions}. Before we state these, let
us briefly recall the relevant notions from plane partition theory
(cf\@. e.g\@. \cite{StanAI} or \cite[Sec.~1]{StemAL}). There are
(at least) three possible equivalent ways to define plane partitions.
Out of the three possibilities,
in this paper, we choose to define a plane partition $\pi$ as 
a subset of the
three-dimensional integer lattice $\mathbb Z_+^3$
(where $\mathbb Z_+$ denotes the set of positive integers), with the property
that if $(i_1,j_1,k_1)$ is an element of $\pi$, then all points
$(i_2,j_2,k_2)$ with $1\le i_2\le i_1$, 
$1\le j_2\le j_1$, 
$1\le k_2\le k_1$ also belong to $\pi$. (In the language of 
partially ordered sets,
$\pi$ is an order ideal of $\mathbb Z_+^3$.) A plane partition $\pi$
is called {\em cyclically symmetric} if for every $(i,j,k)$ in $\pi$ the
point $(j,k,i)$ which results by a
cyclic permutation of coordinates is in $\pi$ as well.

\newbox\symppbox
\setbox\symppbox\hbox{\small A cyclically symmetric plane partition.}
\begin{figure}
\centertexdraw{
\drawdim truecm \linewd.02
\RhombusAsh \RhombusAsh \RhombusAsh \RhombusBsh \RhombusBsh \RhombusAsh \RhombusBsh
\RhombusBsh
\move(-.866025 -.5)
\RhombusAsh \RhombusBsh \RhombusAsh \RhombusBsh \RhombusAsh \RhombusBsh \RhombusAsh \RhombusBsh
\move(-1.732 -1) 
\RhombusBsh \RhombusAsh \RhombusBsh \RhombusAsh \RhombusBsh \RhombusAsh \RhombusAsh \RhombusBsh
\move(-2.598 -1.5)
\RhombusBsh \RhombusAsh \RhombusBsh \RhombusBsh \RhombusAsh \RhombusBsh \RhombusAsh \RhombusAsh
\move(1.732 -1) \RhombusCsh
\move(3.464 -1) \RhombusCsh
\move(3.464 -2) \RhombusCsh
\move(-2.598 -2.5) \RhombusCsh
\move(-2.598 -3.5) \RhombusCsh
\move(-2.598 -4.5) \RhombusCsh \RhombusCsh
\move(0 -5) \RhombusCsh \RhombusCsh
\linewd.1 
\move(.866025 -7.5)
%\rlvec(0 -4) 
\rlvec(3.464 2) \rlvec(0 4)
\rlvec(-3.464 2) \rlvec(-3.464 -2) 
\rlvec(0 -4) \rlvec(3.464 -2) 
%\move(.866025 -3.5)
%\rlvec(3.464 2)
%\move(.866025 -3.5)
%\rlvec(-3.464 2)
}
\caption{\unhbox\symppbox}
\label{fi:sympp} 
\end{figure}

Often, a plane partition is viewed as the corresponding pile of unit
cubes which results when replacing each point $(i,j,k)$ of the plane
partition by the unit cube with center $(i,j,k)$. A three-dimensional
picture of a plane partition, viewed as pile of unit cubes, is shown
in Figure~\ref{fi:sympp} 
(in fact, this example is cyclically symmetric).
As we already mentioned in the Introduction, plane partitions
contained in an $a\times b\times c$ box (i.e., plane partitions $\pi$ with
the property that every $(i,j,k)\in \pi$ satisfies $1\le i\le a$,
$1\le j\le b$, $1\le k\le c$) are in bijection with lozenge tilings
of a hexagon with side lengths $a,b,c,a,b,c$ (see
\cite{David-Tomei}). This bijection 
can be visualised easily on the example in 
Figure~\ref{fi:sympp}. Clearly, under this bijection, cyclically symmetric plane
partitions contained in an $a\times a\times a$ box correspond to 
cyclically symmetric lozenge tilings
of a hexagon with all sides of length $a$. Thus,
Theorem~\ref{thm:cycl-1} with $m=0$ yields
results about certain weighted counts of cyclically symmetric plane 
partitions. We just have to figure out how the weights $\om^{n(T)}$
for lozenge tilings $T$ translate to the plane partition language.

Let $\pi_T$ be the plane partition that corresponds to the lozenge
tiling $T$ under this bijection.
Denote by $m_1(\pi_T)$ the number of elements of the form $(i,i,i)$ in $\pi_T$.
Then there are precisely $m_1(\pi_T)$ unit cubes on the main diagonal of the
pile of unit cubes representing $\pi_T$. Let $v$ be the vertex farthest from 
the origin of the last such unit cube (in the planar rendering of
$\pi_T$ --- for our example, Figure~\ref{fi:sympp} --- 
$v$ is the center of the 
hexagon). A ray through $v$ approaching orthogonally any of the coordinate planes
will cut through precisely $m_1(\pi_T)$ layers of unit thickness. Since
each such cut corresponds to a lozenge side contained in the ray, we see 
that $m_1(\pi_T)$ is precisely the statistic $n(T)$.
We therefore obtain the following corollary of Theorem~\ref{thm:cycl-1}.

\begin{cor} \label{cor:cycl-1}
Let $a$ be a nonnegative integer.
Then the weighted count $\sum \om^{m_1(\pi)}$, summed over all
cyclically symmetric plane partitions $\pi$ contained in an $a\times
a\times a$ box, equals the right-hand
side in \eqref{eq:Zare1} with $m=0$ if $\om=-1$, it
equals the right-hand side in \eqref{eq:om3} with $m=0$ if
$\om$ is a primitive third root of unity, and it 
equals the right-hand side in \eqref{eq:om6} with $m=0$
if $\om$ is a primitive sixth root of unity. \qed
\end{cor}

Weighted enumerations of this sort have been considered
earlier. In fact, the result for $\om=-1$
of Corollary~\ref{cor:cycl-1} had been
conjectured by Stembridge \cite[Case~9 on p.~6]{StemAZ}, and proved
for the first time by Kuperberg \cite[last displayed equation on
p.~27]{KupeAG}. Thus, the $(-1)$-result of
Theorem~\ref{thm:cycl-1} is a generalization of Kuperberg's result.
There are many more conjectures on $(-1)$-enumerations of
cyclically symmetric plane partitions in \cite{StemAZ}. One of these,
the Conjecture on p.~7 of \cite[Case~10]{StemAZ}, asks for the
weighted count $\sum (-1)^{m_6(\pi)}$ of cyclically symmetric plane
partitions in which the statistic $m_6(\pi)$ is defined as the number
of orbits (under cyclic rotation) $\{(i,j,k),(j,k,i),(k,i,j)\}$
of elements of $\pi$ with coordinates that are not all equal.

%\newbox\cycbox
%\setbox\cycbox\hbox{}
\begin{figure}
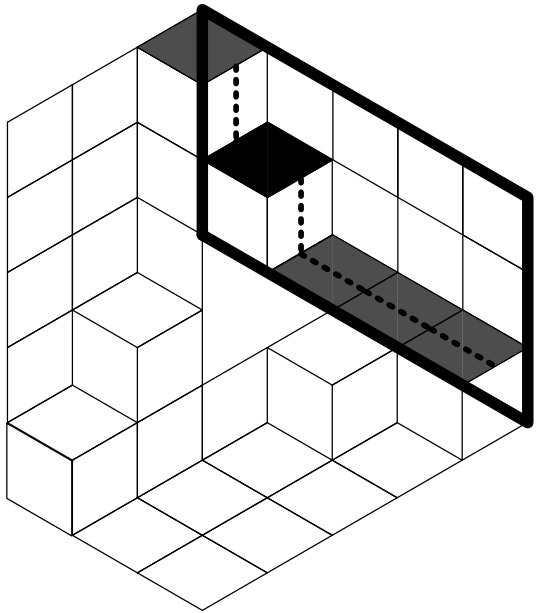
 
\centertexdraw{
\drawdim truecm \linewd.02
\RhombusA \RhombusB \RhombusA\RhombusB \RhombusA \RhombusA\RhombusA \RhombusB
\move(-.866025 -.5)
\RhombusB\RhombusB \RhombusB\RhombusA\RhombusB\RhombusB \RhombusA\RhombusA 
\move(-1.7305 -1)
\RhombusB\RhombusB\RhombusB\RhombusB \RhombusA\RhombusB\RhombusA \RhombusA 
\move(.866025 -4.5)
\RhombusB \RhombusA\RhombusA 
\move(1.7305 -4)
\RhombusA \RhombusB \RhombusA
\move(0 0)  \RhombusC 
\move(0 -1) \RhombusC 
\move(0.866025 -1.5) \RhombusC 
\move(3.464 -2) \rlvec(0 -1)

\move(2.598 -.5) \RhombusC \RhombusC \RhombusC
\move(4.3301 -2.5) \RhombusC
\move(-1.73205 -5) \RhombusC
\move(0.866025 -2.5) \rlvec(0 -1)

\linewd.02
\move(0.866025 0.5) \rlvec(-.866025 -.5) \rlvec(.866025 -.5)
\lfill f:0.3
\rlvec(.866025 .5) \rlvec(-.866025 .5)
\lfill f:0.3

\move(1.73205 -1)  \rlvec(-.866025 -.5) \rlvec(.866025 -.5)
\lfill f:0
\rlvec(.866025 .5) \rlvec(-.866025 .5)
\lfill f:0

\move(2.598 -2.5)  \rlvec(-.866025 -.5) \rlvec(.866025 -.5)
\lfill f:.3
\rlvec(.866025 .5) \rlvec(-.866025 .5)
\lfill f:.3

\move(3.464 -3)  \rlvec(-.866025 -.5) \rlvec(.866025 -.5)
\lfill f:.3
\rlvec(.866025 .5) \rlvec(-.866025 .5)
\lfill f:.3

\move(4.3301 -3.5)  \rlvec(-.866025 -.5) \rlvec(.866025 -.5)
\lfill f:.3
\rlvec(.866025 .5) \rlvec(-.866025 .5)
\lfill f:.3

\linewd.08
\move(1.309037 -.25)
\vdSchritt \hdSchritt \vdSchritt \hdSchritt \hdSchritt \hdSchritt

\linewd.15
\move(.866025 -2.5) \rlvec(0 3)
\rlvec(4.3301 -2.5)  \rlvec(0 -3)
\rlvec(-4.3301 2.5)

}
\caption{The statistic $n_6$ for this tiling is $n_6(T_0)=3$}
\label{fig:n6}
\end{figure}

We prove this conjecture of Stembridge in Theorem~\ref{thm:Case10} below.
In fact, in Theorem~\ref{thm:Case10} we prove a result for 
cyclically symmetric {\em lozenge tilings} of cored
hexagons.
In this result, a cyclically symmetric lozenge tiling $T$ is given a
weight $(-1)^{n_6(T)}$, with the statistic $n_6(T)$ to be described below.
It is defined in a way so that in the case when there is no
core present (i.e., $m=0$) it reduces to $m_6(\pi_T)$,
where again $\pi_T$ denotes the plane partition corresponding
to $T$. 

Let $T$ be a fixed cyclically symmetric lozenge tiling of the cored
hexagon $C_a(m)$ (see
Figure~\ref{fig:n6} for an example with $a=3$ and $m=2$; at this point, all 
thick lines and shadings should be ignored).
We consider
the horizontal lozenges which are at least partially contained in the
top-right fundamental region. (In Figure~\ref{fig:n6} the top-right
fundamental region is framed. The horizontal lozenges which are at
least partially contained in that region are the grey and black
lozenges.) The statistic
$n_6(T)$ is by definition the sum of the vertical distances
between these horizontal lozenges and the lower border of the
fundamental region.
(Thus, for the lozenge tiling $T_0$ in Figure~\ref{fig:n6} we have,
considering the horizontal lozenges in the order from left to right,
$n_6(T_0)=2+1+0+0+0=3$.) 

Suppose now that $m=0$, and view the tiling $T$ as a plane partition
$\pi_T$. The fundamental region of $T$ used in our definition of the
statistic $n_6$ corresponds to a fundamental region of $\pi_T$ with
the main diagonal removed. Since the distances we add up in our definition
of $n_6(T)$ are precisely the heights of the vertical columns of unit
cubes in this fundamental region, we obtain that $n_6(T)$ is equal to
the number of unit cubes contained in it, which is clearly just the number
of orbits of cubes off the main diagonal. This verifies our claim that 
$n_6(T)=m_6(\pi_T)$. 

The weight which is assigned to a tiling $T$
in the theorem below is $(-1)^{n_6(T)}$. An equivalent way to define
this weight is to say that it is the product of the weights of all
lozenges which are, at least partially, contained in the top-right
fundamental region, where the weight of a horizontal lozenge with odd
distance from the lower border of the region is $-1$, the
weight of all other lozenges being 1. (In Figure~\ref{fig:n6} the
black lozenge has weight $-1$, all other lozenges have weight $1$.)
Yet another way to obtain this
weight is through the perfect matchings point of view of lozenge
tilings, elaborated for example in \cite{KupeAA,KupeAG}. In this
setup, the cyclically symmetric lozenge tilings that we consider here
correspond bijectively to perfect matchings in a certain hexagonal
graph (basically, the dual graph of a fundamental region of the
cored hexagon). Assignment of weights to the edges of this graph
so that each face has ``curvature" $-1$ (see \cite[Sec.~II]{KupeAG}) 
generates again (up to a multiplicative constant) 
the above weight for lozenge tilings.

Denote by $L_o^{-1}(R)$ (where the index letter stands for ``orbits'') the 
weighted count of lozenge tilings of region $R$  under the above-defined 
weight.

\begin{theorem} \label{thm:Case10}
Let $a$ and $m$ be nonnegative integers.
Let $R_1(a,m)$ denote the right-hand side of \eqref{eq:Andrews},
let $R_2(a,m)$ denote the right-hand side of \eqref{eq:Zare1}, 
and let $R_3(a,m)$ denote the
right-hand side of \eqref{eq:om6}.
Then the weighted count $\sum (-1)^{n_6(T)}$, summed over all
cyclically symmetric lozenge tilings $T$ of a 
hexagon with side lengths $a,a+m,a,a+m,a,a+m$ with an equilateral triangle of
side length $m$ removed from the center, is given by
\begin{equation} \label{eq:Case10}
L_o^{-1}(C_a(m))=\begin{cases}
|R_3(\tfrac a 2,\tfrac m 2)|^2 \quad &\text {if\/ $a$ is even and $m$ is even,}\\
R_1(\tfrac {a+1}2,\tfrac m 2 -1) R_1(\tfrac {a-1}2,\tfrac m 2 +1) \quad
&\text {if\/ $a$ is odd and $m$ is even,}\\
R_1(\tfrac a 2 ,\tfrac {m-1}2) R_2( \tfrac a 2,\tfrac {m+1}2) \quad &\text
{if\/ $a$ is even and $m$ is odd,}\\
R_1(\tfrac{a+1}2,\tfrac{m-1}2) R_2(\tfrac{a-1}2,\tfrac{m+1}2) \quad &\text
{if\/ $a$ is odd and $m$ is odd.}
\end{cases}
\end{equation}
\end{theorem}

As we show in Section~\ref{sec:outline}, all the above results in
the current section follow from
evaluations of the determinant \eqref{eq:omega} for $\om$ equal to
$1$, to $-1$, to a primitive third root of unity, and to a 
primitive sixth root of unity,
respectively. The corresponding evaluations read as follows, the
evaluation for $\om=1$, given in Theorem~\ref{thm:Andrews} below, 
being due to Andrews \cite[Theorem~8]{AndrAN}. 

\begin{theorem} \label{thm:Andrews}
For any nonnegative integer $a$,
\begin{multline} \label{eq:Andrews}
\det_{0\le i,j\le a-1}\left(\delta_{ij}+\binom {m+i+j}{j}\right)\\
=\cases   \displaystyle
2^{\cl{a/2}}\prod _{i=1} ^{a-2}\(\frac {m} {2}+\cl{i/2}+1\)_{\fl{(i+3)/4}}\\
\displaystyle\times
\frac {\prod _{i=1} ^{a/2}\(\frac {m} {2}+\frac {3a}2-\cl{\frac {3i}2}+\frac 32\)_{\cl{i/2}-1}
                \(\frac {m} {2}+\frac {3a}2-\cl{\frac {3i}2}+\frac 32\)_{\cl{i/2}}}
        {\prod _{i=1} ^{a/2-1}(2i-1)!!\,(2i+1)!!}&\text {if\/ $a$ is
even,}\\
\displaystyle
2^{\cl{a/2}}\prod _{i=1} ^{a-2}\(\frac {m} {2}+\cl{i/2}+1\)_{\cl{(i+3)/4}}\\
\displaystyle\times
\frac {\prod _{i=1} ^{(a-1)/2}\(\frac {m} {2}+\frac {3a}2-\cl{\frac {3i-1}2}+1\)_{\cl{(i-1)/2}}
           \(\frac {m} {2}+\frac {3a}2-\cl{\frac {3i}2}\)_{\cl{i/2}}}
        {\prod _{i=1} ^{(a-1)/2}(2i-1)!!^2}&\text {if\/ $a$ is odd,}
\endcases
\end{multline}
where $(\alpha)_k$ is the standard notation for {\em shifted factorials},
$(\alpha)_k:=\alpha(\alpha+1)\cdots(\alpha+k-1)$, $k\ge1$, and $(\alpha)_0:=1$.
\quad \quad \qed
\end{theorem}

\begin{theorem} \label{thm:Zare1}
For nonnegative integers $a$,
\begin{multline} \label{eq:Zare1}
\det_{0\le i,j\le a-1}\left(-\delta_{ij}+\binom {m+i+j}{j}\right)
\\=
\cases 0,&\text {if\/ $a$ is odd,}\\
(-1)^{a/2}\prod _{i=0} ^{a/2-1}\frac {i!^2\,(\frac {m} {2}+i)!^2\,
(\frac {m} {2}+3i+1)!^2\,(m+3i+1)!^2} {(2i)!\,(2i+1)!\,(\frac {m} {2}+2i)!^2\,
(\frac {m} {2}+2i+1)!^2\,(m+2i)!\,(m+2i+1)!},
&\text {if\/ $a$ is even.}\endcases
\end{multline}
\end{theorem}
The proof of this theorem is given in Section~\ref{sec:Zare1}.

\begin{theorem} \label{thm:om3}
Let $\om$ be a primitive third root of unity. Then 
\begin{multline} \label{eq:om3}
\det_{0\le i,j \le a-1}\(\om \delta_{ij}+ \binom{m+i+j}{j}\)=
\frac {(1+\om)^a\, 2^{\fl{{a} /{2}}}} 
{\prod_{i=1} ^{\fl{ {a} /{2}}}{(2i-1)!!}
\prod_{i=1} ^{\fl{ {(a-1)} /{2}}}{(2i-1)!!}}\\
\times\prod _{i\ge 0} ^{}{\(\tfrac {m} {2}+3i+1\)_{\fl{ {(a-4i)}/{2}}}
\(\tfrac {m} {2}+3i+3\)_{\fl{ {(a-4i-3)}/{2}}}}\\
\cdot
\(\tfrac {m} {2}+a-i+\tfrac {1} {2}\)_{\fl{{(a-4i-1)}/{2}}}
\(\tfrac {m} {2}+a-i-\tfrac {1} {2}\)_{\fl{{(a-4i-2)}/{2}}},
\end{multline}
where, in abuse of notation, by $\fl{\alpha}$ we mean the usual
floor function if $\alpha \ge 0$, however, if $\alpha <0$ then
$\fl{\alpha}$ must be read as $0$, so that the product over $i\ge 0$ is
indeed a finite product. 
\end{theorem}
The proof of this theorem is given in Section~\ref{corsec}.

\begin{theorem} \label{thm:om6}
Let $\om$ be a primitive sixth root of unity. 
Then 
\begin{multline} \label{eq:om6}
\det_{0\le i,j \le a-1}
{\(\om \delta_{ij}+\binom{m+i+j}{j}\)}=
\frac {(1+\om)^a \(\frac {2} {3}\)^{\fl{ {a}/ {2}}}}
{\prod_{i=1}^{\fl{ {a}/ {2}}}{(2i-1)!!}\prod _{i=1} ^{\fl{
{(a-1)}/ {2}}}{(2i-1)!!}}\\
\times\prod _{i\ge 0} ^{}{\(\tfrac {m} {2} +3i+\tfrac {3} {2}\)_{\fl{
{(a-4i-1)}/ {2}}}
\(\tfrac {m} {2}+3i+\tfrac {5} {2}\)_{\fl{ {(a-4i-2)}/ {2}}}}\\
\cdot (\tfrac {m} {2}+a-i)_{\fl{ {(a-4i)}/ {2}}}\,(\tfrac {m} {2}+a-i)_{\fl{
{(a-4i-3)}/{2}}},
\end{multline}
where again, in abuse of notation, by $\fl{\alpha}$ we mean the usual
floor function if $\alpha \ge 0$, however, if $\alpha <0$ then
$\fl{\alpha}$ must be read as $0$, so that the product over $i\ge 0$ is
indeed a finite product. 
\end{theorem}
The proof of this theorem is given in Section~\ref{omsec}.
\end{section}

\begin{section}{Outline of the proofs of
Theorems~\ref{enum}--\ref{thm:Case10}}
\label{sec:outline}
In this section, we give outlines of the proofs of our
enumeration results stated in the Introduction and in Sections~\ref{sec:-1}
and \ref{sec:cyclically}. We fill in the details of these proofs in
later sections.

\begin{proof}[Proof of Theorem~\ref{enum}] 
There is a standard bijection between lozenge tilings and families of
nonintersecting lattice paths. This bijection is explained in
Section~\ref{sec:nonint} (see in particular Figure~\ref{tiling}).
Thus, the problem of enumerating lozenge tilings is converted to
the problem of counting certain families of nonintersecting lattice
paths. 
By the Lindstr\"om--Gessel--Viennot theorem (stated in Lemma~\ref{gv}), 
the number of
such families of paths can be expressed as a determinant 
(see Lemma~\ref{lem:nonint}).
Thus, in principle, we would be done once we evaluate this determinant, 
given in \eqref{mat1}.
However, Lemma~\ref{lem:nonint} applies only if the size $m$ 
of the core is even. 
We show, in Section~\ref{sec:poly}, that it suffices to address this
case, by 
proving that the number of lozenge tilings that we are interested in is  
a polynomial in $m$. 
The evaluation of the determinant \eqref{mat1} for even $m$
is carried out in Section~\ref{sec:det} (see \eqref{det} and
Lemmas~\ref{lem:amgerade} and \ref{lem:amung}).
\end{proof}

\begin{proof}[Proof of Theorem~\ref{-enum}] 
The first steps are identical with those in the preceding proof:
the lozenge tilings are converted into nonintersecting lattice paths,
in the way that is described in Section~\ref{sec:nonint}. Therefore,
Lemma~\ref{gv} yields a determinant for the $(-1)$-enumeration 
that we are interested in. Unlike in the previous proof, this provides
a determinant for our weighted count 
only if the size $m$ of the core is odd (see
Lemma~\ref{lem:nonint}). Again, the considerations in
Section~\ref{sec:poly} show that this number 
is a polynomial in $m$, so it suffices to evaluate the
determinant \eqref{mat1} for odd $m$. This is done 
in Section~\ref{sec:det} (see \eqref{det} and
Lemmas~\ref{lem:-agerade} and \ref{lem:-aung}).
\end{proof}

\begin{proof}[Proof of Theorem~\ref{shifted}] 
Again, we use the strategy from the proof of Theorem~\ref{enum}.
We convert the lozenge tilings into families of nonintersecting 
lattice paths as described in Section~\ref{sec:nonint}. The starting
and ending points are slightly different from the ones used before. They
are given in \eqref{shcoord}. 
Lemma~\ref{gv} yields a determinant for the number we are interested
in for even $m$ (see Lemma~\ref{lem:nonint2}). The considerations of
Section~\ref{sec:poly} still apply, so the number of lozenge tilings
is a polynomial in $m$ and it suffices to evaluate the determinant
\eqref{mat2} for even $m$.
This is accomplished in Section~\ref{sec:det2} (see
\eqref{det2} and
Lemmas~\ref{lem:shamgerade} and \ref{lem:shaome}).
\end{proof}

\begin{proof}[Proof of Theorem~\ref{-shifted}] 
We proceed analogously to the proof of Theorem~\ref{shifted}.
The lozenge tilings are converted into nonintersecting lattice paths,
in the way that is described in Section~\ref{sec:nonint}. Therefore,
Lemma~\ref{gv}
yields a determinant for the $(-1)$-enumeration in the case of odd
$m$ (see
Lemma~\ref{lem:nonint2}). Again, the considerations in
Section~\ref{sec:poly} show that this number 
is a polynomial in $m$, so that it suffices to evaluate the
determinant \eqref{mat2} for odd $m$. This is worked out
in Section~\ref{sec:det2} (see \eqref{det2} and
Lemmas~\ref{lem:shaemo} and \ref{lem:shaomo}).
\end{proof}

\begin{proof}[Proof of Theorem~\ref{thm:cycl-1}] 
We follow the arguments of the proof of Theorem~\ref{thm:cycl}, as
given in \cite[Lemma~3.1]{CiKrAB}. Suppose we are given a cyclically
symmetric lozenge tiling $T$ of our cored hexagon $C_a(m)$. 
It is completely determined by its restriction to a
fundamental region, the lower-left fundamental region, say. (In the
example in Figure~\ref{fig:cycl}, the lower-left fundamental region
is framed.) Some of the lozenges are cut in two by the borders of
the fundamental region. (In Figure~\ref{fig:cycl} these are the
shaded lozenges.) 
We draw lattice paths which connect these ``cut"
lozenges, by ``following" along the other lozenges, as is indicated in
Figure~\ref{fig:cycl} by the dashed lines. To be precise, in each
lozenge in the interior of the fundamental region, 
we connect the midpoints of the sides that run up-diagonal, in case the
lozenge possesses such sides.
Clearly, these paths are {\em nonintersecting}, by which we mean that no
two paths have a common vertex. Since
they determine completely the cyclically symmetric lozenge tiling, we
may as well count all these families of nonintersecting lattice paths, with
respect to the corresponding weight. In fact, as is easy to see,
because of the cyclic symmetry, the statistic $n(T)$ is exactly equal
to $a$ minus the number of paths. If we fix the ``cut" lozenges, say in
positions $i_1,i_2,\dots,i_k$ (counted from inside out,
beginning with $0$; thus, in Figure~\ref{fig:cycl}, the ``cut" lozenges have
positions $0$ and $2$), then, according to Lemma~\ref{gv}, 
the number of families of
nonintersecting lattice paths connecting the fixed ``cut" lozenges is
given by the corresponding 
Lindstr\"om--Gessel--Viennot determinant (the left-hand side of
\eqref{eq:gv}). This determinant turns
out to be the minor of $\(\binom {m+i+j}j\)_{0\le i,j\le a-1}$
consisting of rows and columns with indices
$i_1,i_2,\dots,i_k$. This number must
be multiplied by the common weight $\om^{a-k}$ of these families of
nonintersecting lattice paths. Therefore, in order to obtain the
total weighted count that we are interested in, we have to sum all
these quantities, i.e., take the sum of
$$\textstyle\Big(\text {$(i_1,i_2,\dots,i_k)$-principal minor of }
\(\binom {m+i+j}j\)_{0\le i,j\le a-1}\Big)\times \om^{a-k}$$
over all $k=0,1,\dots,a$ and $0\le i_1<i_2<\dots<i_k\le a-1$.
Clearly, this sum is exactly equal to 
$\det_{0\le i,j\le a-1}\(\om\delta_{ij}+\binom{m+i+j}j\)$, 
which equals the left-hand
side of \eqref{eq:Zare1} if $\om=-1$, the left-hand side of
\eqref{eq:om3} if $\om$ is a primitive third root of unity, and 
the left-hand side of \eqref{eq:om6} if
$\om$ a primitive sixth root of unity. The respective
right-hand sides provide therefore the solution to our enumeration
problem.
\end{proof}

\begin{proof}[Proof of Theorem~\ref{thm:Case10}] 
We adapt the arguments used in the proof of Theorem~\ref{thm:cycl-1}.
(Clearly, here we want to count the same objects, but with respect
to a different weight.) So, again, we draw paths that connect the
lozenges which are cut in two by the borders of the fundamental
region. This time, we choose the top-right region as
the fundamental region. Figure~\ref{fig:n6} shows an example. There,
the top-right fundamental region is framed. As in Figure~\ref{fig:cycl},
paths are indicated by dashed
lines. (In the example in Figure~\ref{fig:n6} there is just one path.)
If we slightly distort 
the underlying lattice, we get orthogonal paths with
positive horizontal and negative vertical steps.
Figure~\ref{fig:n6b} shows the orthogonal path corresponding to the
path in Figure~\ref{fig:n6}. The manner in which we have chosen the coordinate
system ensures that possible starting points of paths are the points
$(0,j)$, $0\le j\le a-1$, and possible ending points are the points
$(m+i,0)$, $0\le i\le a-1$.

\begin{figure}
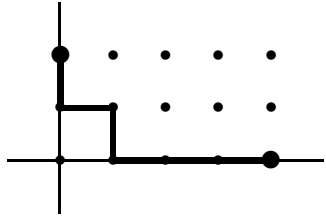

$$
\PfadDicke{2pt}
\Einheit.7cm
\Gitter(5,3)(0,0)
\Koordinatenachsen(5,3)(0,0)
\Pfad(0,1),2\endPfad
\Pfad(0,1),1\endPfad
\Pfad(1,0),2\endPfad
\Pfad(1,0),111\endPfad
\DickPunkt(0,2)
\DickPunkt(4,0)
\hskip2.5cm
$$
\caption{The orthogonal path corresponding to Figure~\ref{fig:n6}}
\label{fig:n6b}
\end{figure}

Now, as before, we fix the positions of the ``cut" lozenges. Then a {\em
weighted version} of the Lindstr\"om--Gessel--Viennot theorem 
(see \cite[Lemma~1]{LindAA} or \cite[Cor.~2]{gv})
can be used to express the weighted count of the corresponding
families of nonintersecting lattice paths in form of a determinant. 
In fact, this weighted
version just says that Lemma~\ref{gv} remains true when the number 
$\P(A\to E)$ of paths from $A$ to $E$ is replaced everywhere by the
{\em weighted count\/} $\sum _{P} ^{}w(P)$ of all paths $P$ 
from $A$ to $E$, where 
$w$ is some weight function on the edges of the square lattice and the weight 
$w(P)$ of a path is the product of the weights of its steps.
Thus, if we repeat the subsequent arguments in the proof of
Theorem~\ref{thm:cycl-1}, then we obtain the determinant
\begin{equation} \label{eq:Pfade}
\det_{0\le i,j\le a-1}\bigg(\delta_{ij}+\sum _{P:(0,j)\to(m+i,0)} ^{}
w(P)\bigg) 
\end{equation}
for the weighted count 
of our families of nonintersecting lattice paths. 

We now choose the weight function $w$ so that the weight of the family
of nonintersecting lattice paths corresponding to a tiling $T$ is equal 
to $(-1)^{n_6(T)}$. To do this, it will be convenient to stick on an extra
initial horizontal step at the beginning of each path, so that now it starts
on the line $x=-1$. Weight the vertical steps on this line by 0, all the
remaining vertical steps by 1, and weight horizontal steps at height $j$ by
$(-1)^j$. Since the height of a horizontal step is equal to the distance
of the corresponding horizontal lozenge to our reference line in the tiling, 
the weight of a family 
$(P_1,P_2,\dots)$
of nonintersecting lattice paths is equal to $(-1)^{n_6(T)}$,
where $T$ is the corresponding tiling. 
On the other hand, it is clearly equal to
$(-1)^{A(P_1)+A(P_2)+\cdots}$, where $A(P)$ denotes the area between
a path $P$ and the $x$-axis.

To find an expression for the entries of the Lindstr\"om--Gessel--Viennot matrix
we use the well-known fact (see \cite[Prop.~1.3.19]{StanAP}) that the
weighted count $\sum _{} ^{}q^{A(P)}$, summed over all lattice paths
$P$ from $(0,c)$ to $(d,0)$, is equal to $\left[\begin{smallmatrix} 
c+d\\c\end{smallmatrix}\right]_q$, where 
$\left[\smallmatrix n\\k\endsmallmatrix\right]_q$ is the standard
$q$-binomial coefficient,
$$\begin{bmatrix} n\\k\end{bmatrix}_q:=
\frac {(1-q^n)(1-q^{n-1})\cdots(1-q^{n-k+1})}
{(1-q^k)(1-q^{k-1})\cdots(1-q)}.$$
Thus, the determinant \eqref{eq:Pfade} becomes
(see also \cite[Lemma~4]{StemAZ})
\begin{equation} \label{eq:n6-det}
\det_{0\le i,j\le a-1}\(\delta_{ij}+
(-1)^j\left[\matrix m+i+j\\j\endmatrix\right]_{-1}\right).
\end{equation}
{}From the $q$-binomial theorem (see \cite[(3.3.6)]{AndrAF}),
$$(1+z)(1+qz)\cdots(1+q^{n-1}z)=\sum _{k=0} ^{n}q^{\binom k2}\bmatrix
n\\k\endbmatrix_q z^k,$$
it is straightforward to extract that
\begin{equation} \label{eq:-qbinom}
\begin{bmatrix} n\\k\end{bmatrix}_{-1}=\cases 0&\text {if $n$ is even and
$k$ is odd,}\\
\displaystyle\binom {\fl{n/2}} {\fl{k/2}}&\text {otherwise.}\endcases
\end{equation}

We have to compute the determinant \eqref{eq:n6-det}. Let us denote
it by $D_0$. We have to distinguish between four cases, depending on the
parities of $m$ and $a$.

First, let $m$ be even. 
We reorder rows and columns
simultaneously, so that the even-numbered rows and columns 
come before the odd-numbered, respectively.
If $a$ is even, then we obtain for $D_0$ the block determinant
$$\det\pmatrix I(\frac {a} {2})+B(\frac {a} {2},\frac {m} {2})&-B(\frac {a} {2},\frac {m} {2})\\
B(\frac {a} {2},\frac {m} {2})&I(\tfrac a2)\endpmatrix,$$
where $I(N)$ is the $N\times N$ identity matrix
and $B(N,m)$ is the $N \times N$ matrix\linebreak 
$\Big(\binom {m+i+j}j\Big)_{0\le i,j\le N-1}$. By a few simple manipulations, 
this determinant can be factored
into a product of two determinants,
\begin{align*} 
D_0&=\det\pmatrix I(\tfrac a2)+B(\frac {a} {2},\frac {m} {2})&-B(\frac {a}
{2},\frac {m} {2})\\B(\frac {a} {2},\frac {m} {2})&I(\tfrac a2)\endpmatrix\\
&=\det\pmatrix I(\tfrac a2)+B(\frac {a} {2},\frac {m} {2})&-B(\frac {a} {2},\frac {m} {2})\\
B(\frac {a} {2},\frac {m} {2})&I(\tfrac a2)\endpmatrix
\det\begin{pmatrix} I(\tfrac a2)&0\\-B(\frac {a} {2},\frac {m}
{2})&I(\tfrac a2)\end{pmatrix}
\\
&=\det\begin{pmatrix} I(\tfrac a2)+B(\frac {a} {2},\frac {m} {2})+B(\frac {a}
{2},\frac {m} {2})^2&-B(\frac {a} {2},\frac {m} {2})\\0&I(\tfrac a2)\end{pmatrix}\\
&=\det\left(I(\tfrac a2)+B(\tfrac {a} {2},\tfrac {m} {2})
+B(\tfrac {a} {2},\tfrac {m} {2})^2\right)\\
&=\det\left(\om I(\tfrac a2)+B(\tfrac {a} {2},\tfrac {m} {2})\right)
\det\left(\overline\om I(\tfrac a2)
+B(\tfrac {a} {2},\tfrac {m} {2})\right),
\end{align*}
where $\om$ is a primitive sixth root of unity, each of which can be
computed by application of
Theorem~\ref{thm:om3}. The result is the first expression in
\eqref{eq:Case10}.

On the other hand, if $a$ is odd, then analogous arguments yield
\begin{equation} \label{eq:Aquer}
D_0=\det\(I(\tfrac {a+1} {2})+B(\tfrac {a+1} {2},\tfrac {m} {2})+ 
B^{(c)}(\tfrac {a+1} {2},\tfrac {m} {2})\, 
B^{(r)}(\tfrac {a+1} {2},\tfrac {m} {2})\),
\end{equation}
where $B^{(c)}(\frac {a+1} {2},\frac {m} {2})$ is the 
$(\frac {a+1} {2})\times (\frac {a-1} {2})$
matrix which arises from $B(\frac {a+1} {2},\frac {m} {2})$ by 
deleting its last column,
while $B^{(r)}(\frac {a+1} {2},\frac {m} {2})$ is the 
$(\frac {a-1} {2})\times (\frac {a+1} {2})$
matrix which arises from $B$ by deleting its last row.

It is easy to check that
$$I(\tfrac {a+1} {2})+B(\tfrac {a+1} {2},\tfrac {m} {2})+ 
B^{(c)}(\tfrac {a+1} {2},\tfrac {m} {2})\, 
B^{(r)}(\tfrac {a+1} {2},\tfrac {m} {2})=(I(\tfrac {a+1} {2})+\overline B)
\left(I(\tfrac {a+1} {2})+B(\tfrac
{a+1} {2},\tfrac {m} {2}-1)\right),$$
where $\overline B$ is the $(\frac {a+1} {2})\times (\frac {a+1}
{2})$-matrix with $(i,j)$-entry $\binom{\frac m 2 +i+j-1}{j-1}$,
$0\le i,j\le (a-1)/2$. (So the first column of $\overline B$ is zero).
We expand $\det(I(\frac {a+1} {2})+\overline B)$ with respect to the first column and get
$\det(I(\frac {a-1} {2})+B(\frac {a-1} {2},\frac {m} {2}+1))$.

Therefore, in the case of even $m$ and odd $a$, we have
$$D_0=\det\left(I(\tfrac {a+1} {2})+B(\tfrac {a+1} {2},\tfrac m 2 -1)\right)
\det\left(I(\tfrac {a-1} {2})+B(\tfrac {a-1} {2},\tfrac m 2 + 1)\right).$$
Both determinants can be evaluated by means of
Theorem~\ref{thm:Andrews}. The result is the second expression in
\eqref{eq:Case10}.

Now let $m$ be odd. We proceed analogously. If $a$ is even, then
reordering rows and columns
according to the parity of the indices gives
\begin{align*}
D_0&=\det \begin{pmatrix}I(\tfrac {a} {2})+B(\frac {a} {2},\frac {m-1} {2}) &
0 \\
B(\frac {a} {2},\frac {m+1} {2}) & I(\tfrac {a} {2})-B(\frac {a} {2},\frac {m+1} {2})
\end{pmatrix}\\
&=\det\left(I(\tfrac {a} {2})+B(\tfrac {a} {2},\tfrac {m-1} {2})\right)
\det\left(I(\tfrac {a} {2})- B(\tfrac{a} {2},\tfrac {m+1} {2})\right).
\end{align*}
The first determinant is evaluated by means of
Theorem~\ref{thm:Andrews}, while the second is evaluated by means of
Theorem~\ref{thm:Zare1}. The result is the third expression in
\eqref{eq:Case10}.

Finally, if $a$ is odd we get
\begin{align*}
D_0&=\det \begin{pmatrix}I(\frac {a+1} {2})+B(\frac {a+1} {2},\frac {m-1} {2}) &
0 \\
B^{(r)}(\frac {a+1} {2},\frac {m+1} {2}) & I(\frac {a-1} {2})-B(\frac {a-1} {2},\frac {m+1} {2})
\end{pmatrix}\\
&=\det\left(I(\tfrac {a+1} {2})+B(\tfrac {a+1} {2},\tfrac {m-1} {2})\right)
\det\left(I(\tfrac {a-1} {2})-B(\tfrac
{a-1} {2},\tfrac {m+1} {2})\right).
\end{align*}
Again, the first determinant is evaluated by means of
Theorem~\ref{thm:Andrews}, while the second is evaluated by means of
Theorem~\ref{thm:Zare1}. The result is the fourth expression in
\eqref{eq:Case10}.
\end{proof}

\end{section}

\begin{section}{Lozenge tilings, nonintersecting lattice paths, and
determinants}
\label{sec:nonint}

The purpose of this section is to derive determinants for the
ordinary and $(-1)$-enumeration of
lozenge tilings of cored hexagons (see Lemmas~\ref{lem:nonint}
and \ref{lem:nonint2}). We find these determinants by first translating
the lozenge tilings to nonintersecting lattice paths, and
subsequently applying the Lindstr\"om--Gessel--Viennot theorem (stated here as 
Lemma~\ref{gv}). 

\medskip

{\em From lozenge tilings to nonintersecting lattice paths.}
There is a well-known translation of lozenge tilings
to families of nonintersecting lattice paths. 
We start with a lozenge
tiling of the cored hexagon (see
Figure~\ref{tiling}.a). We mark the midpoints of the edges along the
sides of length $a$ and $a+m$ and along the side of the triangle
which is parallel to them
(see Figure~\ref{tiling}.b). Now, in the same way
as in the proof of Theorem~\ref{thm:cycl-1} in the preceding section, we
connect these points by paths which ``follow" along the lozenges of the
tiling, as is illustrated in Figure~\ref{tiling}.b. 
Clearly, the resulting
paths are nonintersecting, 
i.e., no two paths have a common vertex. 
If we slightly distort
the underlying lattice, we get orthogonal paths with
positive horizontal and negative vertical steps (see
Figure~\ref{tiling}.c).
In the case that $a$, $b$ and $c$ have the same parity,
we can introduce a coordinate system in a way so that
the coordinates of the
starting points $A_i$ and end points $E_j$ are
\begin{subequations} \label{coord}
\begin{align}
A_i&=(i-1,c+m+i-1),& i&=1,2,\dots ,a, \\
A_i&=\(\frac {a+b}{2}+i-a-1,\frac {a+c} {2} +i-a-1\),&
i&=a+1,a+2,\dots,a+m,\\
E_j&=(b+j-1,j-1),& j&=1,2,\dots, a+m,
\end{align}
\end{subequations}
see Figure~\ref{tiling}.c.

\begin{figure} 
\centertexdraw{
\drawdim cm \setunitscale.7
\linewd.01
\rhombus \rhombus \rhombus \rhombus \rhombus \rhombus \rhombus
\ldreieck
\move(-.866025 -.5)
 \rhombus \rhombus \rhombus \rhombus \rhombus \rhombus \rhombus \rhombus 
\move(-1.73205 -1)
\rhombus \rhombus \rhombus \ldreieck \rmove(.866025 -.5)
\rhombus \rhombus \rhombus 
\move(-1.73205 -1)
\rdreieck \rhombus \rhombus \rhombus \ldreieck \rhombus \rhombus \rhombus 
\move(-1.73205 -2)
\rdreieck \rhombus \rhombus \rhombus \rhombus \rhombus \rhombus 
\move(-1.73205 -3)
\rdreieck \rhombus \rhombus \rhombus \rhombus \rhombus

\linewd.08
\move(0 0)
\RhombusA  \RhombusA \RhombusA \RhombusA \RhombusA \RhombusA \RhombusA
\RhombusB 
\move(-.866025 -.5)
\RhombusA \RhombusA \RhombusA \RhombusA \RhombusA \RhombusA 
\RhombusB \RhombusA 
\move(-1.73205 -1)
\RhombusA \RhombusA \RhombusB \RhombusB \RhombusA \RhombusB \RhombusA 
\RhombusA 
\move(2.598 -3.5)
\RhombusA \RhombusB \RhombusA 
\move(1.73205 -4)
\RhombusB \RhombusA \RhombusA
\move(-1.73205 -1)
 \RhombusC \RhombusC
\move(-1.73205 -2) \RhombusC
\move(-1.73205 -3)\RhombusC \RhombusC \RhombusC
\move(.866025 -1.5) \RhombusC
\move(.866025 -2.5) \RhombusC
\htext(-4 -8){\small 
a. A lozenge tiling of the cored hexagon in Figure~\ref{hex}}

\htext(-1 0.2){$a$}
\htext(4.4 -1){$b+m$}
\htext(6.9 -3.9){$\left. \vbox{\vskip.35cm} \right\} c$}
\htext(-3.8 -4){$c+m \left\{ \vbox{\vskip1.2cm} \right.$}
\htext(0 -6.2){$b$}
\htext(5 -6){$a+m$}
\htext(1.7 -4){$\left. \vbox{\vskip0.77cm} \right\}m$}
\rtext td:60 (4 -1.3) {$\left. \vbox{\vskip2.52cm} \right\} $}
\rtext td:-60 (-.8 0.2){$\left\{ \vbox{\vskip1.02cm} \right. $}
\rtext td:-60 (4.6 -5.3) {$\left. \vbox{\vskip1.75cm} \right\} $}
\rtext td:60 (0.3 -5.8){$\left\{ \vbox{\vskip1.89cm}\right. $}

\move(12 0)
\bsegment
\linewd.05
\move(0 0)
\RhombusA  \RhombusA \RhombusA \RhombusA \RhombusA \RhombusA \RhombusA
\RhombusB 
\move(-.866025 -.5)
\RhombusA \RhombusA \RhombusA \RhombusA \RhombusA \RhombusA 
\RhombusB \RhombusA 
\move(-1.73205 -1)
\RhombusA \RhombusA \RhombusB \RhombusB \RhombusA \RhombusB \RhombusA 
\RhombusA 
\move(2.598 -3.5)
\RhombusA \RhombusB \RhombusA 
\move(1.73205 -4)
\RhombusB \RhombusA \RhombusA
\move(-1.73205 -1)
 \RhombusC \RhombusC
\move(-1.73205 -2) \RhombusC
\move(-1.73205 -3)\RhombusC \RhombusC \RhombusC
\move(.866025 -1.5) \RhombusC
\move(.866025 -2.5) \RhombusC

\ringerl(.433012 .25) 
\hdSchritt \hdSchritt  \hdSchritt \hdSchritt \hdSchritt \hdSchritt \hdSchritt  \vdSchritt  
\ringerl(-.433012 -.25)
\hdSchritt \hdSchritt \hdSchritt \hdSchritt \hdSchritt \hdSchritt  \vdSchritt \hdSchritt  
\ringerl(-1.299037 -.75) \hdSchritt \hdSchritt 
\vdSchritt \vdSchritt  \hdSchritt  \vdSchritt  \hdSchritt \hdSchritt 
\ringerl(3.031 -3.25) \hdSchritt \vdSchritt \hdSchritt
\ringerl(2.165 -3.75) \vdSchritt \hdSchritt \hdSchritt
\ringerl(6.4952 -4.25)
\ringerl(5.6292 -4.75)
\ringerl(4.7632 -5.25) 
\ringerl(3.8971 -5.75)
\ringerl(3.031  -6.25)
\htext(-2 -8){\small b. The corresponding path family}
\esegment
\htext(4 -18){\small c. The path family made orthogonal}}
\vskip-7cm
$$
\Einheit=.7cm
\Gitter(11,7)(0,0)
\Koordinatenachsen(11,7)(0,0)
\Pfad(0,3),11\endPfad
\Pfad(2,1),22\endPfad
\Pfad(2,1),1\endPfad
\Pfad(3,0),2\endPfad
\Pfad(3,0),11\endPfad
\Pfad(1,4),111111\endPfad
\Pfad(7,3),2\endPfad
\Pfad(7,3),1\endPfad
\Pfad(2,5),1111111\endPfad
\Pfad(9,4),2\endPfad
\Pfad(4,1),2\endPfad
\Pfad(4,1),11\endPfad
\Pfad(5,3),1\endPfad
\Pfad(6,2),2\endPfad
\Pfad(6,2),1\endPfad
\Kreis(0,3.02) \Kreis(5,0)
\Label\lo{A_1}(0,3) \Label\ru{E_1}(5,0)
\Kreis(1,4.02) \Kreis(8,3)
\Label\lo{A_2}(1,4) \Label\ru{E_4}(8,3)
\Kreis(2,5.02) \Kreis(9,4)
\Label\lo{A_3}(2,5) \Label\ru{E_5}(9,4)
\Kreis(4,2.02) \Kreis(6,1)
\Label\lo{A_4}(4,2) \Label\ru{E_2}(6,1)
\Kreis(5,3.02) \Kreis(7,2)
\Label\lo{A_5}(5,3) \Label\ru{E_3}(7,2)
\hskip6cm
$$
\vskip1cm
\caption{}
\label{tiling}
\end{figure}

Suppose now that the parity of $a$ is different from that of $b$ and
$c$, which is 
the case in Theorems~\ref{shifted} and \ref{-shifted}. Since
in this case the core is slightly off the ``truly central'' position
(because the triangle in the ``truly central" position
would not be a lattice triangle; see the definitions in the
Introduction), the starting 
points of the lattice paths originating at boundary points of the core
are changed 
slightly as well. The starting and ending points become
\begin{subequations} \label{shcoord}
\begin{align}
A_i&=(i-1,c+m+i-1),& i&=1,2,\dots ,a, \\
A_i&=\(\frac {a+b-1}{2}+i-a-1,\frac {a+c-1} {2} +i-a-1\),&
i&=a+1,a+2,\dots,a+m,\\
E_j&=(b+j-1,j-1),& j&=1,2,\dots, a+m.
\end{align}
\end{subequations}

In either case, the lozenge
tiling can be recovered from the path family, so that it suffices
to count the families of nonintersecting lattice paths with the
above-mentioned starting and end points.

\medskip
{\em From nonintersecting lattice paths to a determinant.} 
In order to count these families of nonintersecting lattice paths, we
make use of a result due to
Lindstr\"om \cite[Lemma~1]{LindAA} and independently to 
Gessel and Viennot 
\cite[Theorem~1]{gv}. In fact, it is the not so well-known general
form of the result which we need here.
In order to state this result, we introduce some lattice path
notation. 
We write $\P (A \to E)$ for the number of
paths starting at $A$ and ending at $E$.
Given two sets $\mathbf A=\{A_1,\dots,A_n\}$ and $\mathbf
E=\{E_1,\dots,E_n\}$  
of lattice points
and a permutation $\sigma$, we write 
$\P ({\mathbf A} \to {\mathbf E_{\sigma}}, \text {nonint.})$ for the
number of families of $n$ nonintersecting paths with the $i$th path  
running from $A_i$ to $E_{\sigma (i)}$, $i=1,2,\dots,n$. 

Now we can state the main result on nonintersecting lattice paths
(see \cite[Lemma~1]{LindAA} or \cite[Theorem~1]{gv}). 
\begin{lemma} \label{gv}
Let $A_1,A_2,\dots, A_n, E_1, E_2,\dots , E_n$ be points of 
the planar integer lattice.
Then the following identity holds:
\begin{equation} \label{eq:gv}
\det_{1\le i,j \le n}{\(\P (A_i \to E_j)\)}=\sum _{\sigma\in S_n} ^{}{
(\sgn \sigma )\cdot \P ({\mathbf A} \to {\mathbf E_\sigma}, \text {\em nonint.})}.
\end{equation}
\quad \quad \qed
\end{lemma}

\begin{remark} The result in \cite{LindAA}, respectively \cite{gv},
is in fact more general, as it is formulated for paths in an arbitrary
oriented graph. But then the graph must satisfy an acyclicity
condition.
We have not mentioned it in the formulation of the above lemma 
as it is automatically satisfied in our more restricted setting.
\end{remark}

Usually, this lemma is applied in the case that the only permutation
for which nonintersecting lattice paths exist is the identity
permutation, so that the sum on the right-hand side reduces to a
single term, which counts all families $(P_1,P_2,\dots,P_n)$ of 
nonintersecting lattice paths, the $i$th path $P_i$ running from $A_i$
to $E_i$, $i=1,2,\dots,n$. 
(The only exceptions that we are aware of, i.e., applications of
the above formula in the case
where the sum on the right-hand side does not reduce to a single
term, can be found in \cite{FiscAA},
\cite{LindAA}, and \cite{SulaAC}.) This is, however, not exactly the
situation that we encounter in our problem.
Therefore, it seems that Lemma~\ref{gv} is
not suited for our problem. However, 
our choice of starting and end points (see Figure~\ref{tiling}.c)
implies that
nonintersecting lattice paths
are only possible if $m$ consecutive end points ($m$ being the side
length of the equilateral triangle removed from the hexagon) 
are paired with the
starting points from the triangle. So the corresponding permutation
$\sigma$, which describes in which order the starting points are
connected to the end points, 
differs from the identity permutation by a composition of 
cycles of length $m+1$.
Thus, if $m$ is even, we have $\sgn \sigma=1$, so that
the right-hand side in Lemma~\ref{gv} counts exactly all
nonintersecting lattice path families and, thus, all the lozenge
tilings that we are interested in.

On the other hand, if $m$ is odd, then the sign of the permutation
$\sigma$ will not be $1$ always. In fact, as is straightforward to
see, the sign of $\sigma$ is $1$ if the number of paths which pass the
core on the right is even, and is $-1$ otherwise. If this
is translated back to the original lozenge tiling, $T$ say, then it follows
that $\sgn\sigma$ is exactly equal to $(-1)^{n(T)}$, with the statistic
$n(.)$ from Section~\ref{sec:-1}. Thus, in the case that $m$ is odd,
the determinant in Lemma~\ref{gv} gives exactly the
$(-1)$-enumeration of our lozenge tilings.

Since the number of paths from $(x_1,y_1)$ to $(x_2,y_2)$ with positive 
horizontal and negative vertical steps equals the binomial coefficient
$\binom{x_2-x_1+y_1-y_2}{x_2-x_1}$, our findings so far can be summarized as
follows.

\begin{lemma} \label{lem:nonint}
Let $a,b,c,m$ be nonnegative integers, $a,b,c$ having the same parity. 
If $m$ is even, then
the number of lozenge tilings of a hexagon with sides
$a,b+m,c,a+m,b,c+m$, with an equilateral triangle of
side length $m$ removed from its center, equals
\begin{equation} \label{mat1}
\det_{1\le i,j\le a+m} \begin{pmatrix} \dbinom{b+c+m}{b-i+j}&
\text {\scriptsize $1\le i\le a$}\\
\dbinom{\frac {b+c} {2}}{\frac {b+a} {2}-i+j}&
\text {\scriptsize $a+1\le i \le a+m$}
\end{pmatrix}.
\end{equation}
If $m$ is odd, then the weighted count $\sum (-1)^{n(T)}$, where $T$
varies through all the above lozenge tilings, is equal to the
above determinant.
\end{lemma}

\begin{lemma} \label{lem:nonint2}
Let $a,b,c,m$ be nonnegative integers, 
$a$ of parity different from the parity of $b$ and $c$.
If $m$ is even, then
the number of lozenge tilings of a hexagon with sides
$a,b+m,c,a+m,b,c+m$, with an equilateral triangle of
side length $m$ removed that is ``central" in the sense 
that was described in the Introduction, equals
\begin{equation} \label{mat2}
\det_{1\le i,j\le a+m} \begin{pmatrix} \dbinom{b+c+m}{b-i+j}&
\text {\scriptsize $1\le i\le a$}\\
\dbinom{\frac {b+c} {2}}{\frac {b+a+1} {2}-i+j}&
\text {\scriptsize $a+1\le i \le a+m$}
\end{pmatrix}.
\end{equation}
If $m$ is odd, then the weighted count $\sum (-1)^{n(T)}$, where $T$
varies through all the above lozenge tilings,
is equal to the
above determinant.
\end{lemma}

\end{section}

\begin{section}{Polynomiality of the number of lozenge tilings}
\label{sec:poly}

The goal of this section is to establish polynomiality in $m$ --- the side 
of the core ---
of the weighted counts of lozenge tilings considered in
Theorems~\ref{enum}, \ref{shifted}, \ref{-enum}, \ref{-shifted},
provided $a,b,c$ are fixed. Below we just address the case that
$a$, $b$ and $c$ have the same parity (i.e., the case considered in
Theorems~\ref{enum} and \ref{-enum}), the other case being completely
analogous. 

We set up a bijection 
between the lozenge tilings of our cored hexagon and 
nonintersecting lattice paths in a manner different from the one in
the preceding section.
We start by extending all sides of the removed triangle to the left
(if viewed from the interior of the triangle; see
Figure~\ref{cycpaths}, where these extensions are marked as thick
segments).
These segments partition the cored hexagon into three regions.
Furthermore, the
segments cut some of the lozenges in two. (In
Figure~\ref{cycpaths} these lozenges are shaded.)
In each of the three regions, we mark the
midpoints of those edges of the ``cut" lozenges and of those edges along
the border of the region that are not parallel to the ``thick" segments
bordering this region (see
Figure~\ref{cycpaths}). Now, in each of the three regions, we connect
the marked points by ``following" along the lozenges of the tiling, in
the same way as in Section~\ref{sec:outline}
(in the proof of Theorem~\ref{thm:cycl-1}), and in Section~\ref{sec:nonint}
(see Figure~\ref{tiling}.b). 
The lozenge tiling can be recovered from the three nonintersecting 
path families. Thus this defines indeed a bijection.

\newbox\cycpaths
\setbox\cycpaths\hbox{\small
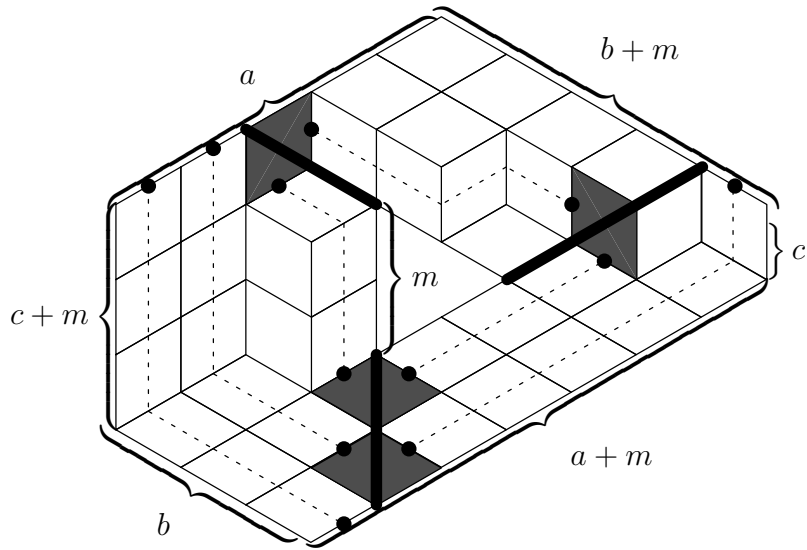
A lozenge tiling of a hexagon with sides $a=5$, $b=3$, $c=1$ and removed triangle 
of side length $m=2$ and the corresponding paths.}
\begin{figure} 
\centertexdraw{
\drawdim truecm \linewd.02
\RhombusA \RhombusA \RhombusA \RhombusA \RhombusB \RhombusA
\move(-.866025 -.5)
\RhombusA \RhombusA \RhombusB \RhombusA \RhombusA \RhombusA 
\move(-1.7305 -1)
\RhombusB \RhombusA \RhombusB \RhombusB \RhombusA \RhombusA 
\move(-2.596 -1.5)
\RhombusB \RhombusB \RhombusB \RhombusA \RhombusA \RhombusA 
\move(-3.462 -2)
\RhombusB \RhombusB \RhombusB \RhombusA \RhombusA \RhombusA 
\move(0 -4)
\RhombusA \RhombusA 
\move(.866025 -3.5)
\RhombusA \RhombusA 
\move(-1.7305 -3)
\rlvec(.866025 -.5)
\move(0 -1)
\rlvec(0 -1)
\move(2.596 -1.5)
\rlvec(0 -1)
\move(4.33 -1.5)
\RhombusC
\move(0 -2)
\rlvec(.866025 -.5)

\move(0 -4)  \rlvec(-.866025 -.5) \rlvec(.866025 -.5)
\lfill f:.3
\rlvec(.866025 .5) \rlvec(-.866025 .5)
\lfill f:.3

\move(0 -5)  \rlvec(-.866025 -.5) \rlvec(.866025 -.5)
\lfill f:.3
\rlvec(.866025 .5) \rlvec(-.866025 .5)
\lfill f:.3

\move(-.866025 -.5)  \rlvec(-.866025 -.5) \rlvec(0 -1)
\lfill f:.3
\rlvec(.866025 .5) \rlvec(0 1)
\lfill f:.3

\move(2.596 -1.5)  \rlvec(.866025 -.5) \rlvec(0 -1)
\lfill f:.3
\rlvec(-.866025 .5) \rlvec(0 1)
\lfill f:.3

\move(-2.165 -1.25)
\vdSchritt \vdSchritt \vdSchritt \hdSchritt \hdSchritt 
%\hhdSchritt \hodaSchritt
\rmove(.866025 0) 
\odaSchritt \odaSchritt \odaSchritt \odaSchritt
\odaSchritt \vdaSchritt
\move(-3.031 -1.75)
\vdSchritt \vdSchritt \vdSchritt \hdSchritt \hdSchritt \hdSchritt
\move(-1.298 -1.75)
%\hvdSchritt 
\hdSchritt \vdSchritt \vdSchritt %\hhdSchritt \hodaSchritt 
\rmove(.866025 0)
\odaSchritt
\odaSchritt \odaSchritt 
%\hvdaSchritt \hhdaSchritt 
\rmove(-.433012 .75)
\hdaSchritt \odSchritt 
\hdaSchritt \hdaSchritt  %\hodSchritt 

\linewd.15
\move(0 -4)
\rlvec(0 -2)
\move(1.73205 -3)
\rlvec(2.598 1.5)
\move(0 -2)
\rlvec(-1.73205 1)
\linewd.02

\ringerl(-1.29904 -1.75)
\ringerl(-2.165 -1.25)
\ringerl(-3.031 -1.75)
\ringerl(-.433012 -4.25)
\ringerl(-.433012 -5.25)
\ringerl(-.433012 -6.25)
\ringerl(.433012 -4.25)
\ringerl(.433012 -5.25)
\ringerl(3.031 -2.75)
\ringerl(4.763 -1.75)
\ringerl(2.598 -2)
\ringerl(-.866025 -1)

\htext(-1.8 -.4){$a$}
\htext(3 -.1){$b+m$}
\htext(5.15 -3){$\left. \vbox{\vskip.5cm} \right\} c$}
\htext(-4.85 -5){$c+m \left\{ \vbox{\vskip1.6cm} \right.$}
\htext(-2.9 -6.4){$b$}
\htext(2.6 -5.5){$a+m$}
\htext(0 -4){$\left. \vbox{\vskip1.1cm} \right\}m$}
\rtext td:60 (3.1 -.9) {$\left. \vbox{\vskip2.75cm} \right\} $}
\rtext td:-60 (-1.6 -.6){$\left\{ \vbox{\vskip2.7cm} \right. $}
\rtext td:-60 (2 -4.7) {$\left. \vbox{\vskip3.7cm} \right\} $}
\rtext td:60 (-2.4 -6){$\left\{ \vbox{\vskip1.7cm}\right. $}
}
\caption{\unhbox\cycpaths}
\label{cycpaths}
\end{figure}

Hence, if we fix the lozenges that are cut in two by the segments, 
the corresponding number of lozenge tilings
which contain these fixed ``cut" lozenges is easily computed by
applying the Lindstr\"om--Gessel--Viennot theorem (Lemma~\ref{gv}) to
each of the three regions separately. This gives a product of three
determinants, one for each region. The total number of lozenge
tilings is then obtained as the sum over all possible choices of
``cut" lozenges (along the segments) of this
product of three determinants.

It is easy to see that each entry in any of the three determinants 
is a binomial coefficient of the
form $\binom{m+x}{y}$, where $x$ and $y$ are independent of $m$.
So the entries are polynomials in $m$, and, hence, the determinants
as well.
The segment which extends the side of the removed
triangle that is parallel to $a$ has length
$\min \{\frac {a+b} {2},\frac {a+c} {2}\}$, which is independent 
of $m$, similarly for the other lines. The total number of
lozenge tilings is thus equal to a sum of polynomials in $m$, where
the range of summation is independent of $m$. Therefore it is itself
a polynomial in $m$, as was claimed.

Basically, the same arguments hold also for $(-1)$-enumeration. The
only difference is that each product of three determinants is
multiplied by a sign,  
depending (according to the definition of our statistic $n$) on the
parity of the number of lozenge sides contained in the northeastern extension 
of the bottom side of the core. However, this number equals the length of this
extension minus the number of lozenges the extension cuts through, and is therefore
again independent of $m$.

\end{section}

\begin{section}{Determinant evaluations, I} \label{sec:det}

In this section we evaluate the determinant in
Lemma~\ref{lem:nonint}. 
The underlying matrix is a mixture of two matrices. 
If we would have to compute the
determinant of just one of the matrices (i.e., if we consider the
case $a=0$ or $m=0$), then the determinant could be easily evaluated
(see \eqref{eq:det-eval}).
However, the mixture is much more difficult to
evaluate.
As it turns out, we
have to distinguish between several cases, depending on the 
parities of $a$ and $m$.

\medskip
It is convenient to take
${(b+c+m)!} \big/{(b+a+m-i)!\,(c+m+i-1)!}$
out of the $i$th row, $i=1,2,\dots,a$, and 
$ {\(\frac {b+c} {2}\)!}\big/{\(\frac {b+3a}
{2}+m-i\)!\(\frac {c-a} {2}+i-1\)!}$
out of the $i$th row, $i=a+1,a+2,\dots,a+m$.
This gives
\begin{multline} \label{det}
\det_{1\le i,j\le a+m} \begin{pmatrix} \dbinom{b+c+m}{b-i+j}&
\text{\scriptsize $1\le i\le a$}\\
\dbinom{\frac {b+c} {2}}{\frac {b+a} {2}-i+j}&
\text{\scriptsize $a+1\le i \le a+m$}
\end{pmatrix}\\
=\prod _{i=1} ^{a}{\frac {(b+c+m)!} {(b+a+m-i)!\,(c+m+i-1)!}}
\prod _{i=a+1} ^{a+m}{\frac {\(\frac {b+c} {2}\)!}{\(\frac {b+3a}
{2}+m-i\)!\(\frac {c-a} {2}+i-1\)!}}\\
\times\det_{1\le i,j\le a+m} \begin{pmatrix} (c+m+i-j+1)_{j-1}\,(b-i+j+1)_{a+m-j}&
\text{\scriptsize $1\le i\le a$}\\
\(\frac {c-a} {2}+i-j+1\)_{j-1}\(\frac {b+a}{2}-i+j+1\)_{a+m-j}&
\text{\scriptsize $a<i\le a+m$}
\end{pmatrix}.
\end{multline}

Thus it suffices to evaluate the determinant on the right-hand side.
The advantage is that this determinant is a polynomial in $b$ and
$c$. This enables us to apply the ``identification of factors"
method, as proposed in \cite[Sec.~2.4]{KratBN}.
The four lemmas below address the four different cases, as $a$ and $m$
vary through all possible parities.

\begin{lemma} \label{lem:amgerade}
Let $a$ and $m$ be both even nonnegative integers.
Then 
\begin{multline} \label{eq:amgerade}
\det_{1\le i,j\le a+m} \begin{pmatrix} (c+m+i-j+1)_{j-1}\,(b-i+j+1)_{a+m-j}&
\text{\scriptsize $1\le i\le a$}\\
\(\frac {c-a} {2}+i-j+1\)_{j-1}\(\frac {b+a}
{2}-i+j+1\)_{a+m-j}&\text{\scriptsize $a<i\le a+m$}
\end{pmatrix}\\
=\frac {\h(a+m)\h(\frac {a} {2})^2\h(\frac {m} {2})^2} {\h(\frac {a+m}
{2})^2\,2^{m(a+m-1)/2}}
 \prod _{k=1} ^{ {m}/ {2}}{\(\tfrac {b} {2}+k\)^2_{ {a}/ {2}}
\(\tfrac {c} {2}+k\)^2_{ {a}/ {2}}}
\prod _{k=0} ^{ {a}/ {2}-1}{({b+c+m}+2k+ {1} )^{a-2k-1}}\\
\times
\prod _{k=1} ^{ {a}/ {2}-1}{({b+c}+2m+2k)^{a-2k}}
\prod _{k= {m}/ {2}+1} ^{m}{({b+c}+2k)^{a+m-k}}
\prod _{k=1} ^{ {m}/ {2}}{({b+c}+2k)^{m-k}}.
\end{multline}
\end{lemma}

\begin{proof}
Let us denote the determinant in \eqref{eq:amgerade} by $D_1(b,c)$.

We proceed in several steps. An outline is as follows. 
The determinant $D_1(b,c)$ is obviously a polynomial in
$b$ and $c$. In Steps~1--5 we show that the right-hand side of
\eqref{eq:amgerade} divides $D_1(b,c)$ as a polynomial in $b$ and
$c$. In Step~6 we show that the degree of $D_1(b,c)$ as a polynomial
in $b$ is at most $\binom {a+m}2$. Of course, the same is true for
the degree in $c$. On the other hand, the
degree of the right-hand side of \eqref{eq:amgerade} as a polynomial
in $b$ is exactly $\binom {a+m}2$. It follows that $D_1(b,c)$ must
equal the right-hand side of \eqref{eq:amgerade} times a
quantity which does not depend on $b$. This quantity must be
polynomial in $c$. But, in fact, it cannot depend on $c$ as well,
because, as we just observed, the degree in $c$ of the right-hand side of
\eqref{eq:amgerade} is already equal to the maximal degree in $c$ of 
$D_1(b,c)$. Thus, this quantity is a constant with respect to $b$ and
$c$. That this constant is equal to 1 is finally shown in
Step~7, by evaluating the determinant $D_1(b,c)$ for $b=c=0$.

Before we begin with the detailed description of the individual
steps, we should explain the odd looking occurrences of ``$e\equiv
a$ mod 2" below (e.g., in Step~1(a)--(d)). Clearly, in the present
context this means ``$e\equiv 0$ mod 2", as $a$ is even by assumption.
However, Steps~1--6 will also serve as a model for the proofs of the
subsequent Lemmas~\ref{lem:amung}--\ref{lem:-aung}. Consequently, 
formulations are
chosen so that they remain valid without change at the corresponding
places. In particular, 
in the context of the proofs of Lemmas~\ref{lem:amung} and
\ref{lem:-aung}, the statement
``$e\equiv a$ mod 2" will mean ``$e\equiv 1$ mod 2".

\smallskip
{\em Step 1. $\prod _{k=1} ^{ {m}/ {2}}{\(\frac {b} {2}+k\)^2_{ {a}/
{2}}}\,{\(\frac {c} {2}+k\)^2_{ {a}/
{2}}}$ divides the determinant.} 
The original
determinant is symmetric in $b$ and $c$ for combinatorial reasons.
The factors which were taken out of the determinant in \eqref{det} are also
symmetric in $b$ and $c$
(this can be seen by reversing all the products involving
$c$). Therefore it suffices to check that the linear factors involving
$b$ divide $D_1(b,c)$, i.e., that the product 
$\prod _{k=1} ^{ {m}/ {2}}{\(\frac {b} {2}+k\)^2_{ {a}/{2}}}$ 
divides $D_1(b,c)$.

We distinguish between four subcases, labeled below as (a), (b),
(c), and (d).

\smallskip
(a)
$(b+e)^e$ divides $D_1(b,c)$ for $1\le e \le \min\{a,m\}$,
$e\equiv a$ mod 2: This follows from the easily verified fact that
$(b+e)$ is a factor of each entry in the first $e$ columns of
$D_1(b,c)$.

\smallskip
(b) $(b+e)^m$ divides $D_1(b,c)$ for $m< e < a$, $e \equiv a$
mod 2: We prove this by finding $m$ ``different" linear combinations
of the columns of $D_1(b,c)$ which vanish for $b=-e$. By the term
``different" we mean that these linear combinations are themselves
linearly independent. (Equivalently, we find $m$ linearly independent
vectors in the kernel of the linear operator defined by the matrix
underlying $D_1(-e,c)$.) See Section~2 of \cite{KratBI}, and in
particular the Lemma in that section, for a formal
justification of this procedure.

To be precise,
we claim that the following equation holds for $s=1,2,\dots,m$,
\begin{equation} \label{eq:lin1b}
\sum _{j=1} ^{e+s-m}{\binom{e-m+s-1}{j-1}\frac
{(c+a-e-s+2m+1)_{e+s-j-m}} {(a-e-s+2m+1)_{e+s-j-m}}\cdot (\text
{column $j$ of $D_1(-e,c)$})}=0.
\end{equation}

Since the entries of $D_1(b,c)$ have a split definition (see
\eqref{eq:amgerade}), for the proof of the above equation we have to
distinguish between two cases. If we restrict \eqref{eq:lin1b} to the
$i$th row, $i\le a$, then \eqref{eq:lin1b} becomes
\begin{multline} \label{eq:lin1b-kleiner}
\sum _{j=1} ^{e+s-m}\binom{e-m+s-1}{j-1}\frac
{(c+a-e-s+2m+1)_{e+s-j-m}} {(a-e-s+2m+1)_{e+s-j-m}}\\
\cdot(c+m+i-j+1)_{j-1}\,(-e-i+j+1)_{a+m-j}=0,
\end{multline}
whereas on restriction to the
$i$th row, $i> a$, equation \eqref{eq:lin1b} becomes
\begin{multline} \label{eq:lin1b-groesser}
\sum _{j=1} ^{e+s-m}\binom{e-m+s-1}{j-1}\frac
{(c+a-e-s+2m+1)_{e+s-j-m}} {(a-e-s+2m+1)_{e+s-j-m}}\\
\cdot\(\tfrac {c-a} {2}+i-j+1\)_{j-1}\(\tfrac {-e+a}
{2}-i+j+1\)_{a+m-j}=0.
\end{multline}

First, let $i\le a$.
Here and in the following, we make use of
the usual hypergeometric notation
\begin{equation} \label{eq:hyper}
{}_r F_s\!\left[\matrix a_1,\dots,a_r\\ b_1,\dots,b_s\endmatrix; 
z\right]=\sum _{k=0} ^{\infty}\frac {\po{a_1}{k}\cdots\po{a_r}{k}}
{k!\,\po{b_1}{k}\cdots\po{b_s}{k}} z^k\ .
\end{equation}
In this notation, the sum on the
left-hand side of \eqref{eq:lin1b-kleiner} reads
\begin{multline*}
\frac{(2 - e - i) _{-1 + a + m} \,
     ( 1 + a + c - e + 2\,m - s)_{-1 + e - m + s} } 
   {(1 + a - e + 2\,m - s) _{-1 + e - m + s} }\\
\times
{} _{3} F _{2} \!\left [ \begin{matrix} { 1 - c - i - m, 1 - e + m - s, 1 -
      a - m}\\ { 1 - a - c - m, 2 - e - i}\end{matrix} ; {\displaystyle
      1}\right ]. 
\end{multline*}
Next we use a transformation formula due to Thomae \cite{Th} (see
also \cite[(3.1.1)]{GaRa}),
\begin{equation} \label{T3207}
{} _{3} F _{2} \!\left [ \begin{matrix} { A, B, -n}\\ { D, E}\end{matrix} ;
   {\displaystyle 1}\right ]  = 
  {\frac{ ({ \textstyle E-B}) _{n} }{({ \textstyle E}) _{n} }}
{} _{3} F _{2} \!\left [ \begin{matrix} { -n, B, D-A}\\ { D, 1 + B - E
       - n}\end{matrix} ; {\displaystyle 1}\right ],
\end{equation}
where $n$ is a nonnegative integer.
This gives
\begin{multline*} 
\frac { (1 + a + c - e + 2\,m - s) _{e-m+s-1} \,
     (1-i-m+s) _{a+m-1} 
} {(1 + a - e + 2 m - s) _{e-m+s-1} }\\
\times 
{} _{3} F _{2} \!\left [ \begin{matrix} { 1 - a - m, 1 - e + m - s, -a + i}\\ { 1 -
      a - c - m, 1 - a + i - s}\end{matrix} ; {\displaystyle 1}\right ].
\end{multline*}
The factor $(1-i-m+s)_{a+m-1}$ vanishes for $i\le a$ and the
denominator is never zero, so the sum in \eqref{eq:lin1b-kleiner} 
equals zero, as desired.

We proceed similarly in order to prove \eqref{eq:lin1b-groesser} for
$i>a$. The hypergeometric form of the sum in \eqref{eq:lin1b-groesser} is
\begin{multline*} 
\frac{(2 + {\frac a 2} - {\frac e 2} - i) _{a+m-1} \,
     ({ 1 + a + c - e + 2\,m - s}) _{e-m+s-1} }
{(1 + a - e + 2\,m - s) _{e-m+s-1} }\\
\times
{} _{3} F _{2} \!\left [ \begin{matrix} 
{ 1 + {\frac a 2} - {\frac c 2} - i, 1
      - a - m, 1 - e + m - s}\\ { 1 - a - c - m, 2 + {\frac a 2} -
{\frac e 2} - i}\end{matrix} ; {\displaystyle 1}\right ] .
\end{multline*}
Using the transformation formula \eqref{T3207} again, we get
\begin{multline*} 
(1+\tfrac {a} {2}+\tfrac {e} {2}-m-i+s)_{a+2m-e-s} \,
({  1 + {\tfrac {3a} 2} - {\tfrac e 2} - i + m}) _{e-m+s-1}\\
\times 
\frac {
({  1 + a + c - e + 2\,m - s}) _{e-m+s-1} } 
{ ({  1 + a - e + 2\,m - s}) _{-1 + e - m + s} }
\\
\times
{} _{3} F _{2} \!\left [ \begin{matrix} { 1 - e + m - s, 1 - a - m,
      {-\frac {3a} 2} - {\frac c 2} + i - m}\\ { 1 - a - c - m, 1 -
      {\frac {3\,a} 2} - {\frac e 2} + i - s}\end{matrix} ; {\displaystyle
      1}\right ] .
\end{multline*}
This expression is zero, because 
the factor $(1+\frac {a} {2}+\frac {e} {2}-m-i+s)_{a+2m-e-s}$
vanishes for $i>a$ (it is here where we need $e\equiv a$ mod 2,
because this
guarantees that $1+\frac {a} {2}+\frac {e} {2}-m-i+s$ is an integer).
So the sum in \eqref{eq:lin1b-groesser} equals zero, as desired.

\smallskip
(c)  $(b+e)^a$ divides $D_1(b,c)$ for $a < e < m$, $e\equiv a$
mod 2: Proceeding in the spirit of case (b), we prove this by finding
$a$ linear combinations
of the columns of $D_1(b,c)$ which vanish for $b=-e$.
To be precise, we 
claim that the following equation holds for $s=1,2,\dots,a$:
\begin{equation} \label{eq:lin1c}
\sum _{j=1} ^{\frac {e-a} {2}+s}{\binom{\frac {e-a} {2}+s-1}{j-1}
\frac {\(\frac {c} {2}+m-\frac {e} {2}+a-s+1\)_{ {(e-a)}/ {2}+s-j}} 
{\(m+\frac {3a-e} {2}-s+1\)_{ {(e-a)}/ {2}+s-j}}\cdot
(\text {column $j$ of $D_1(-e,c)$})}=0.
\end{equation}

In order to prove this equation, we first restrict it to the $i$th
row, $i\le a$. Then, in hypergeometric notation, the left-hand side
reads
\begin{multline*}
\frac {(2 - e - i) _{a+m-1} \,
({1+a+{\frac c  2}-{\frac e  2}+m-s}) _{-1-{\frac a  2}+{\frac e  2} + s}} 
{( 1 + {\frac {3\,a}  2} - {\frac e  2} + m - s) _{-1 -
{\frac a  2}
     + {\frac e  2} + s}}\\
\times
{} _{3} F _{2} \!\left [ \begin{matrix} { 1 - a - m, 1 - c - i - m, 1 +
      {\frac a  2} - {\frac e  2} - s}\\ { 1 - {\frac a  2} - {\frac c  2} - m, 2 - e
      - i}\end{matrix} ; {\displaystyle 1}\right ].
\end{multline*} 
We apply the transformation formula \eqref{T3207} and get 
\begin{multline*}
(1+c-e+m) _{{\frac e 2} -{\frac a  2}+s-1} \,
(1 + a + {\frac c  2} - {\frac e  2} + m - s) _{-1 -
      {\frac a  2} + {\frac e  2} + s} \\
\times
\frac {
     (1 - {\frac a  2} - {\frac e  2} - i + s) _{{\frac {3\,a}  2} -
      {\frac e  2} + m - s}} 
{( 1 + {\frac {3\,a}  2} -{\frac e  2}+m-s) _{-1 -{\frac a  2}
     + {\frac e  2} + s}}\,
{} _{3} F _{2} \!\left [ \begin{matrix} { 1 + {\frac a  2} -
{\frac e  2} - s, 
1 - c - i - m, {\frac a  2} - {\frac c  2}}\\ { 1 - {\frac a  2} -
{\frac c  2} -
      m, 1 + {\frac a  2} - c + {\frac e  2} - m - s}\end{matrix} ; {\displaystyle
      1}\right ] .
\end{multline*}
This expression is zero because the factor 
$( 1 - {\frac a  2} - {\frac e 2} -i+s)_{{\frac {3a}  2} -
  {\frac e  2} + m - s} $ vanishes.

If instead we restrict the left-hand side of \eqref{eq:lin1c} to the
$i$th row, $i>a$, and convert it into hypergeometric form, then we
obtain
\begin{multline*} 
\frac {({ 2 + {\frac a  2} - {\frac e  2} - i}) _{a + m-1} \,
(1 + a + {\frac c  2} - {\frac e  2} + m - s) _{-1 -
      {\frac a  2} + {\frac e  2} + s}} 
{(1 + {\frac {3\,a}  2} - {\frac e  2} + m - s) _{-1 -
{\frac a  2} + {\frac e  2} + s}}\\
\times
{} _{3} F _{2} \!\left [ \begin{matrix} { 1 + {\frac a  2} - {\frac c  2} - i, 1
      - a - m, 1 + {\frac a  2} - {\frac e  2} - s}\\  
   1 - {\frac a  2} - {\frac c  2} - m,  2 + {\frac a  2} -
   {\frac e 2} - i\end{matrix} ; {\displaystyle
      1}\right ] .
\end{multline*}
We apply again the transformation formula \eqref{T3207}. This gives
\begin{multline*} 
{\frac { ({ \textstyle 1 + a + {\frac c 2} - {\frac e 2} + m - s}) _{-1 -
{\frac a 2}
      + {\frac e 2} + s} \,({ \textstyle 1 - i + s}) _{-1 + a + m} } 
   {({ \textstyle 1 + {\frac {3a} 2} - {\frac e 2} + m - s}) _{-1 -
{\frac a 2}
     + {\frac e 2} + s} }}\\
\times
{} _{3} F _{2} \!\left [ \begin{matrix} { 1 + {\frac a 2} - {\frac e 2} - s, 1 - a -
      m, -a + i - m}\\ { 1 - {\frac a 2} - {\frac c 2} - m, 1 - a + i - m -
      s}\end{matrix} ; {\displaystyle 1}\right ].
\end{multline*}
This expression is zero because the factor 
$({ \textstyle 1 - i + s}) _{-1 + a + m}$ vanishes for
$a+1\le i \le a+m$. So the sum in \eqref{eq:lin1c}
equals zero, as desired.

\smallskip
(d) $(b+e)^{a+m-e}$ divides $D_1(b,c)$ for $\max\{a,m\} \le e \le
a+m-1$, $e\equiv a$ mod 2: Still proceeding in the spirit of case
(b), this time we find
$a+m-e$ linear combinations
of the {\em rows} of $D_1(b,c)$ which vanish for $b=-e$.
To be precise,
we claim that the following equation holds for $s=1,2,\dots, a+m-e$:
\begin{multline}\label{eq:lin1d}
\sum _{i=1} ^{s}\binom{s-1}{i-1}
\frac {(-1)^{i}\(\frac {c-e} {2}+1\)_{a+m-s}\(\frac {c-e}
{2}+m\)_{i-1}} {(1+c-e+m)_{a+m-s+i-1}} \\
\cdot(\text {row $(a+m-e-s+i)$
of $D_1(-e,c)$})\\
+(\text {row $\(m+\tfrac {3a} {2}-\tfrac {e} {2}-s+1\)$ of $D_1(-e,c)$})=0.
\end{multline}
In the sum, it is only the first $a$ rows which are involved, whereas
the extra term is a row out of the last $m$ rows of the determinant.
Therefore, by restriction to the $j$th column, we see that it is
equivalent to
\begin{multline}\label{eq:lin1d-explizit}
\sum _{i=1} ^{s}\binom{s-1}{i-1}
\frac {(-1)^{i}\(\frac {c-e} {2}+1\)_{a+m-s}\(\frac {c-e}
{2}+m\)_{i-1}} {(1+c-e+m)_{a+m-s+i-1}} \\
\cdot (a+c+2m-e-s+i-j+1)_{j-1}\,(-a-m+s-i+j+1)_{a+m-j}
\\
+\(\tfrac {c-e} {2}+a+m-s-j+2\)_{j-1}\,(-a-m+s+j)_{a+m-j}=0.
\end{multline}
We treat the cases $j\le a+m-s$ and $j>a+m-s$ separately.
For $j\le a+m-s$ the factor $(-a-m+s-i+j+1)_{a+m-j}$, which appears in the
sum, is zero for all the
summands, as well is the factor $\(-a-m+s+j\)_{a+m-j}$,
which appears in the extra term in \eqref{eq:lin1d-explizit}.

For $j>a+m-s$ we convert the sum in \eqref{eq:lin1d-explizit}
into hypergeometric form and get
\begin{multline*} 
\frac {-(1 + {\frac c 2} - {\frac e 2}) _{a + m - s} \,
      ( 2 + a + c - e - j + 2\,m - s) _{-1 + j} \,
      ( -a + j - m + s) _{a - j + m} } 
{(1 + c - e + m) _{a + m - s} }\\
\times
{} _{2} F _{1} \!\left [ \begin{matrix} { {\frac c 2} - {\frac e 2} + m, 1 +
       a - j + m - s}\\ { 2 + a + c - e - j + 2\,m - s}\end{matrix} ;
       {\displaystyle 1}\right ]. 
\end{multline*}

We can evaluate the ${}_{2}F_{1}$-series by the Chu--Vandermonde summation 
formula (see \cite[(1.7.7), Appendix~(III.4)]{Sl}),
\begin{equation} \label{vand}
{} _{2} F _{1} \!\left [ \begin{matrix} { A, -n}\\ { C}\end{matrix} ;
   {\displaystyle 1}\right ]  = 
  {\frac{({ \textstyle C-A}) _{n} }{({ \textstyle C}) _{n} }},
\end{equation}
where $n$ is a nonnegative integer. Thus we get
\begin{equation} \label{resd}
- (2 + a + {\tfrac c 2} - {\tfrac e 2} - j + m - s) _{j -1} \,
(-a + j - m + s) _{a - j + m}.
\end{equation}
It is easily seen that adding the extra term in \eqref{eq:lin1d-explizit} 
gives zero.

\smallskip
{\em Step 2.
$\prod _{k=0} ^{ {a}/ {2}-1}{(b+c+m+2k+1)^{a-2k-1}}$ divides the
determinant.} We find $e+1$ linear combinations
of the rows of $D_1(b,c)$ which vanish for $b=-c-a-m+1+e$.
To be precise, we claim that the following
equation holds for $0\le e \le a-2$, $s=1,2,\dots,e+1$:
\begin{multline} \label{eq:lin2}
\sum _{i=1} ^{a-e-1}{\(\frac {(c+m+i)_{a-e-i+s-1}}
{(c-e-1+i)_{a-e-i+s-1}}\binom{a-e-2}{i-1}\frac {(s)_{a-e-1}\,(-1)^i}
{(s-i+a-e-1)(a-e-2)!}\right.}\\
\cdot (\text {row $i$ of $D_1(-c-a-m+1+e,c)$})\bigg)\\
+(-1)^{a-e-1} \cdot(\text {row $(a-e-1+s)$ of $D_1(-c-a-m+1+e,c)$})=0.
\end{multline}

Restricted to the $j$th column, and converted into hypergeometric
notation, the sum in \eqref{eq:lin2} reads
\begin{multline*} 
-\frac {({ 1 - a - c + e + j - m}) _{a - j + m} \,
     (1 + c + m) _{a - e + s-2} \,
     (2 + c - j + m) _{j-1} \,
     (s) _{a - e-2}} 
{(1) _{a-e-2} \,
     ( c - e) _{a-e+s-2} }\\
\times
{} _{3} F _{2} \!\left [ \begin{matrix} { 2 - a + e - s, a + c - e - j + m,
      2 - a + e}\\ { 3 - a + e - s, 2 + c - j + m}\end{matrix} ;
      {\displaystyle 1}\right ] .
\end{multline*}

Here we use the Pfaff--Saalsch\"utz summation formula
(see \cite[(2.3.1.3), Appendix~(III.2)]{Sl})
\begin{equation} \label{saal}
{} _{3} F _{2} \!\left [ \begin{matrix} { A, B, -n}\\ { C, 1 + A + B - C -
   n}\end{matrix} ; {\displaystyle 1}\right ]  = 
  {\frac{({ \textstyle C-A}) _{n} \,({ \textstyle C-B}) _{n} } 
    {({ \textstyle C}) _{n} \,({ \textstyle C-A - B}) _{n} }},
\end{equation}
where $n$ is a nonnegative integer.
Thus we get
\begin{multline*} 
(-1)^{a-e-1}({ \textstyle 1 + c + m}) _{-2 + a - e + s} \,
     ({ \textstyle 2 + c - j + m}) _{-1 + j} \\
\times
{\frac {  ( 3 - 2\,a - c + 2\,e + j - m - s) _{-2 + a - e}  } 
   {({ \textstyle c - e}) _{-2 + a - e + s} \,
     ({ \textstyle 1 - c + e}) _{-2 - e + j - m} }}.
\end{multline*}
It is easily verified that adding the $j$th coordinate of the extra
term in \eqref{eq:lin2} gives zero, as desired.
For now, we need equation \eqref{eq:lin2} only for even $e$.

\smallskip
{\em Step 3. $\prod _{k=1} ^{ {a}/ {2}-1}{(b+c+2m+2k)^{a-2k}}$ divides
the determinant.} We find $e$ linear combinations
of the columns of $D_1(b,c)$ which vanish for $b=-c-2m-a+e$. To be
precise, 
we claim that the following equation holds for $0<e\le a$, $e\equiv
a$ mod 2, and $s=1,2,\dots,e$:
\begin{equation} \label{eq:lin3}
\sum _{j=s} ^{a+m+s-e}{\binom{a+m-e}{j-s}\cdot (\text {column $j$ of
$D_1(-c-2m-a+e,c)$})}=0.
\end{equation}

Restricted to the $i$th row, $i\le a$, and converted into hypergeometric
notation, the left-hand side sum in \eqref{eq:lin3} reads
\begin{multline*} 
  ( 1 + c + i + m - s) _{s-1} \,
  ({ \textstyle 1 - a - c + e - i - 2\,m + s}) _{a + m - s} \\
\times {} _{2} F _{1} \!\left [ \begin{matrix} { -c - i - m + s, -a + e - m}\\ { 1 -
   a - c + e - i - 2\,m + s}\end{matrix} ; {\displaystyle 1}\right ].
\end{multline*}
This is summable by the Chu--Vandermonde summation formula
\eqref{vand}. We get
$${\frac {({ \textstyle 1 - a + e - m}) _{a - e + m} \,
     ({ \textstyle 1 + c + i + m - s}) _{s-1} } 
   {({ \textstyle 1 - c + e - i - m}) _{s-e} }}.
$$
This expression equals zero because the factor 
$(1-a+e-m)_{a-e+m}$ vanishes.

On the other hand, if $i>a$, the left-hand side sum in
\eqref{eq:lin3}, restricted to the $i$th row and converted
into hypergeometric from, reads
\begin{multline*} 
({ 1 - {\tfrac a 2} + {\tfrac c 2} + i -
   s}) _{s-1} \,(1-{\tfrac c 2}+{\tfrac e 2}-i-m+s) _{a+m-s} \\
\times
{} _{2} F _{1} \!\left [ \begin{matrix} { {\frac a 2} - {\frac c 2} - i + s, -a
   + e - m}\\ { 1 - {\frac c 2} + {\frac e 2} - i - m + s}\end{matrix} ;
   {\displaystyle 1}\right ].
\end{multline*}
The Chu--Vandermonde summation formula \eqref{vand} turns this
expression into
$${\frac {({ \textstyle 1 - {\frac a 2} + {\frac e 2} - m}) _{a - e + m} \,
     ({ \textstyle 1 - {\frac a 2} + {\frac c 2} + i - s}) _{s-1} } 
   {({ \textstyle 1 + a - {\frac c 2} + {\frac e 2} - i}) _{-e + s} }}.
$$
This expression is zero because the factor 
$(1-\frac {a} {2}+\frac {e} {2}-m)_{a-e+m}$ vanishes
for $e\equiv a$ mod $2$. So the sum in \eqref{eq:lin3} is zero, as
desired.

\smallskip
{\em Step 4. $\prod _{k= {m}/ {2}+1} ^{m}{(b+c+2k)^{a+m-k}}$ divides the
determinant.} We find $a+m-e$ linear combinations
of the columns of $D_1(b,c)$ which vanish for $b=-c-2e$. To be
precise, we claim that the following equation holds for $m/2<e\le m$ 
and $s=1,2,\dots, a+m-e$:
\begin{equation} \label{eq:lin4}
\sum _{j=s} ^{s+e}{\binom{e}{j-s}\cdot
(\text {column $j$ of $D_1(-c-2e,c)$})}=0.
\end{equation}

Restricted to the $i$th row, $i\le a$, and converted into hypergeometric
notation, the left-hand side sum in \eqref{eq:lin4} reads
\begin{equation*} 
  ({ \textstyle 1 + c + i + m - s}) _{s-1} \,
  ({ \textstyle 1 - c - 2e - i + s}) _{a + m - s} \,
{} _{2} F _{1} \!\left [ \begin{matrix} { -c - i - m + s, -e}\\ { 1
   - c - 2e - i + s}\end{matrix} ; {\displaystyle 1}\right ].
\end{equation*}
The result after application of the Chu--Vandermonde summation
formula \eqref{vand} is
$$\frac {({ \textstyle 1 - 2e + m}) _{{e}} \,
     (1 + c + i + m - s) _{s-1} } 
   {(1 + a - c - 2e - i + m) _{-a + {e} - m + s}}.
$$
This expression equals zero because the factor 
$(1-2e+m)_{ e}$ vanishes.

On the other hand, if $i>a$, the left-hand side sum in
\eqref{eq:lin4}, restricted to the $i$th row and converted
into hypergeometric from, reads
\begin{multline*} 
  (1 - {\tfrac a 2} + {\tfrac c 2} + i - s) _{s-1} \,
  (1 + {\tfrac a 2} - {\tfrac c 2} - {e} - i + s) _{a + m - s} \,
{} _{2} F _{1} \!\left [ \begin{matrix} { {\frac a 2} - {\frac c 2} - i + s,
   {-e}}\\ { 1 + {\frac a 2} - {\frac c 2} - {e} - i +s}
\end{matrix} ; {\displaystyle 1}\right ] .
\end{multline*}
Chu-Vandermonde summation \eqref{vand} yields
$$\frac {(1 - {e}) _{{e}} \,
     ({  1 - {\frac a 2} + {\frac c 2} + i - s}) _{s-1} } 
   {(  1 + {\frac {3\,a} 2} - {\frac c 2} - {e} - i + m) _{-a
     + {e} - m + s}}.
$$
This expression is zero because the factor 
$(1-e)_{ e}$ vanishes. So the sum in
\eqref{eq:lin4} is zero, as desired.

\smallskip
{\em Step 5. $\prod _{k=1} ^{ {m}/ {2}}{(b+c+2k)^{m-k}}$ divides the
determinant.} We find $e$ linear combinations
of the rows of $D_1(b,c)$ which vanish for $b=-c-2m+2e$. To be
precise, we claim that the following equation holds for 
$e\le m-1$ and $s=1,2,\dots,e$:
\begin{multline} \label{eq:lin5}
\sum _{i=1} ^{m-s+1}(-1)^i\binom{m-s}{i-1}\frac {\(\frac {c} {2}+\frac
{a} {2}+i\)_{m-s-i+1}} {\(\frac {c} {2}-\frac {a}
{2}-e+i\)_{m-s-i+1}}\\
\cdot (\text {row $(a+i)$ of
$D_1(-c-2m+2e,c)$})=0.
\end{multline}

Restricted to the $j$th row, and converted into hypergeometric
notation, the left-hand side sum in \eqref{eq:lin5} reads
\begin{multline*} 
-\frac {    (1 + {\frac a 2} + {\frac c 2}) _{m - s} \,
      (2 + {\frac a 2} + {\frac c 2} - j) _{j-1} \,
      ( -{\frac {a} 2} - {\frac c 2} + e + j - m) _{a - j + m} } 
{(1 - {\frac a 2} + {\frac c 2} - e) _{m - s} }\\
\times
{} _{2} F _{1} \!\left [ \begin{matrix} { 1 + {\frac a 2} + {\frac c 2} - e -
       j + m, -m + s}\\ { 2 + {\frac a 2} + {\frac c 2} - j}\end{matrix} ;
       {\displaystyle 1}\right ] .
\end{multline*}
After applying Chu--Vandermonde summation \eqref{vand} again, we obtain 
\begin{multline*}-( 1 + {\tfrac a 2} + {\tfrac c 2}) _{m - s} 
      (  1 + e - m) _{m - s}\\
\times
\frac { 
      ( {-\frac {a} 2} - {\frac c 2} + e + j - m) _{a - j + m} 
      (  2 + {\frac a 2} + {\frac c 2} - j + m - s) _{j-m+s-1} } 
{(1 - {\frac a 2} + {\frac c 2} - e) _{m - s}
}.
\end{multline*}
This expression equals zero because the factor $(1+e-m)_{m-s}$
vanishes. So the sum in \eqref{eq:lin5} is zero, as desired.

\smallskip
{\em Step 6. Determination of the degree of $D_1(b,c)$ as a polynomial
in $b$}. Obviously the degree of the $(i,j)$-entry of $D_1(b,c)$ as a
polynomial in $b$ is $a+m-j$. Therefore, if we expand the determinant
$D_1(b,c)$ according to its definition as a sum over permutations,
each term in this expansion has degree $\binom {a+m}2$ in $b$. Hence,
$D_1(b,c)$ itself has degree at most $\binom {a+m}2$ in $b$.

\newbox\bcbox
\setbox\bcbox\hbox{\small
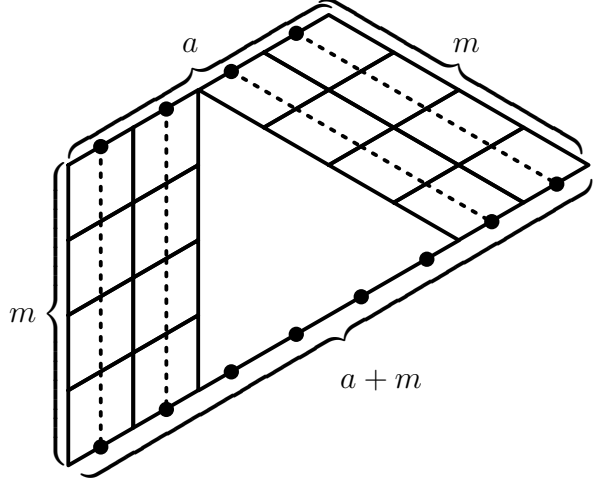
The unique lozenge tiling for $b=c=0$}
\begin{figure} 
\centertexdraw{
\drawdim truecm 
\linewd.05
\move(0 0)
\RhombusA  \RhombusA \RhombusA \RhombusA 
\move(-.866025 -.5)
\RhombusA \RhombusA \RhombusA \RhombusA 
\move(-1.73205 -1)
\RhombusB \RhombusB \RhombusB \RhombusB 
\move(-2.598 -1.5)
\RhombusB \RhombusB \RhombusB \RhombusB
\move(2.598 -2.5)
\rlvec(-3.464 -2)

\ringerl(.433012 .25) 
\hdSchritt \hdSchritt \hdSchritt \hdSchritt   
\ringerl(-.433012 -.25)
\hdSchritt \hdSchritt \hdSchritt \hdSchritt 
\ringerl(-1.299037 -.75) 
\vdSchritt \vdSchritt \vdSchritt \vdSchritt  
\ringerl(-2.165062 -1.25) 
\vdSchritt \vdSchritt \vdSchritt \vdSchritt

\ringerl(3.897 -1.75)
\ringerl(3.031 -2.25) 
\ringerl(2.165062 -2.75)
\ringerl(1.299  -3.25) 
\ringerl(.433012 -3.75)
\ringerl(-.433012 -4.25) 
\ringerl(-1.299 -4.75)
\ringerl(-2.165 -5.25)

\rtext td:-60 (.8 -3.6) {$\left. \vbox{\vskip4.1cm} \right\}$}
\rtext td:-60 (-1.1 -.2){$\left\{ \vbox{\vskip2.2cm} \right.$ }

\rtext td:60 (2.5 -.5) {$\left. \vbox{\vskip2.2cm} \right \}$}
\htext (-3.4 -5.5) {$m \left\{ \vbox{\vskip2.1cm} \right.$} 

\htext(-1.1 0){$a$}
\htext(2.5 0){$m$}
\htext(1 -4.5){$a+m$}
}
\caption{\unhbox\bcbox}
\label{bc0} 
\end{figure}

\smallskip
{\em Step 7. Computation of the multiplicative constant}. As we
observed at the beginning of this proof, Steps~1--6 show that the
determinant $D_1(b,c)$ is equal to the right-hand side of
\eqref{eq:amgerade} up to multiplication by a constant. To determine this
constant, it suffices to compute $D_1(b,c)$ for some particular
values of $b$ and $c$. We choose $b=c=0$. The value of $D_1(0,0)$
is most easily determined by going back, via \eqref{det} and
Lemma~\ref{lem:nonint}, to the origin of
the determinant $D_1(b,c)$, which is enumeration of lozenge tilings.
Figure~\ref{bc0} shows the typical situation for $b=c=0$. As the
figure illustrates, there is exactly {\em one} lozenge tiling of the
region. Hence, by Lemma~\ref{lem:nonint}, it follows that the determinant
\eqref{mat1} must be equal to $1$ for $b=c=0$. If we substitute this into
\eqref{det}, we have evaluated $D_1(b,c)$, which is the determinant
on the right-hand side of \eqref{det}, for $b=c=0$. It is then a
routine task to check that the result agrees exactly with the
right-hand side of \eqref{eq:amgerade} for $b=c=0$.

\medskip
This completes the proof of the theorem.
\end{proof}

\begin{lemma} \label{lem:amung}
Let $a$ and $m$ be nonnegative integers, $a$ odd and $m$ even.
Then 
\begin{multline} \label{eq:amung}
\det_{1\le i,j\le a+m} \begin{pmatrix} (c+m+i-j+1)_{j-1}\,(b-i+j+1)_{a+m-j}&
\text{\scriptsize $1\le i\le a$}\\
\(\frac {c-a} {2}+i-j+1\)_{j-1}\(\frac {b+a}
{2}-i+j+1\)_{a+m-j}&\text{\scriptsize $a<i\le a+m$}
\end{pmatrix}\\
=\frac { \h(a+m) \h(\frac {a-1} {2})\h(\frac {a+1}
{2})\h(\frac {m} {2})^2} 
{\h(\frac {a+m-1} {2})\h(\frac {a+m+1} {2})
\,2^{ {m(a+m-1)}/ {2}}}\hskip8cm\\
\times
\prod _{k=1} ^{ {m}/ {2}}{\(\tfrac {b-1} {2}+k\)_{ {(a+1)}/ {2}}
\(\tfrac {b+1} {2}+k\)_{ {(a-1)}/ {2}}
\(\tfrac {c-1} {2}+k\)_{ {(a+1)}/ {2}}
\(\tfrac {c+1} {2}+k\)_{ {(a-1)}/ {2}}}\\
\times
\prod _{k=0} ^{ {(a-1)}/ {2}-1}{(b+c+m+2k+1)^{a-2k-1}}
\prod _{k=1} ^{ {(a-1)}/ {2}}{(b+c+2m+2k)^{a-2k}}\\
\times\prod _{k= {m}/ {2}+1} ^{m}{(b+c+2k)^{a+m-k}}
\prod _{k=1} ^{ {m}/ {2}}{(b+c+2k)^{m-k}}.
\end{multline}
\end{lemma}

\begin{proof}
We proceed analogously to the proof of Lemma~\ref{lem:amgerade}.
The only difference is the parity of $a$, so we have to read through the
proof of Lemma~\ref{lem:amgerade} and find the places where we used
the fact that $a$ is even. 

As it turns out, the arguments in Steps~1--5 in the proof of
Lemma~\ref{lem:amgerade} can be used here, practically without
change, to establish that the right-hand side of \eqref{eq:amung}
divides the determinant on the left-hand side of \eqref{eq:amung} as
a polynomial in $b$ and $c$.
Differences arise only in the products corresponding to each subcase
(for example, the arguments in Step~3 of the proof of
Lemma~\ref{lem:amgerade} prove that $\prod _{k=1} ^{ {a}/
{2}-1}{(b+c+2m+2k)^{a-2k}}$\linebreak  
divides the determinant $D_1(b,c)$ if $a$ is even,
while for odd $a$ they prove that\linebreak  
$\prod _{k=1} ^{ {(a-1)}/ {2}}{(b+c+2m+2k)^{a-2k}}$
divides $D_1(b,c)$), and in the
fact that in Step~2 
we are now interested in the factors corresponding to odd
values of $e$, $1\le e \le a-2$ (because here the factors with even $e$ are 
covered by Steps~3 and 4). 

Also Step~6, the determination of a degree bound on the determinant, 
can be used verbatim.

\smallskip
\newbox\bcubox
\setbox\bcubox\hbox{\small
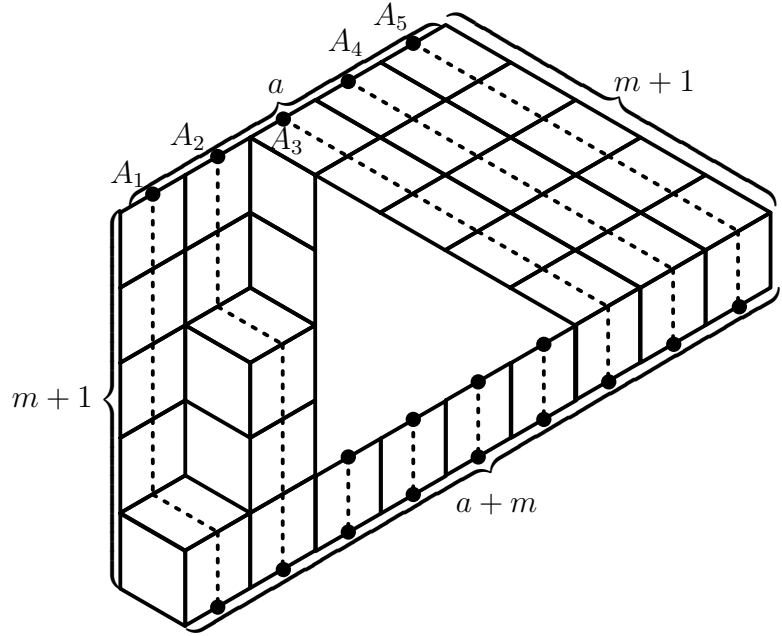
A lozenge tiling and the corresponding path
family for $b=c=1$, $a=5$}
\begin{figure} 
\centertexdraw{
\drawdim truecm 
\linewd.05
\RhombusA \RhombusA \RhombusA \RhombusA \RhombusA \RhombusB
\move(-.866025 -.5) 
\RhombusA \RhombusA \RhombusA \RhombusA \RhombusA \RhombusB
\move(-1.73205 -1)
\RhombusA \RhombusA \RhombusA \RhombusA \RhombusA \RhombusB 
\move(-2.598 -1.5)
\RhombusB \RhombusB \RhombusA \RhombusB \RhombusB \RhombusB
\move(-3.464 -2)
\RhombusB \RhombusB \RhombusB \RhombusB \RhombusA \RhombusB
\move(-1.73205 -1) \RhombusC
\move(-1.73205 -2) \RhombusC
\move(-2.598 -3.5) \RhombusC
\move(-3.4641 -6) \RhombusC
\move(1.73205 -4) \RhombusB
\move(.866025 -4.5) \RhombusB
\move(0 -5) \RhombusB
\move(-.866025 -5.5) \RhombusB
\ringerl(.43301 .25) 
\hdSchritt \hdSchritt \hdSchritt \hdSchritt \hdSchritt \vdSchritt
\ringerl(-.43301 -.25)
\hdSchritt \hdSchritt \hdSchritt \hdSchritt \hdSchritt \vdSchritt
\ringerl(-1.299 -.75)
\hdSchritt\hdSchritt\hdSchritt\hdSchritt\hdSchritt\vdSchritt
\ringerl(-2.165 -1.25)
\vdSchritt\vdSchritt\hdSchritt\vdSchritt\vdSchritt\vdSchritt
\ringerl(-3.031 -1.75)
\vdSchritt\vdSchritt\vdSchritt\vdSchritt\hdSchritt\vdSchritt
\ringerl(4.7631 -3.25)
\ringerl(3.897 -3.75)
\ringerl(3.031 -4.25)

\ringerl(2.165 -4.75) \ringerl(2.165 -3.75) \vdSchritt
\ringerl(1.299 -5.25) \ringerl(1.299 -4.25) \vdSchritt
\ringerl(.43301 -5.75) \ringerl(.43301 -4.75) \vdSchritt
\ringerl(-.43301 -6.25) \ringerl(-.43301 -5.25) \vdSchritt

\ringerl(-1.299 -6.75)
\ringerl(-2.165 -7.25)

\rtext td:-60 (1.2 -5.2) {$\left. \vbox{\vskip4.7cm} \right\}$}
\rtext td:-60 (-1.5 -.5){$\left\{ \vbox{\vskip2.5cm} \right.$ }

\rtext td:60 (3.1 -.8) {$\left. \vbox{\vskip2.7cm} \right \}$}
\htext (-4.9 -7) {$m+1 \left\{ \vbox{\vskip2.7cm} \right.$} 

\htext(-1.5 -.4){$a$}
\htext(1 -6){$a+m$}
\htext(3.1 -.4){$m+1$}

\htext(-3.6 -1.65){$A_1$}
\htext(-2.8 -1.15){$A_2$}
\htext(-1.5 -1.2){$A_3$}
\htext(-.7 .1){$A_4$}
\htext(-.1 .45){$A_5$}

}
\caption{\unhbox\bcubox}
\label{bc1} 
\end{figure}

For the determination of the multiplicative constant 
relating the right-hand and the left-hand side of \eqref{eq:amung}, we
have to modify however the arguments in Step~7 of the proof of
Lemma~\ref{lem:amgerade}. We determine the constant by computing the
determinant for $b=c=1$. Again, this value is most conveniently found
by going back, via \eqref{det} and Lemma~\ref{lem:nonint}, 
to the combinatorial root of the determinant, which is
enumeration of lozenge tilings.

We claim that 
the number of lozenge tilings for $b=c=1$, $a$ odd and $m$ even, equals 
\begin{equation} \label{bc1eq}
2\binom{m+1+\frac {a-1} {2}}{\frac {a-1} {2}}.
\end{equation}
This can be read off Figure~\ref{bc1}, which shows a typical example
of the case $b=c=1$: The path starting at $A_{\frac
{a+1}{2}}$ (see the labeling in Figure~\ref{bc1}; it is derived from
the labeling of starting points of paths in Figure~\ref{tiling}) 
must pass either to the right or to the left of the
triangle. Since the hexagon is symmetric, we can count those path
families where the path
passes to the right, and in the end multiply the resulting number by
two. For those path families, 
the paths starting at points to the right of $A_{\frac
{a+1} {2}}$ are fixed. The paths to the left have all exactly one
South-East step. Suppose that the  South-East step of the path which
starts in $A_i$, $1\le i\le (a-1)/2$, 
occurs as the $h_i$th step. Then we must have
$$m+2\ge h_1\ge h_2\ge\dots \ge h_{\frac {a-1} {2}}\ge 1.$$
So we just have to count monotonously decreasing
sequences of $\frac {a-1} {2}$ numbers between 1 and $m+2$. 
The number is
exactly the binomial coefficient in \eqref{bc1eq}.
It is then a
routine task to check that, on substitution in \eqref{det},
the result agrees exactly with the
right-hand side of \eqref{eq:amung} for $b=c=1$.
\end{proof}

\begin{lemma} \label{lem:-agerade}
Let $a$ and $m$ be nonnegative integers,
$a$ even and $m$ odd. Then 
\begin{multline} \label{eq:-agerade}
\det_{1\le i,j\le a+m} \begin{pmatrix} (c+m+i-j+1)_{j-1}\,(b-i+j+1)_{a+m-j}&
\text{\scriptsize $1\le i\le a$}\\
\(\frac {c-a} {2}+i-j+1\)_{j-1}\(\frac {b+a}
{2}-i+j+1\)_{a+m-j}&\text{\scriptsize $a<i\le a+m$}
\end{pmatrix}\\
=(-1)^{{a} /{2}}\frac {\h(a+m)\h(\frac {a} {2})^2\h(\frac {m-1}
{2})\h(\frac {m+1} {2})} {\h(\frac {a+m-1} {2})\h(\frac
{a+m+1} {2})\,2^{ {m(a+m-1)} /{2}}}\hskip7.8cm\\
\times
\({\tfrac b 2} + \tfrac {1 + m} 2\) _{\frac a 2} \,
\(\tfrac c 2 + {\tfrac {1 + m} 2}\) _{\frac a 2} \,
\prod_{k=1}^{ {(m-1)}/ 2}
     {\( \tfrac b 2 + k\) _{\frac a 2} ^2} 
 \prod_{k=1}^{ {(m-1)}/ 2}
     {\({  {\tfrac c 2} + k}\) _{\frac a 2}^2}
\prod_{k=1}^{ a/ 2-1}
     {\( b + c + 2\,k + m \) }^{a - 2\,k} \\
\times
\prod_{k=1}^{{ a/ 2}-1}
     {{\( b + c + 2\,k + 2\,m \) }^{a - 2\,k}} 
\prod_{k=0}^{{ {(m-1)}/ 2}}
     {{\( 1 + b + c + 2\,k + m \) }^a} 
\prod_{k=1}^{m}{{\( b + c + 2\,k \) }^{m-k}}. 
\end{multline}
\end{lemma}

\begin{proof}
We proceed analogously to the proof of Lemma~\ref{lem:amgerade}.
The only difference is the parity of $m$, so 
we have to check the places in the
proof of Lemma~\ref{lem:amgerade} where we used
the fact that $m$ is even. 

Again, Steps~1--6 can be reused verbatim, except that the
products corresponding to the individual subcases are slightly
different, and 
in Step~2 we are now interested in the factors corresponding to odd
values of $e$, $1\le e \le a-2$ (because the factors with even $e$ are 
covered by Steps~3 and 4). 

The computation of the multiplicative constant 
relating the right-hand and the left-hand side of \eqref{eq:-agerade}
is done analogously to Step~7 in the proof of
Lemma~\ref{lem:amgerade}. I.e., we compute the determinant for $b=c=0$
by going back, via \eqref{det} and Lemma~\ref{lem:nonint}, to the
lozenge tiling interpretation of the determinant. We already
concluded in the proof of Lemma~\ref{lem:amgerade} that for $b=c=0$
there is just one lozenge tiling (see Figure~\ref{bc0}). By
definition, the statistic $n(.)$ attains the value $a/2$ on this
lozenge tiling, so that its weight is $(-1)^{a/2}$. It is then not
difficult to verify that, on substitution of this in \eqref{det}, 
the result agrees exactly with the
right-hand side of \eqref{eq:-agerade} for $b=c=0$.
\end{proof}

\begin{lemma}\label{lem:-aung}
Let $a$ and $m$ be odd nonnegative integers. Then 
\begin{equation} \label{eq:-aung}
\det_{1\le i,j\le a+m} \begin{pmatrix} (c+m+i-j+1)_{j-1}\,(b-i+j+1)_{a+m-j}&
\text{\scriptsize $1\le i\le a$}\\
\(\frac {c-a} {2}+i-j+1\)_{j-1}\(\frac {b+a}
{2}-i+j+1\)_{a+m-j}&\text{\scriptsize $a<i\le a+m$}
\end{pmatrix}=0.
\end{equation}
\end{lemma}

\begin{proof}
Analogously to the previous cases, we can show that the product
\begin{multline*}
\prod _{i=1} ^{ {(m+1)}/ {2}}{(\tfrac{b-1}2+i)_{ {(a+1)}/2}}
\prod _{i=1} ^{ {(m-1)}/ {2}}{(\tfrac{b+1}2+i)_{ {(a-1)}/2}}
\prod _{k=1} ^{ {(a-1)}/2}{(b+c+m+2k)^{a-2k}}\\
\times \prod _{k=1} ^{ {(a-1)}/ {2}}{(b+c+2m+2k)^{a-2k}}
\prod _{k= {(m+1)}/ {2}} ^{m}{(b+c+2k)^{a+m-k}}
\prod _{k=1} ^{ {(m-1)}/ {2}}{(b+c+2k)^{m-k}}
\end{multline*}
divides the determinant as a polynomial in $b$ and $c$. 
Although not completely obvious, this is implied by the linear
combinations of Lem\-ma~\ref{lem:amgerade}, Steps~1--5.
The degree in $b$ of this product is $\binom{a+m}2 +1$ which 
is larger
than the maximal degree $\binom {a+m}2$ of the determinant viewed as a
polynomial in $b$.
So the determinant must be zero.
\end{proof}
\end{section}

\begin{section}{Determinant Evaluations, II} \label{sec:det2}
In this section we evaluate the determinant in
Lemma~\ref{lem:nonint2}.
We proceed analogously to Section~\ref{sec:det} and start by taking
${(b+c+m)!} \big/{(b+a+m-i)!\,(c+m+i-1)!}$
out of the $i$th row, $i=1,2,\dots,a$, and 
$ {\(\frac {b+c} {2}\)!}\big/{\(\frac {b+3a+1}
{2}+m-i\)!\(\frac {c-a-1} {2}+i-1\)!}$
out of the $i$th row, $i=a+1,a+2,\dots,a+m$.
This gives
\begin{multline} \label{det2}
\det_{1\le i,j\le a+m} \begin{pmatrix} \dbinom{b+c+m}{b-i+j}&
\text{\scriptsize $1\le i\le a$}\\
\dbinom{\frac {b+c} {2}}{\frac {b+a+1} {2}-i+j}&
\text{\scriptsize $a+1\le i \le a+m$}
\end{pmatrix}\\
=\prod _{i=1} ^{a}{\frac {(b+c+m)!} {(b+a+m-i)!\,(c+m+i-1)!}}
\prod _{i=a+1} ^{a+m}{\frac {\(\frac {b+c} {2}\)!}{\(\frac {b+3a+1}
{2}+m-i\)!\(\frac {c-a-1} {2}+i-1\)!}}\\
\times\det_{1\le i,j\le a+m} \begin{pmatrix} (c+m+i-j+1)_{j-1}\,(b-i+j+1)_{a+m-j}&
\text{\scriptsize $1\le i\le a$}\\
\(\frac {c-a-1} {2}+i-j+1\)_{j-1}\(\frac {b+a+1}{2}-i+j+1\)_{a+m-j}&
\text{\scriptsize $a<i\le a+m$}
\end{pmatrix}.
\end{multline}

Thus it suffices to evaluate the determinant on the right-hand side.
As in the preceding section,
the advantage is that this determinant is a polynomial in $b$ and
$c$. So we can again apply the ``identification of factors"
method, as proposed in \cite[Sec.~2.4]{KratBN}.
We note that the first $a$ rows of the matrix are identical to those
of \eqref{det}, whereas the other $m$ rows differ only slightly.
Hence we can use many arguments from
Section~\ref{sec:det}.
The four lemmas below address the four different cases, as $a$ and $m$
vary through all combinations of parities.

\begin{lemma} \label{lem:shamgerade}
Let $a$ and $m$ be both even nonnegative integers. Then 
\begin{multline} \label{eq:shamgerade}
\det_{1\le i,j\le a+m} \begin{pmatrix} (c+m+i-j+1)_{j-1}\,(b-i+j+1)_{a+m-j}&
\text{\scriptsize $1\le i\le a$}\\
\(\frac {c-a-1} {2}+i-j+1\)_{j-1}\(\frac {b+a+1}{2}-i+j+1\)_{a+m-j}&
\text{\scriptsize $a<i\le a+m$}
\end{pmatrix}\\
=\frac {\h(a+m)\h(\frac {a} {2})^2\h(\frac {m} {2})^2} 
{\h(\frac {a+m} {2})^2 \,2^{ {m(a+m-1)} /{2}}}
\prod _{k=1} ^{ {m}/ {2}}{(\tfrac{b-1}2+k)_{ a/ 2}\,
(\tfrac{b+1}2+k)_{ {a}/ {2}}}\hskip5.8cm\\
\times\prod _{k=1} ^{ {m}/ {2}}{(\tfrac{c-1}2+k)_{ a/ 2}\,
(\tfrac{c+1}2+k)_{ {a}/ {2}}}
\prod _{k=0} ^{ {a}/ {2}-1}{(b+c+m+2k+1)^{a-2k-1}}\\
\times\prod _{k=1} ^{ {a}/ {2}-1}{(b+c+2m+2k)^{a-2k}}
\prod _{k= {m}/ {2}+1} ^{m}{(b+c+2k)^{a+m-k}}
\prod _{k=1} ^{ {m}/ {2}}{(b+c+2k)^{m-k}}.
\end{multline}
\end{lemma}

\begin{proof}[Proof of Lemma~\ref{lem:shamgerade}]
Let us denote the determinant in \eqref{eq:shamgerade} by $D_2(b,c)$.
We will again proceed in the spirit of the proof of
Lemma~\ref{lem:amgerade}. I.e., we first show,
in Steps~1--5 below, that the right-hand side of
\eqref{eq:shamgerade} divides $D_2(b,c)$ as a polynomial in $b$ and
$c$. Then, in Step~6, we show that the degree of $D_2(b,c)$ as a polynomial
in $b$ is at most $\binom {a+m}2$, the same being true for the degree
in $c$. Analogously to the proof of
Lemma~\ref{lem:amgerade}, we conclude that $D_2(b,c)$ must
equal the right-hand side of \eqref{eq:shamgerade}, times a
constant with respect to $b$ and $c$. 
That this constant is equal to 1 is finally shown in
Step~7, by evaluating the determinant $D_2(b,c)$ for $b=c=1$.

In order to prove (in Steps~1--5) that the right-hand side of
\eqref{eq:shamgerade} divides $D_2(b,c)$, for each linear factor of
\eqref{eq:shamgerade} we exhibit again sufficiently many linear
combinations of columns or rows which vanish. These linear
combinations are almost identical (sometimes they are even identical)
with the corresponding linear combinations in the proof of
Lemma~\ref{lem:amgerade}. Consequently, we will merely state these linear
combinations here, but will not bother to supply their verifications,
because these parallel the verifications
in the proof of Lemma~\ref{lem:amgerade}.

{\em Step~1. $\prod _{k=1} ^{ {m}/ {2}}
{(\tfrac{b-1}2+k)_{ a/ 2}\,
(\tfrac{b+1}2+k)_{ {a}/ {2}}}\,{(\tfrac{c-1}2+k)_{ a/ 2}\,
(\tfrac{c+1}2+k)_{ {a}/ {2}}}$ divides the determinant.}
Unlike in the case of 
the previous determinant $D_1(b,c)$ (see \eqref{eq:amgerade}),
here it is not possible to infer symmetry of $D_2(b,c)$ in $b$ and $c$
directly from the definition. Therefore it will be 
necessary to prove {\em separately} that the factors involving $b$,
respectively $c$, divide the
determinant.

Again, we distinguish between four subcases, labeled below as
(a), (b), (c), and (d).

\smallskip
(a) $(b+e)^e\,(c+e)^e$ divides $D_2(b,c)$ for $1\le e
\le \min\{a,m\}$, $e \not\equiv a$ mod 2:
This follows from the easily verified fact that
$(b+e)$ is a factor of each entry in the first $e$ columns of
$D_2(b,c)$, respectively, that $(c+e)$ is a factor of each entry in the
last $e$ columns of $D_2(b,c)$.

\smallskip
(b) $(b+e)^m\,(c+e)^m$ divides $D_2(b,c)$ for $m<e<a$, $e\not\equiv a$ mod 2:
The following equations hold for $s=1,2,\dots,m$:
\begin{equation} \label{eq:shlin1b}
\sum _{j=1} ^{e+s-m}{\binom{e-m+s-1}{j-1}\frac
{(c+a-e-s+2m+1)_{e+s-j-m}} {(a-e-s+2m+1)_{e+s-j-m}}\cdot (\text
{column $j$ of $D_2(-e,c)$})}=0,
\end{equation}
and 
\begin{multline} \label{eq:shlin1bc}
\sum _{j=1} ^{e+s-m}\binom{e-m+s-1}{j-1}\frac
{(b+a-e-s+2m+1)_{e+s-j-m}} {(a-e-s+2m+1)_{e+s-j-m}}\\
\cdot (\text
{column $(a+m+1-j)$ of $D_2(b,-e)$})=0.
\end{multline}
%This is exactly the same linear combination as in Step~1 of the proof
%of Lemma~\ref{lem:amgerade}, so we already know that the equation holds for
%$i\le a$. It is important to know that we did not need any
%information about the parity of $e$ in this subcase.
%For $i>a$ the hypergeometric form of the sum in \eqref{eq:shlin1b} is
%\begin{multline*} 
%( \textstyle {\frac 3 2} + {\frac {3\,a} 2} -
%{\frac e 2} - i + m) _{-1 + e - m + s} \,
%     ({ \textstyle 1 + a + c - e + 2\,m - s}) _{-1 + e - m + s}\\
%\times
%\frac{ 
%({ \textstyle {\frac 3 2} + {\frac a 2} + {\frac e 2} - i - m + s}) 
%_{a - e + 2\,m - s} } 
%   {({ \textstyle 1 + a - e + 2\,m - s}) _{-1 + e - m + s} }\\
%\times {} _{3} F _{2} \!\left [ \begin{matrix} { 1 - e + m - s, 1 - a - m,
%-{\frac 1  2} - {\frac {3\,a} 2} - {\frac c 2} + i - m}\\ { 1 - a - c - m,
%{\frac 1 2}
% - {\frac {3\,a} 2} - {\frac e 2} + i - s}\end{matrix} ; {\displaystyle
%      1}\right ] .
%\end{multline*}
%The factor $({ \textstyle {\frac 3 2} + {\frac a 2} + {\frac e 2} - i - m + s}) _{a - e
%      + 2\,m - s}$ vanishes for $i>a$ and the denominator is never
%zero, so the sum~\eqref{eq:shlin1b} equals zero, as desired.

\smallskip
(c) $(b+e)^a$ divides $D_2(b,c)$ for $a<e<m$, $e\not\equiv a$ mod  2:
The following equations hold for $s=1,2,\dots,a$:
\begin{multline} \label{eq:shlin1c}
\sum _{j=1} ^{\frac {e-a-1} {2}+s}\binom{\frac {e-a-1} {2}+s-1}{j-1}
\frac {\(\frac {c} {2}+m-\frac {e} {2}+a-s+1\)_{ {(e-a-1)}/ {2}+s-j}} 
{\(m+\frac {3a-e} {2}-s+1\)_{ {(e-a-1)}/ {2}+s-j}}\\
\cdot
(\text {column $j$ of $D_2(-e,c)$})=0,
\end{multline}
and 
\begin{multline} \label{eq:shlin1cc}
\sum _{j=1} ^{\frac {e-a-1} {2}+s}\binom{\frac {e-a-1} {2}+s-1}{j-1}
\frac {\(\frac {b} {2}+m-\frac {e} {2}+a-s+1\)_{ {(e-a-1)}/ {2}+s-j}} 
{\(m+\frac {3a-e} {2}-s+1\)_{ {(e-a-1)}/ {2}+s-j}}\\
\cdot
(\text {column $(a+m+1-j)$ of $D_2(b,-e)$})=0.
\end{multline}

\smallskip
(d) $(b+e)^{a+m-e}$ divides $D_2(b,c)$ for $\max\{a,m\} \le e \le a+m$,
$e \not\equiv a$ mod 2:
The following equations hold for $s=1,2,\dots, a+m-e$:
\begin{multline}\label{eq:shlin1d}
\sum _{i=1} ^{s}\binom{s-1}{i-1}
\frac {(-1)^{i}\(\frac {c-e} {2}+1\)_{a+m-s}\(\frac {c-e}
{2}+m\)_{i-1}} {(1+c-e+m)_{a+m-s+i-1}} \\
\cdot(\text {row $(a+m-e-s+i)$
of $D_2(-e,c)$})\\
+(\text {row $\(m+\tfrac {3a+1} {2}-\tfrac {e} {2}-s+1\)$ of
$D_2(-e,c)$})=0,
\end{multline}
and 
\begin{multline}\label{eq:shlin1dc}
\sum _{i=1} ^{s}\binom{s-1}{i-1}
\frac {(-1)^{i}\(\frac {b-e} {2}+1\)_{a+m-s}\(\frac {b-e}
{2}+m\)_{i-1}} {(1+b-e+m)_{a+m-s+i-1}} \\
\cdot(\text {row $(e+s-m-i+1)$
of $D_2(b,-e)$})\\
+(\text {row $\(\tfrac {a+1} {2}+\tfrac {e} {2}+s\)$ of
$D_2(b,-e)$})=0.
\end{multline}

%The sum contains only row indices $\le a$. Therefore, we can use the
%results of Step~1(d) of the proof of Lemma~\ref{lem:amgerade}. 
%For $j\le a+m-s$ the sum
%is zero and so is the extra term because it contains the factor $(\frac
%{-e+a+1} {2} -i+j+1)_{a+m-j}$. If $j > a+m-s$  the sum equals the
%expression \eqref{resd} and it is again a simple calculation that
%this cancels with the extra term.
%The reason for this is that our new extra row index differs from the
%old one by $\frac {1} {2} $ and this cancels exactly with the changes
%of the matrix entries.

\smallskip
{\em Step~2. $\prod _{k=0} ^{ {a}/ {2}-1}{(b+c+m+2k+1)^{a-2k-1}}$
divides the determinant.}
The following
equation holds for $0\le e \le a-2$, $s=1,2,\dots,e+1$:
\begin{multline} \label{eq:shlin2}
\sum _{i=1} ^{a-e-1}{\(\frac {(c+m+i)_{a-e-i+s-1}}
{(c-e-1+i)_{a-e-i+s-1}}\binom{a-e-2}{i-1}\frac {(s)_{a-e-1}\,(-1)^i}
{(s-i+a-e-1)(a-e-2)!}\right.}\\
\cdot (\text {row $i$ of $D_2(-c-a-m+1+e,c)$})\bigg)\\
+(-1)^{a-e-1} \cdot(\text {row $(a-e-1+s)$ of $D_2(-c-a-m+1+e,c)$})=0.
\end{multline}
Here, we need equation \eqref{eq:shlin2} only for even $e$.
%Here, only rows with index $\le a$ are involved, so this follows
%directly from Step~2 of the proof of Lemma~\ref{lem:amgerade}.

\smallskip
{\em Step~3. $\prod _{k=1} ^{ {a}/ {2}-1}{(b+c+2m+2k)^{a-2k}}$ 
divides the determinant.}
The following equation holds for $0<e\le a$, $e\equiv a$ mod 2, 
and $s=1,2,\dots,e$:
\begin{equation} \label{eq:shlin3}
\sum _{j=s} ^{a+m+s-e}{\binom{a+m-e}{j-s}\cdot (\text {column $j$ of
$D_2(-c-2m-a+e,c)$})}=0.
\end{equation}
%For $i\le a$ this follows directly from Step~3 of the proof of
%Lemma~\ref{lem:amgerade}. 
%If $i>a$ we convert the sum into hypergeometric form and get
%\begin{multline*}
%({ \textstyle {\frac 1 2} - {\frac a 2} + {\frac c 2} + i - s}) _{-1 + s} \,
%  ({ \textstyle {\frac 3 2} - {\frac c 2} + {\frac e 2} - i - m + s}) _{a + m -
%   s} \\
%\times
%{} _{2} F _{1} \!\left [ \begin{matrix} { {\frac 1 2} + {\frac a 2} -
%{\frac c 2}-i+s,-a+e-m}\\ { {\frac 3 2}- {\frac c 2}+ {\frac e 2} - i - m +
%   s}\end{matrix} ; {\displaystyle 1}\right ] .
%\end{multline*}
%The Chu--Vandermonde summation formula \eqref{vand} gives
%$$\frac {({ \textstyle 1 - {\frac a 2} + {\frac e 2} - m}) _{a - e + m} \,
%     ({ \textstyle {\frac 1 2} - {\frac a 2} + {\frac c 2} + i - s}) _{-1 + s} }
%    {({ \textstyle {\frac 3 2} + a - {\frac c 2} + {\frac e 2} - i}) _{-e +
%     s} }.
%$$
%This expression is zero because the factor 
%$(1-\frac {a} {2}+\frac {e} {2}-m)_{a-e+m}$ vanishes
%for $e\equiv a$ mod $2$, so the sum in \eqref{eq:lin3} is zero, as
%desired.

\smallskip
{\em Step~4. $\prod _{k= {m}/ {2}+1} ^{m}{(b+c+2k)^{a+m-k}}$ 
divides the determinant.}
The following equation holds for $m/2<e\le m$ 
and $s=1,2,\dots, a+m-e$:
\begin{equation} \label{eq:shlin4}
\sum _{j=s} ^{s+e}{\binom{e}{j-s}\cdot
(\text {column $j$ of $D_2(-c-2e,c)$})}=0.
\end{equation}
%For $i\le a$ this is analogous to Step~4 of the proof of
%Lemma~\ref{lem:amgerade}.
%For $i>a$ conversion of the sum to hypergeometric form gives
%\begin{multline*}
%  ({ \textstyle {\frac 1 2} - {\frac a 2} + {\frac c 2} + i - s}) _{-1 + s} \,
%  ({ \textstyle {\frac 3 2} + {\frac a 2} - {\frac c 2} - {e} - i + s})
%   _{a + m - s} \\
%\times
%{} _{2} F _{1} \!\left [ \begin{matrix} { {\frac 1 2} + {\frac a 2} -
%{\frac c 2}
%   - i + s, {{-e}}}\\ { {\frac 3 2} + {\frac a 2} - {\frac c 2} -
%{e} - i + s}\end{matrix} ; {\displaystyle 1}\right ] .
%\end{multline*}
%
%Application of the Chu--Vandermonde summation formula \eqref{vand}
%gives
%$$\frac{({ \textstyle 1 - {e}}) _{{e}} \,
% ({ \textstyle {\frac 1 2} - {\frac a 2} + {\frac c 2} + i - s}) _{-1 + s} }
% {({ \textstyle {\frac 3 2} + {\frac {3\,a} 2} - {\frac c 2} -
%{e} - i + m}) _{-a + {e} - m + s} }.
%$$
%This expression is zero because the factor 
%$(1-e)_{ {e}}$ vanishes, so the sum in
%\eqref{eq:shlin4} is zero, as desired.

\smallskip
{\em Step 5. $\prod _{k=1} ^{ {m}/ {2}}{(b+c+2k)^{m-k}}
$ divides the determinant.}
The following equation holds for
$e\le m-1$ and $s=1,2,\dots,e$:
\begin{multline} \label{eq:shlin5}
\sum _{i=1} ^{m-s+1}(-1)^i\binom{m-s}{i-1}\frac {\(\frac {c} {2}+\frac
{a} {2}+i-\frac {1} {2}\)_{m-s-i+1}} {\(\frac {c} {2}-\frac {a}
{2}-e+i-\frac {1} {2}\)_{m-s-i+1}}\\
\cdot (\text {row $(a+i)$ of $D_1(-c-2m+2e,c)$})=0.
\end{multline}

%Conversion to hypergeometric form of the $j$th coordinate gives
%\begin{multline*}
%-  \frac{( \textstyle {\frac 1 2} + {\frac a 2} + {\frac c 2}) _{m - s} \,
%      ({ \textstyle {\frac 3 2} + {\frac a 2} + {\frac c 2} - j}) _{-1 + j} \,
%      ({ \textstyle {\frac 1 2} - {\frac a 2} - {\frac c 2} + e + j - m}) _{a - j
%       + m}  }
%    {({ \textstyle {\frac 1 2} - {\frac a 2} + {\frac c 2} - e}) _{m - s}
%}\\
%{} _{2} F _{1} \!\left [ \begin{matrix} { {\frac 1 2} + {\frac a 2} +
%{\frac c 2} - e - j + m, -m + s}\\ { {\frac 3 2} + {\frac a 2} +
%{\frac c 2} - j}\end{matrix} ; {\displaystyle 1}\right ] .
%\end{multline*}
%We apply the Chu--Vandermonde summation formula \eqref{vand} and
%get
%\begin{multline*}
%-({ \textstyle {\frac 1 2} + {\frac a 2} + {\frac c 2}}) _{m - s} \,
%      ({ \textstyle 1 + e - m}) _{m - s}\\
%\times
%\frac {
%      ({ \textstyle {\frac 1 2} - {\frac a 2} - {\frac c 2} + e + j - m}) _{a - j
%       + m} \,({ \textstyle {\frac 3 2} + {\frac a 2} + {\frac c 2} - j + m - s})
%       _{-1 + j - m + s} } 
%    {({ \textstyle {\frac 1 2} - {\frac a 2} + {\frac c 2} - e}) _{m - s}
%}.
%\end{multline*}

\smallskip
{\em Step~6. Determination of the degree of $D_2(b,c)$ as a polynomial
in $b$.} This is clearly the same degree as for $D_1(b,c)$, that is, 
$\binom {a+m}2$.

\smallskip
{\em Step~7. Computation of the multiplicative constant.}
In analogy to the proof of Lemma~\ref{lem:amung},
we evaluate the determinant 
for $b=c=1$. Again, we do this by going back, via
\eqref{det2} and Lemma~\ref{lem:nonint2}, to the
combinatorial origin of the determinant, which is enumeration of
lozenge tilings. We can still use
Figure~\ref{bc1} for our considerations. 
The number of lozenge tilings is easily seen to be equal to
$\binom{m+1+\frac {a} {2}}{\frac a2}+\binom{m+1+\frac
a2-1}{\frac a2 -1}$.
It is then 
a routine computation to verify that this does indeed give the
multiplicative constant as claimed in \eqref{eq:shamgerade}.
\end{proof}

\begin{lemma} \label{lem:shaome}
Let $a$ and $m$ be nonnegative integers,
$a$ odd and $m$ even. Then 
\begin{multline} \label{eq:shaome}
\det_{1\le i,j\le a+m} \begin{pmatrix} (c+m+i-j+1)_{j-1}\,(b-i+j+1)_{a+m-j}&
\text{\scriptsize $1\le i\le a$}\\
\(\frac {c-a-1} {2}+i-j+1\)_{j-1}\(\frac {b+a+1}{2}-i+j+1\)_{a+m-j}&
\text{\scriptsize $a<i\le a+m$}
\end{pmatrix}\\
=\frac {\h(a+m)\h(\frac {a-1} {2})\h(\frac {a+1} {2})\h(\frac {m}
{2})^2} {\h(\frac {a+m-1} {2})\h(\frac {a+m+1}
{2})\,2^{ {m(a+m-1)}/ {2}}}\hskip8cm\\
\times\prod _{k=1} ^{ {m}/ {2}}{(\tfrac{b}2+k)_{ {(a-1)}/ {2}} \,
(\tfrac{b}2+k)_{ {(a+1)}/ {2}}\,
(\tfrac{c}2+k)_{ {(a-1)}/ {2}}\,
(\tfrac{c}2+k)_{ {(a+1)}/ {2}}}\\
\times
\prod _{k=0} ^{ {(a-3)}/ {2}}{(b+c+m+2k+1)^{a-2k-1}}
\prod _{k=0} ^{ {(a-1)}/ {2}}{(b+c+2m+2k)^{a-2k}}\\
\times
\prod _{k= {m}/ {2}+1} ^{m}{(b+c+2k)^{a+m-k}}
\prod _{k=1} ^{ {m}/ {2}}{(b+c+2k)^{m-k}}.
\end{multline}
\end{lemma}

\begin{proof}
We proceed analogously to the proof of Lemma~\ref{lem:shamgerade}.
The only difference is the parity of $a$, so we have to check the
places in the proof of Lemma~\ref{lem:shamgerade} where we used
the fact that $a$ is even. 

Steps~1, 3--5 can be reused verbatim, but
the corresponding products are slightly different.

In Step~2 we are now interested in the factors corresponding to odd
values of $e$ ($1\le e \le a-2$), because the factors with even $e$ are 
covered by Steps~3 and 4. 

Step~6 can be reused verbatim.

The computation of the multiplicative constant is done analogously to
Step~7 in the proof of Lemma~\ref{lem:amgerade}.
Again using Figure~\ref{bc0}, we see that the number of lozenge
tilings, related to our determinant via \eqref{det2} and
Lemma~\ref{lem:nonint2}, for $b=c=0$ equals 1.
It is then a routine computation to verify that this gives the
multiplicative constant as claimed in \eqref{eq:shaome}.
\end{proof}

\begin{lemma} \label{lem:shaemo}
Let $a$ and $m$ be nonnegative integers,
$a$ even and $m$ odd. Then 
\begin{multline} \label{eq:shaemo}
\det_{1\le i,j\le a+m} \begin{pmatrix} (c+m+i-j+1)_{j-1}\,(b-i+j+1)_{a+m-j}&
\text{\scriptsize $1\le i\le a$}\\
\(\frac {c-a-1} {2}+i-j+1\)_{j-1}\(\frac {b+a+1}{2}-i+j+1\)_{a+m-j}&
\text{\scriptsize $a<i\le a+m$}
\end{pmatrix}\\
=(-1)^{a/2}\frac {\h(a+m)\h(\frac {a}
{2})^2\h(\frac {m-1} {2})\h(\frac {m+1} {2})} 
{\h(\frac {a+m-1} {2})\h(\frac {a+m+1}
{2})\,2^{ {m(a+m-1)} /{2}}}
\prod _{k=1} ^{ {(m+1)}/ {2}}{(\tfrac{b-1}2+k)_{ {a}/ {2}}}
\prod _{k=1} ^{ {(m-1)}/ {2}}{(\tfrac{b+1}2+k)_{ {a}/ {2}}}
\hskip1cm\\
\times
\prod _{k=1} ^{ {(m+1)}/ {2}}{(\tfrac{c-1}2+k)_{ {a}/ {2}}}
\prod _{k=1} ^{ {(m-1)}/ {2}}{(\tfrac{c+1}2+k)_{ {a}/ {2}}}
\prod _{k=1} ^{ {a}/ {2}-1}{(b+c+m+2k)^{a-2k}}\\
\times
\prod _{k=1} ^{ {a}/ {2}-1}{(b+c+2m+2k)^{a-2k}}
\prod _{k= {(m+1)}/ {2}} ^{m}{(b+c+2k)^{a+m-k}}
\prod _{k=1} ^{ {(m-1)}/ {2}}{(b+c+2k)^{m-k}}.
\end{multline}
\end{lemma}
\begin{proof}
We proceed analogously to the proof of Lemma~\ref{lem:shamgerade}.
The only difference is the parity of $m$, so we have to check the
places in the proof of Lemma~\ref{lem:shamgerade} where we used
the fact that $m$ is even. 

Steps~1, 3--5 can be reused verbatim, but
the corresponding products are slightly different.

In Step~2 we are now interested in the factors corresponding to odd
values of $e$ ($1\le e \le a-3$), because the factors with even $e$ are 
covered by Steps~3 and 4. 

Step~6 can be reused verbatim.

The computation of the multiplicative constant is done analogously to
Step~7 in the proof of Lemma~\ref{lem:shamgerade}.
Using again Figure~\ref{bc1}, we see that the $(-1)$-enumeration 
of lozenge tilings, related to our determinant via \eqref{det2} and
Lemma~\ref{lem:nonint2}, for
$b=c=1$ equals $(-1)^{a/2}\binom{m+1+\frac {a} {2}}{\frac a2}
+(-1)^{a/2+1}\binom{m+1+\frac a2-1}{\frac a2 -1}$.
It is then a routine computation to verify that this gives the
multiplicative constant as claimed in \eqref{eq:shaemo}.
\end{proof}

\begin{lemma} \label{lem:shaomo}
Let $a$ and $m$ be odd nonnegative integers.
Then 
\begin{multline} \label{eq:shaomo}
\det_{1\le i,j\le a+m} \begin{pmatrix} (c+m+i-j+1)_{j-1}\,(b-i+j+1)_{a+m-j}&
\text{\scriptsize $1\le i\le a$}\\
\(\frac {c-a-1} {2}+i-j+1\)_{j-1}\(\frac {b+a+1}{2}-i+j+1\)_{a+m-j}&
\text{\scriptsize $a<i\le a+m$}
\end{pmatrix}\\
=(-1)^{(a+1)/2}\frac {\h(a+m)\h(\frac {a-1} {2})
\h(\frac {a+1} {2})\h(\frac {m-1} {2})\h(\frac {m+1} {2})} 
{\h(\frac {a+m} {2})^2
\,2^{ {m(a+m-1)} /{2}+ {1} /{2}}}\hskip5.6cm\\
\times
\prod _{k=1} ^{ {(m+1)}/ {2}}{(\tfrac{b}2+k)_{ {(a-1)}/ {2}}}
\prod _{k=1} ^{ {(m-1)}/ {2}}{(\tfrac{b}2+k)_{ {(a+1)}/ {2}}}\\
\times
\prod _{k=1} ^{ {(m+1)}/ {2}}{(\tfrac{c}2+k)_{ {(a-1)}/ {2}}}
\prod _{k=1} ^{ {(m-1)}/ {2}}{(\tfrac{c}2+k)_{ {(a+1)}/ {2}}}
\prod _{k=1} ^{ {(a-1)}/ {2}}{(b+c+m+2k)^{a-2k}}\\
\times
\prod _{k=1} ^{ {(a-1)}/ {2}}{(b+c+2m+2k)^{a-2k}}
\prod _{k= {(m+1)}/ {2}} ^{m}{(b+c+2k)^{a+m-k}}
\prod _{k=1} ^{ {(m-1)}/ {2}}{(b+c+2k)^{m-k}}.
\end{multline}
\end{lemma}
\begin{proof}
We proceed analogously to the proof of Lemma~\ref{lem:shamgerade}.
The parameters $a$ and $m$ are odd, so we have to check the places in
the proof of Lemma~\ref{lem:shamgerade} where we used
the fact that $a$ or $m$ is even. 

Steps~1--6 can be reused verbatim, but
the corresponding products are slightly different.

The computation of the multiplicative constant is done analogously to
Step~7 in the proof of Lemma~\ref{lem:amgerade}.
Again using Figure~\ref{bc0}, we see that the $(-1)$-enumeration 
of lozenge tilings, related to our determinant via \eqref{det2} and
Lemma~\ref{lem:nonint2}, for $b=c=0$ equals $(-1)^{(a+1)/ {2}}$.
It is then a routine computation to verify that this gives the
multiplicative constant as claimed in \eqref{eq:shaomo}.
\end{proof}
\end{section}

\begin{section}{Proof of Theorem~\ref{thm:Zare1}}
\label{sec:Zare1}

For the proof of Theorem~\ref{thm:Zare1}, we proceed similarly to \cite{MiRRAD}. 
We define determinants $Z_n(x,\mu)$ by
\begin{equation} \label{eq:Z_n}
Z_n(x,\mu)=\det_{0\le i,j\le n-1}\left(-\delta_{ij}+\sum_{t,k=0}
^{n-1}\binom {i+\mu}t\binom kt\binom
{j-k+\mu-1}{j-k}x^{k-t}\right).
\end{equation}
The only difference to the definition of $Z_n(x,\mu)$ in
\cite{MiRRAD} is the minus sign in front of $\delta_{ij}$.

Then an analogue of Theorem~5 of \cite{MiRRAD} is true.
\begin{lemma} \label{lem:Z_n}
Let $n$ be a nonnegative integer. Then $Z_n(x,\mu)=0$
if $n$ is odd. If $n$ is even, then $Z_n(x,\mu)$ factors,
\begin{multline} \label{eq:Zn-Eval}
Z_n(x,\mu)=(-1)^{n/2}\det_{0\le i,j\le n/2-1}\left(\sum _{t=0} ^{n-1}\frac {t+1}
{j+1}\binom {i+\mu}{t-i}\binom {j+1}{t-j}x^{2j+1-t}\right)\\
\times
\det_{0\le i,j\le n/2-1}\left(\sum _{t=0} ^{n-1}\frac {t+\mu+1}
{i+\mu+1}\binom {i+\mu+1}{t-i}\binom {j}{t-j}x^{2j-t}\right).
\end{multline}
\end{lemma}

\begin{proof}As in the proof of Theorem~5 of
\cite{MiRRAD}, define matrices $S,M,U$,
\begin{multline*} 
S=\(\binom {i+\mu}t\)_{0\le i,t\le n-1},\quad 
M=\(\binom {k}tx^{k-t}\)_{0\le t,k\le n-1},\\
U=\(\binom {j-k+\mu-1}{j-k}\)_{0\le k,j\le n-1},
\end{multline*}
and $J$ and $F(x)$,
\begin{equation} \label{eq:JF}
J=\((-1)^{k-i}\binom {\mu}{k-i}\)_{0\le i,k\le n-1},\quad 
F(x)=\(\binom {j-\fl{j/2}}{j-i}(-x)^{j-i}\)_{0\le i,j\le n-1}.
\end{equation}
Thus, $Z_n(x,\mu)$ equals $\det(-I+SMU)$. 
Now, as in \cite{MiRRAD}, 
multiply $Z_n(x,\mu)$ on the left by $\det(F(1)^t)$ and on the right
by $\det(JF(x))$. Subsequently
do the manipulations given in \cite{MiRRAD} (which amount to applying the
Chu--Vandermonde summation formula several times). The result is that 
$$Z_n(x,\mu)=\det_{0\le i,j\le n-1}(-I+SMU)=\det_{0\le i,j\le
n-1}(-V(x,\mu)+W(x,\mu)),$$
where
\begin{align} \label{eq:vw}
V(x,\mu)_{2i+r,2j+s}&=\sum _{t=0} ^{n-1}(-1)^{r+s}\binom
{i+r+\mu}{t-i}\binom {j+s}{t-j}x^{2j+s-t},\\
W(x,\mu)_{2i+r,2j+s}&=\sum _{t=0} ^{n-1}\binom
{i+\mu}{t-i-r}\binom {j}{t-j-s}x^{2j+s-t},
\end{align}
where $r$ and $s$ are restricted to be $0$ or $1$, as in
\cite{MiRRAD}.

It is straightforward to check that $V_{2i,2j}=W_{2i,2j}$. Hence,
each entry of the
matrix $-V+W$ in an even-numbered row and even-numbered column is 0.
This implies that $\det(-V+W)$ must be 0 whenever the size of the
matrix, $n$, is odd. In the case that $n$ is even it implies
the factorization
\begin{multline*}
Z_n(x,\mu)=\det(-V(x,\mu)+W(x,\mu))\\
=(-1)^{n/2}\det_{0\le i,j\le n/2-1}(-V_{2i,2j+1}+W_{2i,2j+1})
\det_{0\le i,j\le n/2-1}(-V_{2i+1,2j}+W_{2i+1,2j}).
\end{multline*}
As is easily verified, this equation is exactly equivalent to
\eqref{eq:Zn-Eval}.
\end{proof}

\begin{proof}[Proof of Theorem~\ref{thm:Zare1}]
Now choose $x=1$, $\mu=m/2$, $n=a$ in Lemma~\ref{lem:Z_n}. 
Then all the sums appearing in
\eqref{eq:Zn-Eval} can be
evaluated by means of the Chu--Vandermonde summation \eqref{vand}. The result is
\begin{align*} 
Z_a(1,m/2)&=\det_{0\le i,j\le a-1}\left(-\delta_{ij}+\binom {m+i+j}{j}\right)\\
&=\det_{0\le i,j\le a/2-1}\left((3i+m+1)
\frac {(i+j+m/2)!}{(2i-j+m/2)!\,(2j-i+1)!}\right)\\
&\quad \quad \quad \times
\det_{0\le i,j\le a/2-1}\left((3j+m/2+1)
\frac {(i+j+m/2)!}{(2i-j+m/2+1)!\,(2j-i)!}\right).
\end{align*}
Both determinants on the right-hand side of this identity can be
evaluated by means of Theorem~10 in \cite{KratBD}, which reads
\begin{multline} \label{eq:Th10} 
\det_{0\le i,j\le n-1}\left(\frac {(x+y+i+j-1)!}
{(x+2i-j)!\,(y+2j-i)!}\right)\\
=\prod _{i=0} ^{n-1}\frac {i!\,(x+y+i-1)!\,(2x+y+2i)_i\,(x+2y+2i)_i}
{(x+2i)!\,(y+2i)!}.
\end{multline}
This completes the proof of the theorem.
\end{proof}

\end{section}

\begin{section}{Proof of Theorem~\ref{thm:om3}} \label{corsec}

We prove Theorem~\ref{thm:om3} by counting the lozenge tilings of a
hexagon with side lengths $a,a+m,a,a+m,a,a+m$ and removed central triangle of
side length $m$ in two different ways.

\newbox\cycbox
\setbox\cycbox\hbox{\small
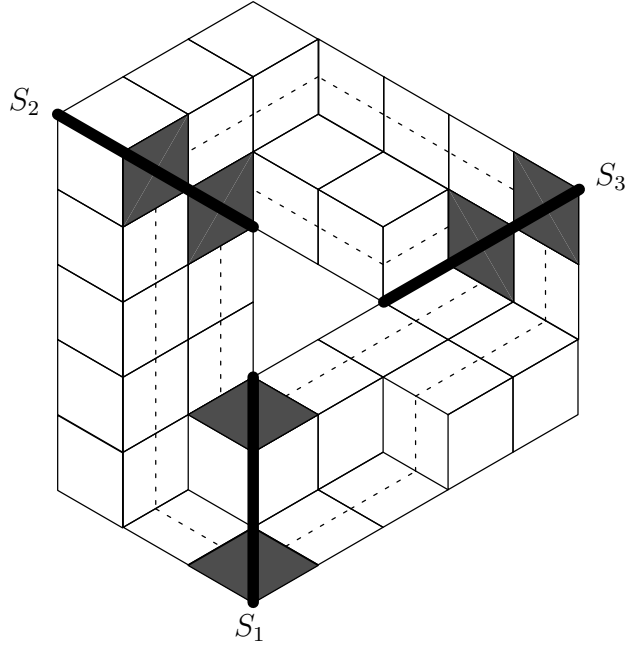
A lozenge tiling of a hexagon with $a=b=c=3$
and core of size $m=2$.}
\begin{figure} 
\centertexdraw{
\drawdim truecm \linewd.02
\RhombusA \RhombusB \RhombusA\RhombusA \RhombusB \RhombusA\RhombusA \RhombusB
\move(-.866025 -.5)
\RhombusA\RhombusB \RhombusB\RhombusB\RhombusB\RhombusA \RhombusB\RhombusA 
\move(-1.7305 -1)
\RhombusA\RhombusB\RhombusB\RhombusB \RhombusB\RhombusB\RhombusA \RhombusA 
\move(.866025 -4.5)
\RhombusA \RhombusB\RhombusA 
\move(1.7305 -4)
\RhombusA \RhombusA \RhombusB
\move(0.866025 -1.5) \RhombusC \RhombusC
\move(1.732 0) \RhombusC \RhombusC \RhombusC \RhombusC
\move(-1.73205 -1) \RhombusC
\move(-1.73205 -2) \RhombusC
\move(-1.73205 -3) \RhombusC
\move(-1.73205 -4) \RhombusC
\move(-1.73205 -5) \RhombusC
\move(3.464 -2) \RhombusC \RhombusC
\move(-4 -1.5)
%\ravec(1.5 .5)
\move(1.299 -4.75)
%\hodaSchritt 
\odaSchritt \odaSchritt \odaSchritt
%\hvdaSchritt \hhdaSchritt 
\rmove(-.433012 .75)
\odSchritt 
\hdaSchritt \hdaSchritt %\hodSchritt \hvdSchritt 
\rmove(-.433012 -.75)
\vdSchritt
\vdSchritt %\hhdSchritt
\rmove(.866025 -2)
%\hodaSchritt 
\odaSchritt \odaSchritt \vdaSchritt \odaSchritt 
\odaSchritt \vdaSchritt 
%\hvdaSchritt \hhdaSchritt 
\rmove(-.433012 .75)
\hdaSchritt
\hdaSchritt \hdaSchritt \odSchritt \odSchritt 
%\hodSchritt \hvdSchritt 
\rmove(-.433012 -.75)
\vdSchritt \vdSchritt \vdSchritt \vdSchritt
\hdSchritt %\hhdSchritt

\move(0 -1)  \rlvec(-.866025 -.5) \rlvec(0 -1)
\lfill f:.3
\rlvec(.866025 .5) \rlvec(0 1)
\lfill f:.3

\move(.866025 -1.5)  \rlvec(-.866025 -.5) \rlvec(0 -1)
\lfill f:.3
\rlvec(.866025 .5) \rlvec(0 1)
\lfill f:.3

\move(.866025 -4.5)  \rlvec(-.866025 -.5) \rlvec(.866025 -.5)
\lfill f:.3
\rlvec(.866025 .5) \rlvec(-.866025 .5)
\lfill f:.3

\move(.866025 -6.5)  \rlvec(-.866025 -.5) \rlvec(.866025 -.5)
\lfill f:.3
\rlvec(.866025 .5) \rlvec(-.866025 .5)
\lfill f:.3

\move(3.464 -2)  \rlvec(.866025 -.5) \rlvec(0 -1)
\lfill f:.3
\rlvec(-.866025 .5) \rlvec(0 1)
\lfill f:.3

\move(4.3301 -1.5)  \rlvec(.866025 -.5) \rlvec(0 -1)
\lfill f:.3
\rlvec(-.866025 .5) \rlvec(0 1)
\lfill f:.3

\linewd.15
\move(.866025 -4.5) \rlvec(0 -3)
\move(.866025 -2.5) \rlvec(-2.598 1.5)
\move(2.598 -3.5) \rlvec(2.598 1.5)

\htext(-2.4 -1){$S_2$}
\htext(.6 -8){$S_1$}
\htext(5.4 -2){$S_3$}
%\htext(-.8 0){$a$}
%\htext(4 -1){$b+m$}
%\htext(6.9 -3.9){$\left. \vbox{\vskip.5cm} \right\} c$}
%\htext(-3.1 -4){$c+m \left\{ \vbox{\vskip1.6cm} \right.$}
%\htext(0 -6){$b$}
%\htext(5 -6){$a+m$}
%\htext(1.7 -4){$\left. \vbox{\vskip1.1cm} \right\}m$}
%\rtext td:60 (4 -1.3) {$\left. \vbox{\vskip3.6cm} \right\} $}
%\rtext td:-60 (-.8 0){$\left\{ \vbox{\vskip1.6cm} \right. $}
%\rtext td:-60 (4.6 -5.3) {$\left. \vbox{\vskip2.5cm} \right\} $}
%\rtext td:60 (0.3 -5.6){$\left\{ \vbox{\vskip2.7cm}\right. $}
}
\caption{\unhbox\cycbox}
\label{fig:asy}
\end{figure}

First, we already know that this number equals \eqref{eq:enum} with
$a=b=c$. On the other hand, we claim that it equals $\det(I+B^3)$,
where, as before in the proof of Theorem~\ref{thm:Case10} in
Section~\ref{sec:outline}, $B=B(a,m)$ is
the $a\times a$ matrix with entries $\binom{m+i+j}{j}$, 
$0\le i,j \le a-1$, and
$I=I(a)$ is the $a\times a$ identity matrix.

To prove this claim,
we first note that $\det(I+B^3)$ is the sum of all principal minors of
$B^3$. Next we consider the construction used in
Section~\ref{sec:poly} in order to prove polynomiality in $m$ of the number
of lozenge tilings of a cored hexagon. 
I.e., we extend all sides of the removed triangle to the left (if
viewed from the interior of the triangle), as is indicated by the
thick segments, labeled as $S_1$, $S_2$, and $S_3$, in
Figure~\ref{fig:asy}. These segments cut the cored hexagon into three
regions. In particular, they cut some of the lozenges in two. (In
Figure~\ref{fig:asy}, these lozenges are shaded.) Subsequently, in
each of the three regions, we connect the ``cut" lozenges by paths,
by ``following" along the lozenges of the tiling, as is illustrated in
Figure~\ref{fig:asy} by the dashed lines. (Note the
difference between Figures~\ref{fig:asy} and \ref{cycpaths}.
In our special case $a=b=c$ {\em all\/} the paths form cycles.)

Let us number the possible
positions of the ``cut" lozenges, from inside to outside, by
$0,1,\dots,a-1$. Thus, the positions of the ``cut" lozenges on the
segment $S_1$ are 0 and 2, they are 0 and 1 on $S_2$, and they are 1
and 2 on $S_3$. 
The number of paths in the lower left region which start at position $i$
on $S_1$ and end at position $j$ on $S_2$ is $\binom{m+i+j}
{j}$, which is the $(i,j)$-entry of $B$. 
The rotational symmetry of the cored hexagon guarantees that
an analogous fact is true for the other regions.
Thus, the number of paths starting at position $i$ on $S_1$, then
running around the removed triangle, and finally 
ending at position $j$ on $S_1$,
equals the $(i,j)$-entry of $B^3$. If we have a {\em family} of
paths starting and ending at positions $i_1,i_2,\dots, i_k$, the
Lindstr\"om--Gessel--Viennot theorem (see Lemma~\ref{gv}) implies that
the number of these paths is the minor consisting of rows and
columns with indices $i_1,i_2,\dots,i_k$ of the matrix $B^3$. 
Thus, the number of these families of paths is the sum of all 
principal minors of $B^3$, which we
have already found to be equal to $\det (I+B^3)$.

Now we use the factorization 
\begin{equation*} 
I+B^3=(I+B)(\om I + B)({\overline \om} I+B),
\end{equation*}
where $\om$ is a primitive third root of unity.
Thus we have 
\begin{equation} \label{eq:det3}
\det (I+B^3)=\det(I+B) \cdot |\det(\om I+B)|^2.
\end{equation}
The left-hand
side equals \eqref{eq:enum} with $a=b=c$ by the above considerations,
and the
determinant $\det(I+B)$ has been computed by Andrews
\cite[Theorem~8]{AndrAN}, restated here as Theorem~\ref{thm:Andrews}.

Thus, a combination of \eqref{eq:det3}, Theorem~\ref{thm:Andrews} and
\eqref{eq:enum} with $a=b=c$ will give
$\det(\om I+B)$, the determinant that we want to compute, 
up to a complex factor of modulus
1. We note that the determinant is a polynomial in $m$. It is a
routine computation to verify that the determinant is the expression
claimed in Theorem~\ref{thm:om3}, up to this multiplicative constant.

In order to compute the multiplicative constant, we compute the
leading coefficient of the determinant as a polynomial in $m/2$, and
compare the result with the leading coefficient of the right-hand
side of \eqref{eq:om3}. Unfortunately, the leading coefficient of the
determinant cannot be determined straightforwardly by extracting the
leading coefficient of each of the entries and computing the
corresponding determinant, for the result would be zero. Therefore we
have to perform some manipulations of the matrix first to
avoid cancellation of leading terms. We use the strategy 
from \cite{MiRRAD}, which we have already used in
the proof of Lemma~\ref{lem:Z_n}. Instead of the determinant
$Z_n(x,\mu)$, we consider here the slightly different determinant
\begin{equation} \label{eq:Z_nom}
\widetilde{Z_n}(x,\mu)=\det_{0\le i,j\le n-1}\left(\om\delta_{ij}+\sum_{t,k=0}
^{n-1}\binom {i+\mu}t\binom kt\binom
{j-k+\mu-1}{j-k}x^{k-t}\right),
\end{equation}
where $\om$ is a primitive third root of unity.

Now we proceed analogously to the proof of Lemma~\ref{lem:Z_n}, i.e.,
we multiply $\widetilde{Z_n}(x,\mu)$ on the left by 
$\det(F(1)^t)$ and on the right by $\det(JF(x))$, where the matrices
$F(x)$ and $J$ are given in \eqref{eq:JF}, and use Chu--Vandermonde
summation several times. This yields
$$\widetilde{Z_n}(x,\mu)=\det_{0\le i,j \le n-1}(\om V(x,\mu)
+W(x,\mu)),$$
where $V(x,\mu)$ and $W(x,\mu)$ are the matrices defined in
equation~\eqref{eq:vw}.

Now let $x=1$, $\mu=m/2$, $n=a$, and $V=V(1,m/2)$, $W=W(1,m/2)$.
Again using Chu--Vandermonde summation,
we can express the desired determinant in terms of $V=(V_{ij})_{0\le
i,j\le a-1}$ and $W=(W_{ij})_{0\le i,j\le a-1}$:
\begin{equation} \label{eq:VW3}
\det(\om I+B)=\det(\om V+W),
\end{equation}
where
\begin{equation} \label{eq:V}
V_{2i+r,2j+s}=(-1)^{r+s}\binom{i+j+r+s+m/2} {s+2j-i}
\end{equation}
and 
\begin{equation} \label{eq:W}
W_{2i+r,2j+s}=\binom{i+j+m/2}{s+2j-i-r},
\end{equation}
where $r$ and $s$ are restricted to be $0$ or $1$.
Next we extract the leading coefficients of all the entries of
$\om V+W$, viewed as polynomials in $m/2$, and compute the
corresponding determinant. If we should obtain something nonzero, then this
must be the leading coefficient of the determinant $\det(\om V+W)$,
and hence of $\det(\om I+B)$, as a
polynomial in $m/2$.
Thus, we have to compute the determinant of the matrix
$L=(L_{ij})_{0\le i,j\le a-1}$, where
$$L_{2i+r,2j+s}=\begin{cases}
\frac {(-1)^{s+1}} {(s+2j-i)!}\om \quad &\text {if $r=1$,}\\
\frac {(-1)^{s}} {(s+2j-i)!}\om +\frac {1} {(s+2j-i)!} \quad &\text
{if $r=0$.}
\end{cases}$$
We add row 1 of $L$ to row 0, row 3 to row 2, etc. In that manner, we obtain
the matrix $L'=(L'_{ij})_{o\le i,j\le a-1}$, where
$$L'_{2i+r,2j+s}=\begin{cases}
\frac {1} {(s+2j-i)!} \quad  &\text {if  $r=0$, $2i\not= a-1$,}\\
\frac {(-1)^{s+1}} {(s+2j-i)!}\om \quad &\text {if } r=1,\\
\frac {((-1)^{s}\om+1)} {(s+2j-i)!} \quad &\text {if } 2i=a-1.
\end{cases}
$$
Clearly, we have $\det L=\det L'$, 
and we can take out $\om$ from all the rows
of $L'$ with odd row index. We get
$$\det L=\om^{\fl{\frac {a} {2}}}\det L'',$$
with the matrix $L''=(L''_{ij})_{0\le i,j\le a-1}$ defined by
$$L''_{2i+r,2j+s}=\begin{cases}
\frac {1} {(s+2j-i)!} \quad  &\text {if $r=0$, $2i\not= a-1$,}\\
\frac {(-1)^{s+1}} {(s+2j-i)!} \quad &\text {if } r=1,\\
\frac {((-1)^{s}\om+1)} {(s+2j-i)!} \quad &\text {if } 2i=a-1.
\end{cases}
$$
Now we add row 0 of $L''$ to row 1, row 2 to row 3, etc.
We obtain the matrix $L'''=(L'''_{ij})_{0\le i,j\le a-1}$, where
$$L'''_{2i+r,2j+s}=\begin{cases}
\frac {1} {(s+2j-i)!}\quad & \text {if $r=0$, $2i\not=a-1$,}\\
\frac {2} {(s+2j-i)!}\quad &\text {if $r=1$, $s=1$,}\\
0\quad  & \text {if $r=1$, $s=0$,}\\
\frac {((-1)^{s}\om+1)} {(s+2j-i)!} \quad &\text {if } 2i=a-1.
\end{cases}
$$
We rearrange the rows and columns simultaneously, so that 
the odd-numbered rows and columns 
come before the even-numbered, respectively.
Now we have obtained a block matrix with one
block formed by the rows and columns with odd indices and the other
one formed by the rows and columns with even indices.
Consequently, we have
\begin{multline*}
\det L =\om^{\fl{{a} /{2}}} 
 \det_{0\le i,j \le \fl{\frac {a-2} {2}}}{\(\frac {2} {(1+2j-i)!}\)}
\det_{0\le i,j \le \fl{\frac {a-1} {2}}}{
\begin{pmatrix} \frac
{1} {(2j-i)!} \quad &\text {\scriptsize $i\not=\frac {a-1} {2}$}\\
(\om+1)\frac {1} {(2j-i)!} \quad  &\text {\scriptsize $i=\frac {a-1} {2}$}
\end{pmatrix}
}
\\
=\om^{\fl{{a} /{2}}} (\om+1)^{\chi(a \text { odd})}
2^{\fl{ {a}/ {2}}}
\prod _{j=0} ^{\fl{ {(a-2)}/ {2}}}{\frac {1} {(2j+1)!}}
\prod _{j=0} ^{\fl{ {(a-1)}/ {2}}}{\frac {1} {(2j)!}}\\
\times\det_{0\le i,j \le \fl{ {(a-2)}/ {2}}}{\((2j-i+2)_{i}\)}
\det_{0\le i,j \le \fl{ {(a-1)}/ {2}}}{\( (2j-i+1)_i \)}  ,
\end{multline*} 
where we used the notation $\chi(\mathcal A)$=1 if $\mathcal A$ is
true and $\chi(\mathcal A)$=0 otherwise.
The two determinants can be evaluated by special cases of a variant of the
Vandermonde determinant evaluation which we state in Lemma~\ref{vanddet}
below.
After application of this lemma and some simplification we get
\begin{equation} \label{lc-even}
2^{l^2}\om^l\frac {\h(l)^2} {\h(2l)} 
\end{equation}
if $a$ is even, $a=2l$, and
\begin{equation} \label{lc-odd}
2^{l^2+l}\om^l(\om+1)\frac {\h(l)\h(l+1)} {\h(2l+1)} 
\end{equation}
if $a$ is odd, $a=2l+1$.

It is routine to check that the leading coefficient of the right-hand
side of \eqref{eq:om3}, viewed as a polynomial in $m/2$, is exactly
the same.

This finishes the proof of the theorem.\qed

\begin{lemma} \label{vanddet}
Let $p_i$ be a monic polynomial of degree $i$, $i=0,1,\dots,n$. Then
$$\det_{0\le i,j \le n}(p_i(X_j))=\prod _{0\le i <j \le n}
(X_j-X_i).$$
\quad \quad \qed
\end{lemma}
\end{section}

\begin{section}{Proof of Theorem~\ref{thm:om6}} \label{omsec}
If $a$ is even, $a=2l$ say, the formula can be derived analogously to
Theorem~\ref{thm:om3}. (The derivation of the latter was the subject
of the preceding section.) Here, the starting point is to do the
{\em $(-1)$-enumeration} (as opposed to ``ordinary" enumeration) 
of all the lozenge tilings of a
hexagon with side lengths $a,a+m,a,a+m,a,a+m$ and removed central triangle of
side length $m$ in two different ways.

First, the $(-1)$-enumeration of these lozenge tilings is given by 
\eqref{eq:-enum} with $a=b=c$. On the other hand, the arguments given
at the beginning of the preceding section, suitably modified, show
that it also equals $\det (-I+B^3)$, where $B$ is again the matrix
from the preceding section.

Now we use the factorization
\begin{equation} \label{eq:det6}
\det (-I+B^3)= \det(-I+B)\cdot |\det(\om I+B)|^2,
\end{equation}
where $\om$
is a primitive sixth root of unity. (Note that this equation is the
analogue of \eqref{eq:det3} in the present context.)
By the above considerations, the left-hand
side equals \eqref{eq:-enum} with $a=b=c$, and the
determinant $\det(-I+B)$ is computed in Theorem~\ref{thm:Zare1}.
This determines $\det(\om I+B)$ up to a multiplicative constant 
of modulus 1.
It is then a
routine computation to check that the result agrees with the
expression at the right-hand side of \eqref{eq:om6}, 
up to a factor of modulus 1.

In order to determine the multiplicative constant, one proceeds as in
the preceding section. In fact, the determination of the leading
coefficient of the determinant as a polynomial in $m/2$ given there
can be used here verbatim, because we treated $\om$ like an
indeterminate in the respective computations. Thus, the
leading coefficient is the expression in \eqref{lc-even}, 
with $\om$ now a
primitive sixth root of unity. It is routine to check that 
for $a=2l$ the
right-hand side of \eqref{eq:om6} has the same leading coefficient as
polynomial in $m/2$.

\smallskip
Now let us suppose that $a$ is odd, $a=2l+1$ say.
Unfortunately, the above strategy of determining the value of 
$\det(\om I+B)$ through equation \eqref{eq:det6} fails miserably
here, because $\det (-I+B^3)$ as well as $\det (-I+B)$ are zero
in the case of odd $a$ (compare Theorems~\ref{-enum} and
\ref{thm:Zare1}). Therefore we have to find a different line of
attack. We approach the evaluation of $\det(\om I+B)$, for odd $a$, 
by first transforming the determinant in the way we have already done
in the proofs of Lemma~\ref{lem:Z_n} and of Theorem~\ref{thm:om3}, and by
then applying once again the ``identification of factors" method to
evaluate the obtained determinant.

In fact, the manipulations explained in the preceding section that proved
\eqref{eq:VW3} (which
are based on multiplying the relevant matrix to the left and right by
suitable matrices, as elaborated in the proof of Lemma~\ref{lem:Z_n}
in Section~\ref{sec:Zare1}) remain valid in the present context,
again, because there $\om$ is treated like an indeterminate. Therefore
we have
$$\det(\om I+B)=\det(\om V+W),$$
where the matrices $V=(V_{ij})_{0\le i,j\le 2l}$ and
$W=(W_{ij})_{0\le i,j\le 2l}$ are again the matrices defined by
\eqref{eq:V} and \eqref{eq:W}.

Our goal is now to evaluate the determinant of
the matrix $\om V+W$. We denote this matrix by $X(2l+1,m/2)$.
The determinant $\det X(2l+1,m/2)$
is a polynomial in $m$, so we can indeed use the
``identification of factors" method to compute this determinant.
Again, there are several steps to be performed.
In Steps~1--4 below we prove that the right-hand side of \eqref{eq:om6}
does indeed divide the determinant as a polynomial in $m$.
In Step~5 we determine the maximal degree of the determinant as a
polynomial in $m$. It
turns out to be $(a^2-1)/4$, which is exactly the degree of the
right-hand side of \eqref{eq:om6} (for odd $a$, of course). Therefore
the determinant must be equal to the right-hand side of
\eqref{eq:om6}, up to a multiplicative constant. This multiplicative
constant is finally found to be 1 in Step~6.

\medskip
{\em Step 1.  $\prod _{i=0} ^{\fl{l/2}-1}(\frac {m}
{2}+2l-i+1)_{l-2i-1}$ divides the determinant $\det X(2l+1,m/2)$.}
Proceeding in the spirit of Step~1(b) in the proof of
Lemma~\ref{lem:amgerade}, we prove this by finding, for each linear
factor of the product, a linear combination
of the columns of $X(2l+1,m/2)$ which vanishes if the factor vanishes.
To be precise,
we claim that for $m/2=-3l+k+3d,$ $d\ge 0$ and $1\le k\le l-2d-1$
the following equation holds:
\begin{multline} \label{x1}
\sum _{j=0} ^{k-1}\binom{k-1}{j}\Big(\om \cdot(\text {column
$(2l-2d-2j-1)$ of $X(2l+1,-3l+k+3d)$})\\
+(\text {column $(2l-2d-2j-2)$ of
$X(2l+1,-3l+k+3d)$})\Big)=0.
\end{multline}
If we restrict the left-hand side of this equation to the $(2i)$th row, 
and simplify a little bit, it becomes
\begin{equation} \label{x1-2i}
\sum _{j=0} ^{k-1}\binom{k-1}{j}\left(\binom {i-2l+2d-j+k}
{2l-2d-2j-i-1} + \binom {i-2l+2d-j+k-1} {2l-2d-2j-i-2}
\right).
\end{equation}
It becomes $(\om-1)$ times the same expression if we restrict to the
$(2i+1)$th row. 

As is seen by inspection, 
the expression \eqref{x1-2i} vanishes trivially for $k=1$. 
{}From now on, let $k>1$. In order to establish that 
\eqref{x1-2i} vanishes in that case as well,
we first rewrite the sum \eqref{x1-2i} in hypergeometric
notation \eqref{eq:hyper}:
\begin{multline} \label{eq:x1-hyp}
{\frac {\left( k-1  \right)   
     ({ \textstyle 2 + 4 d + 2 i + k - 4 l}) _{-2 - 2 d - i + 2 l} }
     {(2l - 2 d - i -1)! }}\\
\times
{} _{4} F _{3} \!\left [ \matrix { 1 - k, {\frac 4 3}
      - {\frac k 3}, 1 + d + {\frac i 2} - l, {\frac 1 2} + d +
{\frac i 2} - l}\\
      { {\frac 1 3} - {\frac k 3}, 1 - 2 d - i - k + 2 l, 2 + 4 d + 2 i + k
      - 4 l}\endmatrix ; {\displaystyle 4}\right ].
\end{multline}
The hypergeometric 
summation formula which is relevant here, and as well in the
subsequent steps, is the following ``strange" evaluation of a
$_7F_6$-series, due to
Gessel and Stanton \cite[(1.7)]{GeStAA} (see also
\cite[(3.8.14), $c=1$, $a\to q^A$, etc., $q\to 1$]{GaRa}):
\begin{multline*} 
{} _{7} F _{6} \!\left [ \matrix { A, 1 + {\frac A 3}, B, 1 - B,
{\frac F 2},
   {\frac 1 2} + A - {\frac F 2} + n, -n}\\ { {\frac A 3}, 1 + {\frac A 2} -
   {\frac B 2}, {\frac 1 2} + {\frac A 2} + {\frac B 2}, 1 + A - F, -A + F - 2 n,
   1 + A + 2 n}\endmatrix ; {\displaystyle 1}\right ]  \\
=   {\frac {({ \textstyle 1 + A}) _{2 n} \,
      ({ \textstyle 1+\frac {A} {2}-\frac {B} {2}-\frac {F} {2}}) _{n} \,
      ({ \textstyle \frac {1} {2}+\frac {A} {2}+\frac {B} {2}-\frac {F} {2}}) _{n} } 
    {      ({ \textstyle 1 + A - F}) _{2 n} \,
      ({ \textstyle 1+\frac {A} {2}-\frac {B} {2}}) _{n} \,
      ({ \textstyle \frac {1} {2}+\frac {A} {2}+\frac {B} {2}}) _{n}
}},
\end{multline*}
where $n$ is a nonnegative integer.
If in this formula we let $B$ tend to infinity, we obtain
\begin{equation} \label{eq:5F4}
{} _{5} F _{4} \!\left [ \matrix { A, 1 + {\frac A 3}, {\frac F 2},
{\frac 1 2} +
   A - {\frac F 2} + n, -n}\\ { {\frac A 3}, 1 + A - F, -A + F - 2 n, 1 + A +
   2 n}\endmatrix ; {\displaystyle 4}\right ]  = 
  {\frac {({ \textstyle 1 + A}) _{2 n} } {({ \textstyle 1 + A - F}) _{2 n}
}}.
\end{equation}
In particular, this formula allows us to deduce that the left-hand side
of \eqref{eq:5F4} must be zero whenever $A$ is a negative integer.
This is seen as follows: Multiply both sides of \eqref{eq:5F4} by
\begin{equation} \label{eq:term}
(-A + F - 2 n)_{-A} \,(1 + A + 2 n)_{-A}.
\end{equation}
Then, for a fixed
negative integer $A$, the left-hand side becomes polynomial in $n$. The
right-hand side is zero for all $n$ larger than $-A/2$ because of the
presence of the term $({ \textstyle 1 + A}) _{2 n}$. The term
\eqref{eq:term} is nonzero for these values of $n$, therefore the
left-hand side of \eqref{eq:5F4} must be zero for these $n$. Since
these are infinitely many $n$, the left-hand side of \eqref{eq:5F4} 
must be in fact zero for all $n$. (An alternative way to see that
the left-hand side of \eqref{eq:5F4} vanishes for all negative $A$ is
by setting $c=1$ in \cite[(5.13)]{GaRaAB} or 
\cite[(3.8.11)]{GaRa}, then replace $a$ by
$q^A$, etc., and finally let $q\to1$ and $B\to \infty$.)

If we use \eqref{eq:5F4} with $A=1-k$, $F=2d+i-2l+2$, $n=2d+i+k-2l$, 
together with the above remarks, then we get immediately
that the $_4F_3$-series in \eqref{eq:x1-hyp} vanishes for $k>1$. (It
should be noted that, for this choice of parameters, the $_5F_4$-series
in \eqref{eq:5F4} reduces to the $_4F_3$-series in \eqref{eq:x1-hyp}.) 
Thus, equation \eqref{x1} is established.

\smallskip
{\em Step 2. $\prod _{i=0} ^{\fl{l/2}}(m/2+2l-i)_{l-2i}$ divides
the determinant.}
We claim that for $m/2=-3l+k+3d-1$, $d\ge 0$ and 
$1\le k\le l-2d$ the following equation holds:
\begin{multline} \label{x2}
\sum _{j=0} ^{k-1}\binom{k-1}{j}\Big((\text {column $(2l-2d-2j)$ of
$X(2l+1,-3l+k+3d-1)$})\\
+(2\om-1) \cdot(\text {column $(2l-2d-2j-1)$ of
$X(2l+1,-3l+k+3d-1)$})\\
+ (\om-1) \cdot(\text {column
$(2l-2d-2j-2)$ of $X(2l+1,-3l+k+3d-1)$})\Big)=0.
\end{multline} 
Restricted to the $(2i)$th row, 
the left-hand side of this equation becomes, after a little
simplification,
\begin{equation} \label{x2-2i}
(1+\om)\sum _{j=0} ^{k-1}\binom{k-1}{j}\left(\binom {i-2l+2d-j+k-1}
{2l-2d-2j-i} + \binom {i-2l+2d-j+k-2} {2l-2d-2j-i-1}
\right).
\end{equation}
Clearly, this expression vanishes 
for $k=1$. If $k>1$, we write \eqref{x2-2i}
in hypergeometric notation, to obtain
\begin{multline} \label{eq:x2-hyp}
(1+\om){\frac {\left( k-1  \right) 
     ({ \textstyle 4 d + 2 i + k - 4 l}) _{-1 - 2 d - i + 2 l} } 
   {(2l-2 d - i )! }}\\
\times
{} _{4} F _{3} \!\left [ \matrix { 1 - k, {\frac 4 3}
      - {\frac k 3}, {\frac 1 2} + d + {\frac i 2} - l, d + {\frac i 2} - l}\\ {
      {\frac 1 3} - {\frac k 3}, 2 - 2 d - i - k + 2 l, 4 d + 2 i + k -
      4 l}\endmatrix ; {\displaystyle 4}\right ].
\end{multline}
This time we use \eqref{eq:5F4} with $A=1-k$, $F=2d+i-2l+1$,
$n=2d+i+k-2l-1$. Together with the remarks accompanying
\eqref{eq:5F4}, this implies
immediately that the $_4F_3$-series in \eqref{eq:x2-hyp} vanishes for
$k>1$. 

On the other hand, restricted to the 
$(2i+1)$th row, 
the left-hand side of \eqref{x2} becomes, after a little
simplification,
\begin{multline*} 
(\om-1)\sum _{j=0} ^{k-1}\binom{k-1}{j}\left(\binom {i-2l+2d-j+k}
{2l-2d-2j-i-1} + \binom {i-2l+2d-j+k-1} {2l-2d-2j-i-2}
\right)\\
-\om\sum _{j=0} ^{k-1}\binom{k-1}{j}\left(\binom {i-2l+2d-j+k-1}
{2l-2d-2j-i} + \binom {i-2l+2d-j+k-2} {2l-2d-2j-i-1}
\right).
\end{multline*}
That the first sum vanishes was already shown in Step~1 (compare
\eqref{x1-2i}), that the second sum vanishes was shown just
above (compare \eqref{x2-2i}). Thus, equation \eqref{x2} is
established.

A short argument shows that the linear combinations of Step~1 are
independent of the linear combinations of Step~2. Let us denote the
columns of $X(2l+1,m/2)$ by $C_0, C_1, \dots, C_{2l}$.
In Step~1 we have
linear combinations of vectors of the form
$\om C_{2k+1}+C_{2k}$, whereas in Step~2  we have always 
linear combinations of vectors
of the form $C_{2k+2}+(2\om-1)C_{2k+1}+(\om-1)C_{2k}$. 
If these linear combinations were dependent we could
use the identity
$$(C_{2k+2}+(2\om-1)C_{2k+1}+(\om-1)C_{2k})-(\om+1)(\om
C_{2k+1}+C_{2k})=C_{2k+2}-2C_{2k},$$
and get a linear combination of vectors of the form
$\om C_{2k+1}+C_{2k}$ equal to a nonzero {\em real}
linear combination of the $C_i$'s, which is clearly impossible.

\smallskip
{\em Step 3. $\prod _{i=0} ^{\fl{l/2}-1}(m/2+3i+5/2)_{l-2i-1}$
divides the determinant.}
We claim that for $m/2=-k-\frac{3}{2}d$, $d$ odd, $d\ge1$, and $1\le k\le
l-d$ the following equation holds:
\begin{multline} \label{x3}
\sum _{i=0} ^{k-1}\binom{k-1}{i}\Big((\text {row $(2i+2d)$ of $X(2l+1,-k-\tfrac
{3} {2}d)$})\\
+\om\cdot (\text {row $(2i+2d+1)$ of $X(2l+1,-k-\tfrac {3}
{2}d)$})\Big)=0.
\end{multline}
Restricted to the $(2j)$th column, 
the left-hand side of this equation becomes, after a little
simplification,
\begin{equation} \label{x3-2j}
\sum _{i=0} ^{k-1}\binom{k-1}{i}\left(\binom {i-d/2+j-k+1}
{2j-i-d} + \binom {i-d/2+j-k} {2j-i-d} 
\right).
\end{equation}
It becomes $(\om-1)$ times the same expression if we restrict to the
$(2j+1)$th column.

Again, the expression \eqref{x3-2j} vanishes trivially
for $k=1$. In order to establish that 
\eqref{x3-2j} vanishes for $k>1$ as well,
we reverse the order of summation, and then write the sum
in hypergeometric notation. Thus we obtain
\begin{equation} \label{eq:x3-hyp}
{{\left( -1 \right) }^k}{\frac {\left( 1 - k \right) \,
     ({ \textstyle d - 2 j}) _{k-1} } 
   {(2j-d )! \,
     ({ \textstyle -{\frac {d} 2} + j}) _{d - 2 j + k} }}
     {} _{4} F _{3} \!\left [ \matrix { 1 - k, {\frac 4 3} - {\frac k 3}, 
    {\frac 1 2} - {\frac d 4} + {\frac j
      2} - {\frac k 2}}, 1 - {\frac d 4} + {\frac j 2} - {\frac k 2}\\ 
    { {\frac 1 3} - {\frac k 3},  1 + {\frac d 2} - j, 2 - d
      + 2 j - k}\endmatrix ; {\displaystyle 4}\right ].
\end{equation}
By \eqref{eq:5F4} with $A=1-k$, $F=1-d/2+j-k$,
$n=j-d/2$, together with the remarks accompanying \eqref{eq:5F4}, 
this implies
immediately that the $_4F_3$-series in \eqref{eq:x3-hyp} vanishes for
$k>1$. Thus, equation \eqref{x3} is established.

\smallskip
{\em Step 4. $\prod _{i=0} ^{\fl{l/2}}(m/2+3i+3/2)_{l-2i}$
divides the determinant.}
We claim that for $m/2=-k-\frac{3}{2}d-\frac {1} {2}$,
$d$ even, $d\ge 0$, and $1\le k\le l-d$ the following equation
holds:
\begin{multline} \label{x4}
\sum _{i=0} ^{k-1}\binom{k-1}{i}\Big((\text {row
$(2i+2d)$ of $X(2l+1,-k-\tfrac {3} {2}d-\tfrac {1} {2})$})\\
+(2-\om)\cdot
(\text {row $(2i+2d+1)$ of $X(2l+1,-k-\tfrac {3} {2}d-\tfrac {1} {2})$})\\
-\om\cdot (\text {row $(2i+2d+2)$ of $X(2l+1,-k-\tfrac {3} {2}d
-\tfrac {1} {2})$})\Big)=0.
\end{multline}
Restricted to the $(2j)$th column, 
the left-hand side of this equation becomes, after a little
simplification,
\begin{equation} \label{x4-2j}
(1-2\om)\sum _{i=0} ^{k-1}\binom{k-1}{i}\left(\binom {i-d/2+j-k+1/2}
{2j-i-d-1} + \binom {i-d/2+j-k-1/2} {2j-i-d-1} 
\right).
\end{equation}
Again, this expression vanishes trivially
for $k=1$. If $k>1$,
after reversion of summation, the hypergeometric form of \eqref{x4-2j} is
\begin{multline} \label{eq:x4-2j-hyp}
{{\left( -1 \right) }^k}{\frac {\left( 1 - k \right)  
       ({ \textstyle 1 + d - 2 j}) _{-1 + k} } 
   {(2j - d-1 )! \,
     ({ \textstyle -{\frac 1 2} - {\frac d 2} + j}) _{1 + d - 2 j + k} }}
\\
\times
     {} _{4} F _{3} \!\left [ \matrix { 1 - k, {\frac 4 3} - {\frac k 3},
       {\frac 1 4} - {\frac d 4} + {\frac j 2} - {\frac k 2}, 
     {\frac 3 4} - {\frac d 4} + {\frac j 2} - {\frac k 2}}\\ 
    { {\frac 1 3} - {\frac k 3}, {\frac 3 2} +
      {\frac d 2} - j, 1 - d + 2 j - k}\endmatrix ; {\displaystyle 4}\right
].
\end{multline}
Now we use \eqref{eq:5F4} with $A=1-k$, $F=1/2-d/2+j-k$,
$n=j-d/2-1/2$. Together with the remarks accompanying \eqref{eq:5F4}, 
this implies 
immediately that the $_4F_3$-series in \eqref{eq:x4-2j-hyp} vanishes for
$k>1$. 

On the other hand, restricted to the $(2j+1)$th column, 
the left-hand side of \eqref{x4} becomes, after a little
simplification,
\begin{multline} \label{x4-2j+1}
(\om-1)\sum _{i=0} ^{k-1}\binom{k-1}{i}\left(\binom {i-d/2+j-k+1/2}
{2j-i-d-1} + \binom {i-d/2+j-k-1/2} {2j-i-d-1} 
\right)\\
+\sum _{i=0} ^{k-1}\binom{k-1}{i}\left(\binom {i-d/2+j-k+3/2}
{2j-i-d+1} + \binom {i-d/2+j-k+1/2} {2j-i-d+1} 
\right).
\end{multline}
It was already shown just before that the first sum in
\eqref{x4-2j+1} vanishes (compare \eqref{x4-2j}). 
The second sum certainly vanishes for $k=1$. To see that it
vanishes for $k>1$ as well, we reverse the order of summation and
then convert the sum into hypergeometric notation,
\begin{multline} \label{eq:x4-2j+1-hyp}
{{\left( -1 \right) }^k}{\frac {\left( 1 - k \right) \,
       ({ \textstyle -1 + d - 2 j}) _{-1 + k} } 
   {(2j - d +1)! \,
     ({ \textstyle {\frac 1 2} - {\frac d 2} + j}) _{-1 + d - 2 j + k}
}}\\
\times
     {} _{4} F _{3} \!\left [ \matrix { 1 - k, {\frac 4 3} - {\frac k 3},
       {\frac 3 4} - {\frac d
      4} + {\frac j 2} - {\frac k 2}},{\frac 5 4} - {\frac d 4} +
{\frac j 2} - 
     {\frac k 2}\\ { {\frac 1 3} - {\frac k 3}, {\frac 1 2} +
      {\frac d 2} - j, 3 - d + 2 j - k}\endmatrix ; {\displaystyle 4}\right
].
\end{multline}
Again, by \eqref{eq:5F4}, this time with $A=1-k$, $F=3/2-d/2+j-k$,
$n=j-d/2+1/2$, together with the remarks accompanying \eqref{eq:5F4}, 
it follows 
immediately that the $_4F_3$-series in \eqref{eq:x4-2j+1-hyp} vanishes for
$k>1$. Thus, equation \eqref{x4} is established.

The linear combinations of Steps~3 and 4 are independent by the
argument used at the end of Step~2.

{\em Step~5. Determination of the degree of $\det X(2l+1,m/2)$ as a polynomial
in $m$.}
The $(i,j)$-entry of $X(2l+1,m/2)$, viewed as polynomial in $m$, 
has the degree $j-\fl{i/2}$.
Therefore, the determinant of $X(2l+1,m/2)$ has degree at most
$$\sum _{j=0} ^{2l}{j}-\sum _{i=0} ^{2l}{\fl{\frac {i}
{2}}}=l(l+1)=\frac{a^2-1}{4}$$
as a polynomial in $m$. 

{\em Step 6. Computation of the multiplicative constant}.
It suffices to compute the leading coefficient of the determinant
$\det X(2l+1,m/2)$ as a polynomial in $m/2$.
This leading coefficient can be computed as the determinant of the
leading coefficients of the individual entries. In fact, we already
did such a computation at the end of the proof of Theorem~\ref{thm:om3}
in the preceding section, with $\om$ a primitive third root of unity
instead of a primitive sixth root of unity. However,
since $\om$ was treated there as an indeterminate, everything can be used 
here as well. Thus we obtain the expression \eqref{lc-odd}, with $\om$
a primitive sixth root of unity. It is then routine to check that 
for $a=2l+1$ the
right-hand side of \eqref{eq:om6} has the same leading coefficient as
a polynomial in $m/2$.
\qed

\end{section}

\begin{section}{Comments and open problems}
\label{sec:Comments}

1) {\em Conjectured further enumeration results}. There is
overwhelming evidence (through computer supported empirical
calculations) that there are also ``nice" formulas for the
number of lozenge tilings of a cored hexagon for at
least two further locations of the core. 

\begin{figure} 
\centertexdraw{
\drawdim truecm \setunitscale.7
\linewd.02
\rhombus \rhombus \rhombus \rhombus \rhombus \rhombus \rhombus
\ldreieck
\move(-.866025 -.5)
 \rhombus \rhombus \rhombus \rhombus \rhombus \rhombus \rhombus \rhombus 
\move(-1.7305 -1)
\rhombus \rhombus \rhombus \rhombus \rhombus
\rhombus \rhombus \rhombus 
\move(-1.7305 -1)
\rdreieck \rhombus \rhombus \ldreieck  \rmove(.866025 -.5)
\rhombus \rhombus \rhombus 
\move(-1.7305 -2)
\rdreieck \rhombus \rhombus \ldreieck \rhombus \rhombus \rhombus 
\move(-1.7305 -3)
\rdreieck \rhombus \rhombus \rhombus \rhombus \rhombus 
\htext(-.8 .25){$a$}
\htext(4 -.7){$b+m$}
\htext(6.9 -4){$\left. \vbox{\vskip.5cm} \right\} c$}
\htext(-3.9 -3.95){$c+m \left\{ \vbox{\vskip1.2cm} \right.$}
\htext(0 -6){$b$}
\htext(4.9 -6){$a+m$}
 \rmove(.866025 -.5)
\htext(.833975 -4.45){$\left. \vbox{\vskip.8cm} \right\}m$}
\rtext td:60 (4 -1.3) {$\left. \vbox{\vskip2.6cm} \right\} $}
\rtext td:-60 (-.8 0.15){$\left\{ \vbox{\vskip1.2cm} \right. $}
\rtext td:-60 (4.6 -5.2) {$\left. \vbox{\vskip1.9cm} \right\} $}
\rtext td:60 (0.3 -5.7){$\left\{ \vbox{\vskip1.9cm}\right. $}

\move(10 0)
\bsegment
\rhombus \rhombus \rhombus \rhombus \rhombus \rhombus \rhombus
\ldreieck
\move(-.866025 -.5)
 \rhombus \rhombus \rhombus \rhombus \rhombus \rhombus \rhombus \rhombus 
\move(-.866025 -.5)
\rdreieck
\rhombus \rhombus \rhombus \rhombus \rhombus \rhombus \rhombus 
\move(-.866025 -1.5)
\rdreieck \rhombus \ldreieck \rmove(.866025 -.5)
\rhombus \rhombus \rhombus 
\move(-.866025 -2.5)
\rdreieck  \rhombus \ldreieck \rhombus \rhombus \rhombus 
\move(-1.7305 -3)
\htext(-.4 .45){$a$}
\htext(4 -.6){$b+m$}
\htext(6.9 -4){$\left. \vbox{\vskip.5cm} \right\} c$}
\htext(-3 -3.5){$c+m \left\{ \vbox{\vskip1.2cm} \right.$}
\htext(.8 -5.6){$b$}
\htext(5.2 -5.8){$a+m$}
\htext(.833975 -4.4){$\left. \vbox{\vskip.8cm} \right\}m$}
\rtext td:60 (4 -1.3) {$\left. \vbox{\vskip2.5cm} \right\} $}
\rtext td:-60 (-.4 0.3){$\left\{ \vbox{\vskip.8cm} \right. $}
\rtext td:-60 (5 -4.9) {$\left. \vbox{\vskip1.5cm} \right\} $}
\rtext td:60 (1.1 -5.2){$\left\{ \vbox{\vskip1.9cm}\right. $}
\esegment
\htext(-3 -8){\small a. Removal of the triangle which is}
\htext(8 -8){\small b. Removal of the triangle which is}
\htext(-3 -8.6){\small \hphantom{a. }off-center by one ``unit"}
\htext(8 -8.6){\small \hphantom{b. }off-center by 3/2 ``units"}
}

\caption{}
\label{fig:shifted}
\end{figure}

First, let $a$, $b$ and $c$ have the same parity, and consider a
hexagon with side lengths $a,b+m,c,a+m,b,c+m$ from which an
equilateral triangle of side length $m$ is removed which is
off-center by ``one
unit". To be more precise, let again $s_a$ be the side 
of the triangle which is parallel to the borders of the
hexagon of lengths $a$ and $a+m$, and similarly for $s_b$ and $s_c$.
Then the distance of $s_a$ to the border of length $a+m$ is the same
as the distance of the vertex of the triangle opposite to $s_a$ to
the border of length $a$.
The distance of $s_b$ to the border of length $b+m$ exceeds
the distance of the vertex of the triangle opposite to $s_b$ to
the border of length $b$ by two units.
Finally, the distance of $s_c$ to the border of length $c+m$ is 
two units
less than the distance of the vertex of the triangle opposite to $s_c$ to
the border of length $c$. See Figure~\ref{fig:shifted}.a for an example.
Then the following seems to be true.

\begin{conj} \label{conj1}
Let $a,b,c,m$ be nonnegative integers, $a,b,c$ having the same parity. 
The number of lozenge tilings of a hexagon with sides
$a,b+m,c,a+m,b,c+m$, with an equilateral triangle of
side length $m$ removed from the position that was described above 
(see Figure~\ref{fig:shifted}.a), equals
\begin{multline} \label{eq:1-step}
\frac {1} {4}\frac {\h(a + m)\h(b + m)\h(c + m)\h(a + b + c + m)
} 
{\h(a + b + m)\h(a + c + m)\h(b + c + m)
}\\
\times
\frac {\h(m + \left \lceil {\frac{a + b + c}{2}} \right \rceil)
\h(m + \left \lfloor {\frac{a + b + c}{2}} \right \rfloor)
} {\h({\frac{a + b}{2}} + m+1)     \h({\frac{a + c}{2}} + m-1)\h({\frac{b + c}{2}} + m)
}
\\
\times\frac {\h(\left \lceil {\frac{a}{2}} \right \rceil)
\h(\left \lceil {\frac{b}{2}} \right \rceil)
     \h(\left \lceil {\frac{c}{2}} \right \rceil)
     \h(\left \lfloor {\frac{a}{2}} \right \rfloor)\,
     \h(\left \lfloor {\frac{b}{2}} \right \rfloor)\,
     \h(\left \lfloor {\frac{c}{2}} \right \rfloor)\,
} 
{\h({\frac{m}{2}} + \left \lceil {\frac{a}{2}} \right \rceil)\,
     \h({\frac{m}{2}} + \left \lceil {\frac{b}{2}} \right \rceil)\,
     \h({\frac{m}{2}} + \left \lceil {\frac{c}{2}} \right \rceil)\,
\h({\frac{m}{2}} + \left \lfloor {\frac{a}{2}} \right \rfloor)\,
     \h({\frac{m}{2}} + \left \lfloor {\frac{b}{2}} \right \rfloor)\,
     \h({\frac{m}{2}} + \left \lfloor {\frac{c}{2}} \right \rfloor)\,
     }\\
\times
\frac {\h(\frac{m}{2})^2 \h({\frac{a + b + m}{2}})^2 
\h({\frac{a + c + m}{2}})^2 \h({\frac{b + c +
m}{2}})^2
} 
{\h({\frac{m}{2}} + \left \lceil {\frac{a + b + c}{2}} \right \rceil)
\h({\frac{m}{2}} + \left \lfloor {\frac{a + b + c}{2}} \right \rfloor)
\h({\frac{a + b}{2}}-1)\h({\frac{a + c}{2}}+1)\h({\frac{b + c}{2}})
               }P_1(a,b,c,m),
\end{multline}
where $P_1(a,b,c,m)$ is the polynomial given by
$$P_1(a,b,c,m)=\begin{cases} (a+b)(a+c)+2am&\text {if $a$ is
even,}\\
(a+b)(a+c)+2(a+b+c+m)m&\text {if $a$ is odd.}\end{cases}$$
\end{conj}

The reader should notice that the only differences between formulas
\eqref{eq:1-step} and \eqref{eq:enum} are in some hyperfactorials 
involving $(a+b)/2$ and $(a+c)/2$,
in the polynomial $P_1(a,b,c,m)$, which does not appear in \eqref{eq:enum}, 
and in the factor $1/4$ in front of \eqref{eq:1-step}.

The second case needs $a$ to have a parity different from $b$ and
$c$. Consider a
hexagon with side lengths $a,b+m,c,a+m,b,c+m$ from which an
equilateral triangle of side length $m$ is removed which is
off-center by ``3/2
units". To be more precise, with $s_a$, $s_b$, $s_c$
the sides of the triangle as above,
the distance of $s_a$ to the border of length $a+m$ is the same
as the distance of the vertex of the triangle opposite to $s_a$ to
the border of length $a$,
the distance of $s_b$ to the border of length $b+m$ exceeds
the distance of the vertex of the triangle opposite to $s_b$ to
the border of length $b$ by three units, and
the distance of $s_c$ to the border of length $c+m$ is 
three units
less than the distance of the vertex of the triangle opposite to $s_c$ to
the border of length $c$. See Figure~\ref{fig:shifted}.b for an
example.
Then the following seems to be true.

\begin{conj} \label{conj2}
Let $a,b,c,m$ be nonnegative integers, 
$a$ of parity different from the parity of $b$ and $c$.
The number of lozenge tilings of a hexagon with sides
$a,b+m,c,a+m,b,c+m$, with an equilateral triangle of
side length $m$ removed from the position that was described above 
(see Figure~\ref{fig:shifted}.b), equals
\begin{multline} \label{eq:3/2-step}
\frac {1} {16}\frac 
{\h\(a + m\)\h\(b + m\)\h\(c + m\)\h\(a + b + c + m\)} 
{\h\(a + b + m\)\h\(a + c + m\)\h\(b + c + m\)}\\
\times
\frac {\h\(\frac m 2\)^2\h\(\cl{\frac{a} {2}}\)\h\(\cl{\frac{b} {2}}\)
    \h\(\cl{\frac{c} {2}}\)\h\(\fl{\frac{a} {2}}\)
\h\(\fl{\frac{b} {2}}\)\h\(\fl{\frac{c} {2}}\)  } 
{\h\(\frac{m} {2} + \cl{\frac{a} {2}}\)\h\(\frac{m} {2} +
\cl{\frac{b} {2}}\)\h\(\frac{m} {2} + \cl{\frac{c} {2}}\)
\h\(\frac{m} {2} + \fl{\frac{a} {2}}\)\h\(\frac{m} {2} +
\fl{\frac{b} {2}}\)\h\(\frac{m} {2} + \fl{\frac{c} {2}}\)
}\\
\times
\frac 
{\h\(\cl{\frac {a + b} 2} + \frac{m} {2}\)
\h\(\fl{\frac{a + b}{2}} + \frac{m} {2}\)
\h\(\fl{\frac{a + c} {2}} + \frac{m} {2}\)\h\(\cl{\frac{a + c}
{2}} + \frac{m} {2}\)
    \h\(\frac{b + c} {2} + \frac{m} {2}\)^2
} 
{\h\(\frac{m} {2} + \cl{\frac{a + b + c} {2}}\)
\h\(\frac{m} {2} + \fl{\frac{a + b + c} {2}}\)
\h\(\fl{\frac{a + b} {2}}-1\)\h\(\cl{\frac{a + c}
{2}}+1\) \h\(\frac{b + c} {2}\)
}\\
\times
\frac {\h\(m + \cl{\frac{a + b + c} {2}}\)\h\(m + \fl{\frac{a + b + c} {2}}\)} 
{\h\(\fl{\frac{a + c} {2}}+m-1\)\h\(\frac{b + c} {2} +
m\)\h\(\cl{\frac{a+b}{2}}+m+1\)}
P_2(a,b,c,m),
\end{multline}
where the polynomial $P_2(a,b,c,m)$ is given by
$$P_2(a,b,c,m)=\begin{cases} \hbox to 9.5cm{$((a + b)^2-1) ((a + c)^2-1)+
4 a m (a^2 + 2 a b + b^2 + 2 a c + 3 b c + c^2$\hss}\\
\quad \quad  + 2 a m + 3 b m + 3 c m + 
    2 m^2-1)&\text {if $a$ is even,}\\
\hbox to 9.5cm{$((a + b)^2-1) ((a + c)^2-1)+
4 (a + b + c + m)m(a^2 + b c-1)$\hss}\\
&\text {if $a$ is odd.}
\end{cases}$$
\end{conj}

Again, the reader should notice that the only differences between formulas
\eqref{eq:3/2-step} and \eqref{eq:shifted} are in some hyperfactorials 
involving $(a+b)/2$ and $(a+c)/2$,
in the polynomial $P_2(a,b,c,m)$, which does not appear in
\eqref{eq:shifted}, 
and in the factor $1/16$ in front of \eqref{eq:3/2-step}.

Conjectured results about the $(-1)$-enumeration of the above two
families of lozenge tilings could be easily worked out as well, and
would have similar appearance, i.e., the result would be a
quotient of products of many ``nice" factors times an irreducible polynomial
of small degree. However, if one moves the triangle
farther away from the center, then, for both ordinary and
$(-1)$-enumeration, the irreducible polynomial factor
seems to grow rather quickly in degree, and is therefore difficult to predict
in general.

For a proof of Conjectures~\ref{conj1} and \ref{conj2}, one might go
through considerations analogous to those in
Section~\ref{sec:nonint}, i.e., convert the lozenge tilings into
families of nonintersecting lattice paths, and, by means of the
Lindstr\"om--Gessel--Viennot theorem (Lemma~\ref{gv}), obtain a determinant
for the number of lozenge tilings. This determinant, which then must
be evaluated, is
\begin{equation} \label{det-ep}
\det_{1\le i,j\le a+m} \begin{pmatrix} \dbinom{b+c+m}{b-i+j}&
\text {\scriptsize $1\le i\le a$}\\
\dbinom{\frac {b+c} {2}}{\frac {b+a} {2}-i+j+\varepsilon}&
\text {\scriptsize $a+1\le i \le a+m$}
\end{pmatrix},
\end{equation}
with $\varepsilon=1$ and $\varepsilon=3/2$, respectively. 
(The determinants in Lemmas~\ref{lem:nonint} and \ref{lem:nonint2} are
the respective special cases $\varepsilon=0$ and $\varepsilon=1/2$ of
\eqref{det-ep}.)

\smallskip
2) {\em A multidimensional analogue of Watson's $_3F_2$-summation,
and some variants}.
There is another possible way to approach the evaluation of the
determinants in Lemmas~\ref{lem:nonint} and \ref{lem:nonint2}. 
This approach consists of applying Laplace expansion to these
determinants. More precisely, we write an $(a+m)\times(a+m)$
determinant (such as the determinant in Lemma~\ref{lem:nonint} or
\ref{lem:nonint2}) as a (signed) sum of products of a minor formed of
elements of the first $a$ rows times the complementary minor formed
of elements of the last $m$ rows. That is, given an $(a+m)\times(a+m)$
matrix $M$, we write
\begin{equation} \label{eq:Laplace}
\det M=\sum _{K} ^{}(-1)^{s(K)}\(\det M^K\) \(\det M_{K'}\),
\end{equation}
where the sum is over all $a$-element subsets $K$ of
$\{1,2,\dots,a+m\}$, where $s(K)=\sum _{k\in K} ^{}k-\binom {a+1}2$,
$M^K$ denotes the submatrix of $M$ determined by the first $a$ rows and the 
columns with indices in $K$, $K'$ denotes the complement of $K$ in
$\{1,2,\dots,a+m\}$,
and $M_{K'}$ denotes the submatrix of $M$ determined by
the last $m$ rows and the columns with indices in $K'$.

The gain in applying \eqref{eq:Laplace} to our determinants in
Lemmas~\ref{lem:nonint} and \ref{lem:nonint2} is that the entries of 
the resulting minors
which then appear on the right-hand side of \eqref{eq:Laplace}
have now a uniform definition
(in contrast to the original
determinants), and can in fact easily be evaluated in closed form,
by means of the determinant evaluation 
\begin{equation} \label{eq:det-eval}
\det_{1\le i,j\le n}\(\pmatrix
A\\L_j-i\endpmatrix\)
=\frac {\prod _{1\le i<j\le n}
^{}(L_j-L_i)} {\prod _{i=1} ^{n}(A-L_i+n)!}
\frac {\prod _{i=1} ^{n}(A+i-1)!} {\prod _{i=1} ^{n}(L_i-1)!}.
\end{equation}
(This determinant evaluation is easily proved, e.g., by means of a general
determinant lemma from \cite[Lemma~2.2]{KratAM}; see
also \cite[Sec.~2.2 and (3.12)]{KratBN}). 
Thus, on the right-hand side of \eqref{eq:Laplace} we obtain a
multiple (hypergeometric) series for our determinants. If an evaluation of 
this multiple sum would appear in the existing literature, then we
would be immediately done. Unfortunately, this does not seem to be
the case. On the other hand, we did evaluate the determinants in
Sections~\ref{sec:det} and \ref{sec:det2}. Thus, comparison of the results
with the right-hand side in \eqref{eq:Laplace} 
establishes summation theorems for multiple hypergeometric series.
The summation theorem that results, after some replacement of
parameters, from the evaluations in Section~\ref{sec:det} 
of the determinant in Lemma~\ref{lem:nonint} is the following.

\begin{theorem} \label{thm:Watson1}
Let $a$ be a positive integer and $M$ be a nonnegative integer.
The multiple series
\begin{equation} \label{eq:Watson1-series}
\sum _{1\le k_1<k_2<\dots<k_a} ^{}\prod _{1\le i<j\le a}
^{}(k_i-k_j)^2\prod _{i=1} ^{a}\frac {(-M)_{k_i}\,(C)_{k_i}\,(B)_{k_i}} 
{k_i!\,\(\frac {a} {2}-\frac {M} {2}+\frac {C}
{2}\)_{k_i}\,(2B+a-1)_{k_i}}
\end{equation}
equals
\begin{multline} \label{eq:Watson1-GG}
{{\left( -1 \right) }^{{a/ 2}}}
{2^{a^2-a - a M}}
{\frac { {{M!}^a} 
      \prod_{i = 1}^{a}({ \textstyle B}) _{i-1}    
      } 
   {{{\left({ \textstyle {\frac a 2} + {\frac {{C}} 2} - {\frac M
2}}\right)^a
        _{ {M/ 2}-{{a}/ 2}} }}}}\\
\times
\prod_{i = 1}^{{a/ 2}}
        {\frac {{{\left( i-1 \right) !}^2}\,
            {{\left({ \textstyle {\frac 1 2} + {\frac {{C}} 2}}\right)^2 _{i-1} }}\,
            \left({ \textstyle B - {\frac {{C}} 2} + i-1}\right) 
  _{ {M/ 2} - {a/ 2} +
             1} \,\left({ \textstyle B - {\frac {{C}} 2} +
i}\right)
             _{ {M/ 2}-{{a}/ 2}} } 
          {\left( {\frac M 2} -i\right) !\,
            \left( {\frac M 2} - i + 1 \right) !\,
            {{\left({ \textstyle {\frac a 2} + B-{\frac 1 2}}\right)^2 
  _{ {M/ 2} - i + 1} }}
  \,{{\left({ \textstyle {\frac a 2} + B}\right)^2 _{i-1} }}\,
            \left({ \textstyle 1 + {\frac {{C}} 2} - i + {\frac M
2}}\right) _{2i-1} }} 
\end{multline}
if $a$ and $M$ are even, it equals
\begin{multline} \label{eq:Watson1-GU}
{{\left( -1 \right) }^{{a/ 2}}}
{2^{a^2-a - a M}}
{\frac { {{M!}^a} 
      \prod_{i = 1}^{a}({ \textstyle B}) _{i-1}    
 } 
   {{{\left({ \textstyle {\frac a 2} + {\frac {{C}} 2} - {\frac M
2}}\right)^a _{  {M/ 2}- {a/ 2} + {1/
        2}} }}}}
\prod _{i=1} ^{a/2}
\frac {{{\left( i-1 \right) !}^2}} 
{{{\left( {\frac M 2} - i + {\frac 1 2} \right) !}^2} }\\
\times
\prod_{i = 1}^{{a/ 2}}
        {\frac {
            {{\left({ \textstyle {\frac {{C}} 2}}\right) _{i-1} }}\,
            {{\left({ \textstyle {\frac {{C}} 2}}\right) _{i} }}\,
            {{\left({ \textstyle B - {\frac {{C}} 2} +
i-{\frac 1 2}}\right)^2
                _{ {M/ 2} - {a/ 2} + {1/ 2}} }}} 
          {
            \left({ \textstyle {\frac a 2} + B-{\frac 1 2}}\right) 
  _{ {M/ 2} - i + {1/ 2}
             } \,\left({ \textstyle  {\frac a 2} +
B-{\frac 1 2}}\right)
             _{ {M/ 2} - i + {3/ 2}} \,
            {{\left({ \textstyle {\frac a 2} + B}\right) ^2_{i-1} }}\,
            \left({ \textstyle 1 + {\frac {{C}} 2} - i + {\frac M
2}}\right) _{2i-1} }} 
\end{multline}
if $a$ is even and $M$ is odd, it equals
\begin{multline} \label{eq:Watson1-UG}
{{\left( -1 \right) }^{{M/ 2}}}
{2^{a^2-a - a M}}
{\frac {     {{M!}^a} \,
     \left({ \textstyle B - {\frac {{C}} 2}+{\frac a 2}}\right) 
  _{ {M/ 2} -
      {a/ 2} + {1/ 2}} 
      \prod_{i = 1}^{a}({ \textstyle B}) _{i-1}    
 } 
   {\left({\frac M 2}\right)!\,\left({ \textstyle {\frac a 2} + B}\right) 
   _{ {M/ 2}} \,
     {{\left({ \textstyle {\frac a 2} + {\frac {{C}} 2} - {\frac M
2}}\right) ^a_{ {M/ 2} - {a/ 2} + {1/
         2}} }}}}\\
\times
\prod_{i = 1}^{{(a-1)/ 2}}
        {\frac {{{\left( i-1 \right) !}}\,{{i!}}\,
            {{\left({ \textstyle {\frac {{C}} 2}}\right) ^2_{i} }}\,
            {{\left({ \textstyle B - {\frac {{C}} 2} +
i-{\frac 1 2}}\right)^2
                _{ {M/ 2} - {a/ 2} + {1/ 2}} }}} 
          {{{\left( {\frac M 2} -i\right) !}^2}\,
            {{\left({ \textstyle {\frac a 2} + B-{\frac 1 2}}\right) ^2_{i} }}\,
            {{\left({ \textstyle {\frac a 2} + B}\right)
    ^2_{ {M/ 2}-i} }}\,
            \left({ \textstyle {\frac 1 2} + {\frac {{C}} 2} - i +
{\frac M
2}}\right)
             _{2 i} }} 
\end{multline}
if $a$ is odd and $M$ is even, and it vanishes if both $a$ and $M$
are odd.
\end{theorem}

There are two interesting features of this summation theorem to be
observed. First, if we set $a=1$, the theorem reduces to a
terminating case of Watson's
$_3F_2$-summation (see \cite[(2.3.3.13); Appendix
(III.23)]{Sl}),
$$
{} _{3} F _{2} \!\left [ \matrix { A, C, B}\\ { {\frac {1 + A + C} 2}, 2
   B}\endmatrix ; {\displaystyle 1}\right ]  =
  \frac {\Gamma \( {\frac 1 2}\) \Gamma \( {\frac 1 2} + B\) \Gamma \(
{\frac 1 2} + {\frac A 2} +
   {\frac C 2}\) \Gamma \( {\frac 1 2} - {\frac A 2} - {\frac C 2} + B\)} 
  {\Gamma \( {\frac 1 2} + {\frac A
   2}\) \Gamma \( {\frac 1 2} + {\frac C 2}\) \Gamma \( {\frac 1 2} -
{\frac A 2} + B\) \Gamma \( {\frac 1 2} -
   {\frac C 2} + B\)},
$$
which is a summation formula which is not so often met. Second,
however, the above theorem is an {\em unusual\/} 
multidimensional analogue of Watson's $_3F_2$-summation, because of
the term $\prod _{1\le i<j\le a}
^{}(k_i-k_j)^2$ appearing in the summand. Whereas for series containing 
a term like $\prod _{1\le i<j\le a}
^{}(k_i-k_j)$ (i.e., the same term, but without the square) there
is now an extensive theory of summation and transformation
formulas (such a series is called a {\em hypergeometric series in $U(a)$}
or an {\em $A_a$ hypergeometric series}), mainly thanks to Milne and
Gustafson (see for example \cite{GustAA,MilnAK,MiLiAB,MiNeAA,SchlAB},
and the references contained therein), it is only occasionally that
series containing the square $\prod _{1\le i<j\le a}
^{}(k_i-k_j)^2$ appear. Most of the time, they arise from series
featuring Schur functions (see \cite[Theorem~6]{KratBM} for such an
example). However, our Theorem~\ref{thm:Watson1} does not seem to 
extend to a ``Schur function theorem."

The summation theorem that results from 
the evaluations in Section~\ref{sec:det2} 
of the determinant in Lemma~\ref{lem:nonint2} is a variant of the
preceding theorem. 

\begin{theorem} \label{thm:Watson2}
Let $a$ be a positive integer and $M$ be a nonnegative integer.
The multiple series
\begin{equation} \label{eq:Watson2-series}
\sum _{1\le k_1<k_2<\dots<k_a} ^{}\prod _{1\le i<j\le a}
^{}(k_i-k_j)^2\prod _{i=1} ^{a}\frac {(-M)_{k_i}\,(C)_{k_i}\,(B)_{k_i}} 
{k_i!\,\(\frac {a} {2}-\frac {M} {2}+\frac {C}
{2}+\frac {1} {2}\)_{k_i}\,(2B+a-2)_{k_i}}
\end{equation}
equals
\begin{multline}
{{\left( -1 \right) }^{{a/ 2}}}\,{2^{a^2-a - a M}}
{\frac { {{M!}^a} 
      \prod_{i = 1}^{a}({ \textstyle B}) _{i-1}    
 } 
   {{{\left({ \textstyle {\frac 1 2} + {\frac a 2} + {\frac C 2} -
{\frac M 2}}\right)
        ^a_{ {M/ 2}-{{a}/ 2}} }}}}
\\
\times
\prod_{i = 1}^{{a/ 2}}
\frac {{{\left( i-1 \right) !}^2\,
\left({ \textstyle {\frac C 2}}\right) _{i-1} \,
            \left({ \textstyle {\frac C 2}}\right) _{i} }} 
{\left( {\frac M 2} -i\right) !\,
            \left( {\frac M 2} -i+1\right) !\,
\left({ \textstyle {\frac a 2} + B-1}\right) _{i-1} \,
            \left({ \textstyle {\frac a 2} + B-1}\right) _{i} }
\\
\times
\prod_{i = 1}^{{a/ 2}}
        {\frac { \left({ \textstyle B - {\frac C 2} + i-{\frac 3 2}}\right)
_{ {M/ 2} - {a/ 2}
             + 1} \,\left({ \textstyle B - {\frac C 2} + i-{\frac 1 2}}\right)
             _{ {M/ 2}-{{a}/ 2}} } 
          { \left({ \textstyle {\frac a 2} + B-{\frac 1 2}}\right) _{
{M/ 2}-i} \,
            \left({ \textstyle {\frac a 2} + B-{\frac 1 2}}\right) _{
{M/ 2} - i + 1}
             \,\left({ \textstyle {\frac 1 2} + {\frac C 2} - i +
{\frac M 2}}\right)
             _{2 i} }} 
\end{multline}
if $a$ and $M$ are even, it equals
\begin{multline}
{{\left( -1 \right) }^{{a/ 2}}}\,{2^{a^2-a - a M}}
{\frac {{{M!}^a}\,
      \prod_{i = 1}^{a}({ \textstyle B}) _{i-1}  
   } 
   {{{\left({ \textstyle {\frac 1 2} + {\frac a 2} + {\frac C 2} -
{\frac M 2}}\right)
        ^a_{ {M/ 2} - {a/ 2} -{1/ 2} } }}}}
 \prod_{i = 1}^{{a/ 2}}
\frac {{{\left( i-1 \right) !}^2}} 
{{{\left( {\frac M 2} - i + {\frac 1 2} \right) !}^2}}
\\
\times
 \prod_{i = 1}^{{a/ 2}}
        {\frac {{{\left({ \textstyle {\frac 1 2} + {\frac C 2}}\right) ^2_{i-1} }}\,
            {{\left({ \textstyle B - {\frac C 2} + i-1}\right) ^2_{
{M/ 2} - {a/
                2} + {1/ 2}} }}} 
          {\left({ \textstyle {\frac a 2} + B-1}\right) _{i-1} \,
            \left({ \textstyle {\frac a 2} + B-1}\right) _{i} \,
            {{\left({ \textstyle {\frac a 2} + B-{\frac 1 2}}\right)
^2_{ {M/ 2} - i +
                {1/ 2}} }}\,
            \left({ \textstyle {\frac 1 2} + {\frac C 2} - i + {\frac M 2}}\right) _{2 i} }
          } 
\end{multline}
if $a$ is even and $M$ is odd, it equals
\begin{multline}
{{\left( -1 \right) }^{{M/ 2}}}\,{2^{a^2-a - a M}}
{\frac {{{M!}^a}\,
     \left({ \textstyle B - {\frac C 2}+{\frac a 2}-{\frac 1 2}}\right)
_{ {M/ 2} -
      {a/ 2} + {1/ 2}} \,
     } {\left( {\frac C 2} + {\frac M 2} \right) \,\left({\frac M
2}\right)!\,
     \left({ \textstyle {\frac a 2} + B-1}\right) _{ {M/
      2} - {a/ 2} + {1/ 2}} }}
\\
\times
\frac { \prod_{i = 1}^{a}({ \textstyle B}) _{i-1}  }
{{{\left({ \textstyle {\frac 1 2} + {\frac a 2} + {\frac C 2} -
{\frac M
         2}}\right) ^a_{ {M/ 2} - {a/ 2} -{1/ 2} } }}}
\prod_{i = 1}^{{(a-1)/ 2}}
\frac {\left( i-1 \right) !\,i!\,
\left({ \textstyle {\frac 1 2} + {\frac C 2}}\right) _{i-1} \,
            \left({ \textstyle {\frac 1 2} + {\frac C 2}}\right) _{i}} 
{{{\left( {\frac M 2} -i\right) !}^2\,
{{\left({ \textstyle {\frac a 2} + B-1}\right) ^2_{ {M/ 2} - i + 1} }}}}
\\
\times
\prod_{i = 1}^{{(a-1)/ 2}}
        {\frac {{{\left({ \textstyle B - {\frac C 2} + i-1}\right) ^2_{
{M/ 2} - {a/
                2} + {1/ 2}} }}} 
          {\left({ \textstyle {\frac a 2} + B-{\frac 1 2}}\right) _{i-1} \,
            \left({ \textstyle {\frac a 2} + B-{\frac 1 2}}\right) _{i} \,
            \left({ \textstyle {\frac C 2} - i + {\frac M 2}}\right) _{2i+1} }} 
\end{multline}
if $a$ is odd and $M$ is even, and it equals
\begin{multline}
{{\left( -1 \right) }^{{M/ 2}-{1/ 2}}}\,{2^{a^2-a - a M}}
{\frac { {{M!}^a}\,\left({ \textstyle B - {\frac C 2}-{\frac 1 2}}\right)
_{ {M/ 2} - {a/ 2} +
      1} 
     } {\left( {\frac C 2} + {\frac M 2} \right) \,
     \left( {\frac M 2} -{\frac {a} 2}\right) !\,
     \left({ \textstyle {\frac a 2} + B-1}\right) _{ {M/ 2} +
{1/ 2}} 
     }}
\\
\times
\frac { \prod_{i = 1}^{a}({ \textstyle B}) _{i-1} } 
{{{\left({ \textstyle {\frac 1 2} + {\frac a 2} + {\frac C 2} -
{\frac M 2}}\right)
         ^a_{ {M/ 2}-{{a}/ 2}} }}}
\prod_{i = 1}^{{(a-1)/ 2}}
\frac {\left( i-1 \right) !\,i!\,
{{\left({ \textstyle {\frac C 2}}\right) ^2_{i} }}} 
{{{\left( {\frac M 2} - i + {\frac 1 2} \right) !}^2\,
{{\left({ \textstyle {\frac a 2} + B-1}\right) ^2_{ {M/
                2} - i + {1/ 2}} }}}}
\\
\times
\prod_{i = 1}^{{(a-1)/ 2}}
        {\frac { \left({ \textstyle B - {\frac C 2} + i-{\frac 1 2}}\right)
_{ {M/
2}-{{a}/ 2}} \,\left({ \textstyle B - {\frac C 2} + i-{\frac 1 2}}\right)
             _{ {M/ 2} - {a/ 2} + 1} } 
          {\left({ \textstyle {\frac a 2} + B-{\frac 1 2}}\right)
_{i-1 } \,\left({ \textstyle {\frac a 2} + B-{\frac 1 2}}\right) _{i} \,
            \left({ \textstyle {\frac C 2} - i + {\frac M 2}}\right) _{2i+1} }} 
\end{multline}
if both $a$ and $M$ are odd.
\end{theorem}

In fact, the evaluations in Section~\ref{sec:det2} 
of the determinant in Lemma~\ref{lem:nonint2} establish even a
further variant of Theorem~\ref{thm:Watson1}. 
This variant is obtained as follows. Recall (see the
Introduction) that the
determinant in Lemma~\ref{lem:nonint2} arose from the case when
the parity of $a$ was different from that of $b$ and $c$, so that, in
order to have a well-defined enumeration problem,
we had to adjust the definition of a ``central" triangle of the
hexagon. 
What we did was to shift the {\em really central\/} triangle by half
a unit in parallel to the sides of the hexagon of length $a$ and
$a+m$. Now let us suppose that, 
unlike in that case, it is $b$ that has 
parity different from that of $a$ and $c$, so that the
``central" triangle in the sense of the Introduction is
the {\em really central\/} triangle shifted by half a unit in
parallel to the sides of the hexagon of length $b$ and $b+m$.
Clearly, our enumeration results in Theorems~\ref{shifted} and
\ref{-shifted} can be still used, we just have to
interchange the roles of $a$ and $b$. On the other hand, if we go
through the considerations in Section~\ref{sec:nonint} ({\em without\/}
interchange of the roles of $a$ and $b$, i.e., starting and end points
of the lattice paths are chosen on the sides of the hexagon of length
$a$ and $a+m$ and on the side of the triangle which is parallel),
then we obtain a certain determinant, which differs slightly from the
determinants in Lemmas~\ref{lem:nonint} and \ref{lem:nonint2}. Comparison of the
enumeration results with Laplace expansion \eqref{eq:Laplace} 
of the determinant establishes the following summation theorem.

\begin{theorem} \label{thm:Watson3}
Let $a$ be a positive integer and $M$ be a nonnegative integer.
The multiple series
\begin{equation} \label{eq:Watson3-series}
\sum _{1\le k_1<k_2<\dots<k_a} ^{}\prod _{1\le i<j\le a}
^{}(k_i-k_j)^2\prod _{i=1} ^{a}\frac {(-M)_{k_i}\,(C)_{k_i}\,(B)_{k_i}} 
{k_i!\,\(\frac {a} {2}-\frac {M} {2}+\frac {C}
{2}\)_{k_i}\,(2B+a-2)_{k_i}}
\end{equation}
equals
\begin{multline}
{{\left( -1 \right) }^{{a/ 2}}}\,{2^{a^2-a - a M}}
{\frac {{{M!}^a}\,
      \prod_{i = 1}^{a}({ \textstyle B}) _{i-1}  
 } 
   {\left({ \textstyle {\frac a 2} + B-1}\right) _{ {a/ 2}} \,
     {{\left({ \textstyle {\frac a 2} + {\frac C 2} - {\frac M 2}}\right)
^a_{ {M/
2}-{{a}/ 2}} }}}}
\prod_{i = 1}^{{a/2}}
\frac {{{\left( i-1 \right) !}^2}} 
{\left( {\frac M 2} -i\right) !\,
            \left( {\frac M 2}-i+1 \right) !}
\\
\times
\prod_{i = 1}^{{a/2}}
        {\frac {{{\left({ \textstyle {\frac 1 2} + {\frac C 2}}\right) ^2_{i-1} }}\,
            \left({ \textstyle B - {\frac C 2} + i-1}\right) _{
{M/
             2}-{{a}/ 2}} \,\left({ \textstyle B - {\frac C 2} +
i-1}\right) 
_{ {M/ 2} - {a/ 2} + 1} } 
          {{{\left({ \textstyle {\frac a 2} + B-1}\right) ^2_{i-1} }}\,
            \left({ \textstyle {\frac a 2} + B-{\frac 1 2}}\right) _{
{M/
2}-i} \,
            \left({ \textstyle {\frac a 2} + B-{\frac 1 2}}\right) _{
{M/ 2} - i + 1}
             \,\left({ \textstyle 1 + {\frac C 2} - i + {\frac M 2}}\right) _{2i-1} }}
\end{multline}
if $a$ and $M$ are even, it equals
\begin{multline}
{{\left( -1 \right) }^{{a/ 2}}}\,{2^{a^2-a - a M}}
{\frac {{{M!}^a}\,
      \prod_{i = 1}^{a}({ \textstyle B}) _{i-1}  
 } 
   {{{\left({ \textstyle {\frac a 2} + {\frac C 2} - {\frac M 2}}\right)
^a_{ {M/ 2} -
        {a/ 2} + {1/ 2}} }}}}
\prod_{i = 1}^{{a/2}}
\frac {{{\left( i-1 \right) !}^2}} 
{{{\left( {\frac M 2} - i + {\frac 1 2} \right) !}^2}}
\\
\times
\prod_{i = 1}^{{a/2}}
        {\frac {\left({ \textstyle {\frac C 2}}\right) _{i-1} \,
            \left({ \textstyle {\frac C 2}}\right) _{i} \,
            \left({ \textstyle B - {\frac C 2} + i-{\frac 3 2}}\right)
_{ {M/ 2} -
             {a/ 2} + {1/ 2}} \,
            \left({ \textstyle B - {\frac C 2} + i-{\frac 1 2}}\right)
_{ {M/ 2} -
             {a/ 2} + {1/ 2}} } 
          {\left({ \textstyle {\frac a 2} + B-1}\right) _{i-1} \,
            \left({ \textstyle {\frac a 2} + B-1}\right) _{i} \,
            {{\left({ \textstyle {\frac a 2} + B-{\frac 1 2}}\right)
^2_{ {M/ 2} - i +
                {1/ 2}} }}\,
            \left({ \textstyle 1 + {\frac C 2} - i + {\frac M 2}}\right) _{2i-1} }}
\end{multline}
if $a$ is even and $M$ is odd, it equals
\begin{multline}
{{\left( -1 \right) }^{{M/ 2}}}\,{2^{a^2-a - a M}}
{\frac {{{M!}^a}\,
     \left({ \textstyle B - {\frac C 2}-{\frac 1 2}}\right) _{
{M/ 2} - {a/ 2} +
      {1/ 2}}  \prod_{i = 1}^{a}({ \textstyle B}) _{i-1} 
          } 
   {\left({\frac M 2}\right)!\,\left({ \textstyle {\frac a 2} + B-1}\right) _{
{M/ 2} - {a/ 2}
      + {1/ 2}} \,{{\left({ \textstyle {\frac a 2} + {\frac C 2} -
{\frac M 2}}\right)
         ^a_{ {M/ 2} - {a/ 2} + {1/ 2}} }}}}
\\
\times
 \prod_{i = 1}^{{(a-1)/ 2}}
\frac {\left( i-1 \right) !\,i!\,
{{\left({ \textstyle {\frac C 2}}\right) ^2_{i} }}} 
{{{\left( {\frac M 2} -i\right) !}^2\,
{{\left({ \textstyle {\frac a 2} + B-1}\right) ^2_{ {M/ 2} - i + 1} }}}}
\\
\times
 \prod_{i = 1}^{{(a-1)/ 2}}
        {\frac { {{\left({ \textstyle B - {\frac C 2} + i-{\frac 1 2}}\right)
^2_{ {M/ 2} -
                {a/ 2} + {1/ 2}} }}} 
          { \left({ \textstyle {\frac a 2} + B-{\frac 1 2}}\right) _{i-1} \,
            \left({ \textstyle {\frac a 2} + B-{\frac 1 2}}\right) _{i} \,
            \left({ \textstyle {\frac 1 2} + {\frac C 2} - i + {\frac M 2}}\right) _{2 i} }
          } 
\end{multline}
if $a$ is odd and $M$ is even, and it equals
\begin{multline}
{{\left( -1 \right) }^{{M/ 2}+{1/ 2}}}\,{2^{a^2-a - a M}}
{\frac { {{M!}^a}\,\left({ \textstyle B - {\frac C 2}+{\frac a 2}-{\frac 1 2}}\right)
      _{ {M/ 2}-{{a}/ 2}} \,
      \prod_{i = 1}^{a}({ \textstyle B}) _{i-1}  
 } 
   {\left( {\frac M 2}-{\frac {a} 2} \right) !\,
     \left({ \textstyle {\frac a 2} + B-1}\right) _{ {M/ 2} +
{1/ 2}} \,
     {{\left({ \textstyle {\frac a 2} + {\frac C 2} - {\frac M 2}}\right)
^a_{ {M/
2}-{{a}/ 2}} }}}}
\\
\times
\prod_{i = 1}^{{(a-1)/ 2}}
\frac {\left( i-1 \right) !\,i!\,
\left({ \textstyle {\frac 1 2} + {\frac C 2}}\right) _{i-1} \,
            \left({ \textstyle {\frac 1 2} + {\frac C 2}}\right) _{i} } 
{{{\left( {\frac M 2} - i + {\frac 1 2} \right) !}^2\,
{{\left({ \textstyle {\frac a 2} + B-1}\right) ^2_{ {M/
                2} - i + {1/ 2}} }}}}
\\
\times
\prod_{i = 1}^{{(a-1)/ 2}}
        {\frac { \left({ \textstyle B - {\frac C 2} + i-1}\right) _{
{M/
             2}-{{a}/ 2}} \,\left({ \textstyle B - {\frac C 2} +
i-1}\right)
_{ {M/ 2} - {a/ 2} + 1} } 
          {\left({ \textstyle {\frac a 2} + B-{\frac 1 2}}\right) _{i-1
             } \,\left({ \textstyle {\frac a 2} + B-{\frac 1 2}}\right) _{i} \,
            \left({ \textstyle {\frac 1 2} + {\frac C 2} - i + {\frac M 2}}\right) _{2 i} }
          } 
\end{multline}
if both $a$ and $M$ are odd.
\end{theorem}

The reader should observe that, by similar considerations, i.e.,
by applying Laplace expansion \eqref{eq:Laplace}
to \eqref{det-ep}, Conjectures~\ref{conj1}
and \ref{conj2} are equivalent to further variations of
Theorem~\ref{thm:Watson1}. To be precise, Conjectures~\ref{conj1} and
\ref{conj2} could be proved by establishing summation theorems for
the multiple series 
\begin{equation} \label{eq:Watson-allg}
\sum _{1\le k_1<k_2<\dots<k_a} ^{}\prod _{1\le i<j\le a}
^{}(k_i-k_j)^2\prod _{i=1} ^{a}\frac {(-M)_{k_i}\,(C)_{k_i}\,(B)_{k_i}} 
{k_i!\,\(\frac {a} {2}-\frac {M} {2}+\frac {C}
{2}+\varepsilon\)_{k_i}\,(2B+a-1-2\varepsilon)_{k_i}},
\end{equation}
with $\varepsilon=1$ and $\varepsilon=3/2$, respectively.

\smallskip
3) {\em Are there $q$-analogues of our results?} By
``$q$-analogue", we mean, as usual, that objects $x$ are counted
with respect to a weight $q^{w(x)}$, where $w(x)$ is some statistic
defined on the objects. 
The question of whether there is a $q$-analogue, say of 
Theorems~\ref{enum} and \ref{shifted}, is
motivated by two facts: In the case of $m=0$ of Theorems~\ref{enum}
and \ref{shifted}, i.e., if one counts
lozenge tilings of a hexagon with no triangle removed, or,
equivalently, plane partitions contained in a given box, there is a
well-known $q$-analogue due to MacMahon
\cite[Sec.~429; proof in Sec.~494]{MM}, in which every plane
partition $P$ is given the weight $q^{|P|}$, where $|P|$ denotes the
number of ``boxes" (points, according to our definition of plane
partitions in Section~\ref{sec:cyclically}) of $P$. The result is
the $q$-analogue of formula \eqref{box} which is obtained by
replacing all factorials in \eqref{box} by the respective $q$-factorials.
Similarly, in the case $m=1$, $q$-analogues of Theorems~\ref{enum}
and \ref{shifted} can be gleaned from \cite[Theorem~3]{OkKr}, 
by setting $x_i=q^i$, $i=1,2,\dots,n+1$, respectively 
$x_i=q^i$, $i=1,2,\dots,n$, $x_{n+1}=0$, and using the hook-content
formula for the principal specialization of Schur functions (see
\cite[I, Sec.~3,  Ex.~1]{MacdAC}, \cite[Ex.~A.30, (ii)]{FuHaAA}).
The question of whether there are $q$-analogues for arbitrary $m$
remains open. 
Furthermore,
it would be particularly interesting if there were a
$q$-analogue of Theorem~\ref{thm:Andrews} 
that would specialize for $m=0$ to the
the statement of the Macdonald (ex)conjecture on cyclically symmetric
plane partitions (cf\@. \cite{MiRRAA}).

\end{section}


\begin{thebibliography}{10}

\bibitem{AndrAF}
G. E. Andrews, {\em The Theory of Partitions},
Encyclopedia of Math\@. and its Applications, Vol.~2,
Addison--Wesley, Reading, 1976.

\bibitem{AndrAN}
G. E. Andrews, {\em Plane partitions (III): The weak Macdonald
conjecture}, Invent\@. Math. {\bf 53} (1979), 193--225.

\bibitem{Ba}
W.~N.~Bailey, {\em Generalized hypergeometric series},
Cambridge University Press, Cambridge, 1935.

\bibitem{CiucAI}
M.    Ciucu, {\em Enumeration of lozenge tilings of punctured
hexagons}, J.~Combin\@. Theory Ser.~A {\bf 83} (1998), 268--272.

\bibitem{CiucAH}
M.    Ciucu, {\em Plane partitions I: A generalization of MacMahon's
formula}, preprint, {\tt math/9808017}.

\bibitem{CiKrAB} 
M.    Ciucu and C. Krattenthaler, {\em Plane partitions II: $5\frac {1}{2}$ 
symmetry classes}, 
in: Combinatorial Methods in Representation
Theory, M.~Kashiwara, K.~Koike, S.~Okada, I.~Terada, H.~Yamada, eds.,
Advanced Study in Pure Mathematics, RIMS, Kyoto (to appear);
{\tt math/9808018}.

\bibitem{David-Tomei}
G. David and C. Tomei, {\em The problem of the calissons}, 
Amer\@. Math\@. Monthly {\bf 96} (1989), 429--431.

\bibitem{FiscAA}
I. Fischer, {\em Enumeration of rhombus tilings of a hexagon which
contain a fixed rhombus in the centre}, preprint, {\tt math/9906102}.

\bibitem{FuHaAA} W.    Fulton and J. Harris,
{\it Representation Theory},
 Sprin\-ger--Ver\-lag,
New York, 1991.

\bibitem{GaRa}
G.~Gasper and M.~Rahman, {\em Basic hypergeometric series},
Encyclopedia of Mathematics And Its Applications 35, Cambridge
University Press, Cambridge, 1990.

\bibitem{GaRaAB} 
G.    Gasper and M. Rahman, {\em An indefinite bibasic summation formula 
and some quadratic, cubic and quartic summation and transformation
formulae}, Can\@. J. Math\@. {\bf 42} (1990), 1--27.

\bibitem{GeStAA} 
I. M. Gessel and D. Stanton, {\em Strange evaluations of hypergeometric 
series}, SIAM J. Math\@. Anal\@. {\bf 13} (1982), 295--308.

\bibitem{gv}
I. M.~Gessel and X.~Viennot, {\em Determinant, paths and plane
partitions}, preprint, 1989.

\bibitem{GustAA}
R. A. Gustafson,
{\em The Macdonald identities for affine root systems of classical type 
and hypergeometric series very-well-poised on semisimple Lie
algebras}, Ramanujan International Symposium on Analysis 
(Dec.~26th to 28th, 1987, Pune, India), N.~K.~Thakare, ed., 1989,
187--224.

\bibitem{HeGeAA} 
H.    Helfgott and I. M. Gessel,
{\em Exact enumeration of tilings of diamonds and hexagons with
defects}, Electron\@. J. Combin\@. {\bf 6} (1) (1999), \#R16, 26~pp.

\bibitem{KratAM} 
C.~Krattenthaler, {\em Generating functions for plane partitions of a given
shape\/},  Manuscripta Math.\ {\bf 69}, (1990), 173--202.

\bibitem{KratBD}
C.~Krattenthaler, {\em Determinant identities and a generalization of the
number of totally symmetric self-complementary plane partitions\/},
Electron.\ J.\ Combin.\ {\bf 4}(1) (1997), \#R27, 62~pp.

\bibitem{KratBI}
C.~Krattenthaler, {\em An alternative evaluation of the Andrews--Burge
determinant\/}, in: Mathematical Essays in Honor of Gian-Carlo Rota,
B.~E.~Sagan, R.~P.~Stanley, eds., Progress in Math., vol.~161, 
Birkh\"auser, Boston, 1998, pp.~263--270.

\bibitem{KratBN}
C.    Krattenthaler, {\em Advanced determinant calculus}, S\'eminaire 
Lotharingien Combin\@. {\bf 42} (``The Andrews Festschrift") (1999),
paper~B42q, 67~pp.

\bibitem{KratBM}
C.    Krattenthaler, {\em Schur function identities and the number of 
perfect matchings of holey Aztec rectangles}, Contemporary Math\@.
(to appear).

\bibitem{KupeAA}
G.    Kuperberg, {\em Symmetries of plane partitions and the permanent 
determinant method}, J.~Combin\@. Theory Ser.~A (1994) {\bf 68},
115--151.

\bibitem{KupeAG}
G.    Kuperberg, {\em An exploration of the permanent-determinant
method}, Electron\@. J. Combin\@. {\bf 5} (1998), \#R46, 34~pp.

\bibitem{LindAA}
B.    Lindstr\"om, {\em On the vector representations of induced
matroids}, Bull\@. London Math\@. Soc\@. {\bf 5} (1973), 85--90.

\bibitem{MacdAC} I. G. Macdonald,
{\em Symmetric Functions and Hall Polynomials},
 second edition,
 Oxford University Press, New York/Lon\-don, 1995.

\bibitem{MM}
P. A.~MacMahon, {\em Combinatory Analysis\/}, vol.~2, Cambridge University
Press, 1916; reprinted by Chelsea, New York, 1960.

\bibitem{MiRRAA}
W. H. Mills, D. H. Robbins and H. Rumsey, {\em Proof of the Macdonald 
conjecture}, Invent\@. Math\@. {\bf 66} (1982), 73--87.

\bibitem{MiRRAD}
W. H. Mills, D. H. Robbins and H. Rumsey, {\em Enumeration of a symmetry 
class of plane partitions}, Discrete Math\@. {\bf 67} (1987), 43--55.

\bibitem{MilnAK}
S. C. Milne, {\em $q$-analog of a Whipple's transformation for hypergeometric 
series in $U(n)$}, Adv\@. in Math\@. {\bf 108} (1994), 1--76.

\bibitem{MiLiAB}
S. C. Milne and G. M. Lilly, {\em Consequences of the $A_l$ and $C_l$ 
Bailey transform and Bailey lemma}, Discrete Math\@. {\bf 139} (1995),
319--346.

\bibitem{MiNeAA}
S. C. Milne and J. W. Newcomb, {\em $U(n)$ very-well-poised $_{10}\phi_9$ 
transformations}, J. Comput\@. Appl\@. Math\@. {\bf 68} (1996),
239--285.

\bibitem{OkKr}
S.~Okada and C.~Krattenthaler, {\em The number of rhombus tilings of a
``punctured" hexagon and the minor summation formula}, Adv\@. Appl\@.
Math\@. {\bf 21} (1998), 381--404.

\bibitem{Propp}
J.~Propp, {\em Twenty open problems on enumeration of matchings\/},
manuscript, 1996, {\tt math/9801060}.

\bibitem{Propp2}
J.~Propp, {\em Enumeration of matchings: Problems and progress},
in: "New Perspectives in Algebraic
 Combinatorics," L.~Billera, A.~Bj\"orner, C.~Greene, R.~Simion, and
R.~P.~Stanley, eds.,
 Mathematical Sciences Research Institute Publications, vol.~38, 
Cambridge University Press, Cambridge, 1999, pp.~255--291.

\bibitem{SchlAB}
M.    Schlosser, {\em Multidimensional matrix inversion and $A_r$ and 
$D_r$ basic hypergeometric series}, The Ramanujan J. {\bf 1} (1997),
243--276.

\bibitem{Sl}
L.~J.~Slater, {\em Generalized hypergeometric functions},
Cambridge University Press, Cambridge, 1966.

\bibitem{StanAI}
R. P. Stanley, {\em Symmetries of plane partitions},
J. Combin\@. Theory A {\bf 43} (1986), 103--113; Erratum: ibd\@. {\bf 44}
(1987), 310.

\bibitem{StanAP}
R. P. Stanley, {\em Enumerative Combinatorics},
Vol.~1, Wadsworth \& Brooks/Cole, Pacific Grove, California, 1986;
reprinted by Cambridge University Press, Cambridge, 1998.

\bibitem{StemAL}
J. R. Stembridge, {\em Some hidden relations involving the ten symmetry 
classes of plane partitions}, J.~Combin\@. Theory Ser\@. A {\bf 68}
(1994), 372--409.

\bibitem{StemAI}
J. R. Stembridge, {\em On minuscule representations, plane partitions and 
involutions in complex Lie groups}, Duke Math\@. J. {\bf 73} (1994),
469--490.

\bibitem{StemAZ}
J. R. Stembridge, {\em Strange enumerations of CSPP's and TSPP's},
unpublished manuscript, 1993, available at {\tt
http://www.math.lsa.umich.edu/\~{}jrs/other.html}.

\bibitem{SulaAC}
R. A. Sulanke, {\em A determinant for $q$-counting lattice paths}, 
Discrete Math\@. {\bf 81} (1990), 91--96.

\bibitem{Th}
J.~Thomae, {\em \"Uber die Funktionen, welche durch Reihen von der
Form dargestellt werden: \dots},
J.~f\"ur Math\@. (Crelle) {\bf 87} (1879), 26--73.

\end{thebibliography}
\end{document}